\documentclass[12pt]{article}
\usepackage{graphicx}
\usepackage{amssymb,amsmath,amsfonts,relsize}
\usepackage[usenames]{color}
\usepackage{cite}
\usepackage{esvect}
\usepackage{txfonts}

\makeatletter
\renewcommand{\fnum@figure}{Fig. \thefigure}
\makeatother
\newtheorem{theorem}{Theorem}[section]
\newtheorem{definition}[theorem]{Definition}
\newtheorem{corollary}[theorem]{Corollary}
\newtheorem{proposition}[theorem]{Proposition}
\newtheorem{remark}[theorem]{Remark}
\newtheorem{lemma}[theorem]{Lemma}

\newtheorem{problem}[theorem]{Problem}

\newtheorem{conjecture}[theorem]{Conjecture}
\newtheorem{assumption}[theorem]{Assumption}
\newtheorem{claim}[theorem]{Claim}
\newcommand {\Ac}      {{\mathcal A}}

\newcommand {\Ec}      {{\mathcal E}}

\newcommand {\Hc}      {{\mathcal H}}
\newcommand {\Ic}      {{\mathcal I}}

\newcommand {\Kc}      {{\mathcal K}}

\newcommand {\Mc}      {{\mathcal M}}

\newcommand {\Pc}      {{\mathcal P}}

\newcommand {\Rc}      {{\mathcal R}}

\newcommand {\Tc}      {{\mathcal T}}

\newcommand {\Wc}      {{\mathcal W}}
\newcommand {\R}       {{\bf R}}

\newcommand {\RN}      {\R^n}
\newcommand {\RT}      {\R^2}

\newcommand {\tA}      {\widetilde{A}}
\newcommand {\tB}      {\widetilde{B}}
\newcommand {\tC}      {\widetilde{C}}

\newcommand {\tF}      {\widetilde{F}}
\newcommand {\tG}      {\widetilde{G}}
\newcommand {\tH}      {\widetilde{H}}
\newcommand {\tI}      {\widetilde{I}}
\newcommand {\tJ}      {\widetilde{J}}
\newcommand {\tK}      {\widetilde{K}}

\newcommand {\tS}      {\widetilde{S}}
\newcommand {\tT}      {\widetilde{T}}

\newcommand {\tAc}     {\widetilde{{\mathcal A}}}

\newcommand {\tMc}     {\widetilde{{\mathcal M}}}

\newcommand {\tKc}     {\widetilde{\Kc}}
\newcommand {\ta}      {\tilde{a}}
\newcommand {\tb}      {\tilde{b}}

\newcommand {\tf}      {\tilde{f}}

\newcommand {\tilh}    {\tilde{h}}

\newcommand {\tr}      {\tilde{r}}
\newcommand {\ts}      {\tilde{s}}
\newcommand {\tilt}    {\tilde{t}}
\newcommand {\tu}      {\tilde{u}}
\newcommand {\tv}      {\tilde{v}}
\newcommand {\tw}      {\tilde{w}}

\newcommand {\ty}      {\tilde{y}}

\newcommand {\tgm}     {\tilde{\gamma}}
\newcommand {\tlm}     {\tilde{\lambda}}
\newcommand {\trh}     {\tilde{\rho}}

\newcommand {\tell}    {\tilde{\ell}}

\newcommand {\hM}      {\widehat{M}}
\newcommand {\mfM}      {{\mathfrak M}}
\newcommand {\mfC}      {{\mathfrak C}}

\newcommand {\mfR}      {{\mathfrak R}}
\newcommand {\ve}      {\varepsilon}

\newcommand {\MR}      {(\Mc,\rho)}
\newcommand {\MS}      {\mfM}

\newcommand {\BS}      {B}
\newcommand {\vl}      {\vec{\lambda}}
\newcommand {\vf}      {\varphi}
\newcommand {\emp}     {\emptyset}

\newcommand {\CX}      {\mfC(\X)}
\newcommand {\CRT}     {\Conv(\RT)}
\newcommand {\RCT}     {\mfR(\RT)}
\newcommand {\SO}      {{\bf S}_1}
\newcommand {\RL}      {\Rc_{F}}

\newcommand {\dq}      {\tilde{\dt}}

\newcommand {\HR}      {\Hc}
\newcommand {\dl}      {\delta}
\newcommand {\AF}      {\vf_F}
\newcommand {\AG}      {\vf_G}
\newcommand {\DV}      {n^{\bot}}
\newcommand {\slbig}   {\mathlarger{\mathlarger{/}}}

\newcommand {\cbig}    {\mathlarger{\mathlarger{\cap}}}
\newcommand {\cbg}    {\,\mathlarger{\mathlarger{\cap}}\,}

\newcommand {\cupbig}    {\mathlarger{\mathlarger{\cup}}}
\newcommand {\smed}    {\mathlarger{\sum}}

\newcommand {\blbig}    {\mathlarger{\mathlarger{\bullet\hspace{0.2mm}}}}

\newcommand {\BX}      {B_{\X}}
\newcommand {\BXR}     {I_0}

\newcommand {\LTI}     {\ell^2_\infty}
\newcommand {\LNI}     {\ell^n_\infty}
\newcommand {\LTT}     {\ell^2_2}
\newcommand {\LNT}     {\ell^n_2}
\newcommand {\X}       {X}
\newcommand {\TI}      {\tilde{I}}

\newcommand {\KM}      {\Kc_m(\X)}
\newcommand {\CNV}     {\Conv(\RT)}
\newcommand {\DXY}     {\Delta_n(x,y)}
\newcommand {\DXP}     {\Delta_n(x,x')}
\newcommand {\DYP}     {\Delta_n(y,y')}
\newcommand {\DGXY}    {\Delta_g(x,y)}
\newcommand {\DGXP}    {\Delta_g(x,x')}
\newcommand {\DGYP}    {\Delta_g(y,y')}
\newcommand {\HPL}     {\Hc\Pc(\RT)}
\newcommand {\SXP}     {\sin{\AF(x,x')}}
\newcommand {\SYP}     {\sin{\AF(y,y')}}
\newcommand {\SXPG}    {\sin{\AG(x,x')}}
\newcommand {\SYPG}    {\sin{\AG(y,y')}}
\newcommand {\TSXP}    {\sin{\vf_{\tF}(X,X')}}
\newcommand {\TSYP}    {\sin{\vf_{\tF}(Y,Y')}}

\newcommand {\BAL}[3]  {{\cal BR\,}[{#1}\!:\!{#2};{#3}]}
\newcommand {\ip}[1]   {\langle{#1}\rangle}

\newcommand {\al}      {\alpha}
\newcommand {\VST}     {\vspace*{1mm}}

\newcommand {\Lip}     {\operatorname{Lip}}

\newcommand {\ST}      {\operatorname{St}}

\newcommand {\dhf}     {\operatorname{d_H}}
\newcommand {\ds}      {\operatorname{d}}
\newcommand {\dt}      {\operatorname{d}}

\newcommand {\Prj}     {\operatorname{Pr}}
\newcommand {\diam}    {\operatorname{diam}}
\newcommand {\dist}    {\operatorname{dist}}
\newcommand {\length}  {\operatorname{length}}

\newcommand {\cntr}    {\operatorname{cntr}}

\newcommand {\Conv}    {\operatorname{Conv}}
\newcommand {\sign}    {\operatorname{sign}}
\newcommand {\smsk}    {\smallskip}
\newcommand {\msk}     {\medskip}
\newcommand {\bsk}     {\bigskip}
\newcommand {\bx}      {\hspace{10mm}$\blacksquare$}
\newcommand {\rbx}     {\hspace{10mm}$\vartriangleleft$}

\newcommand {\nn}      {\nonumber}
\newcommand {\rf}[1]    {(\ref{#1})}      
\newcommand {\reff}[1] {\ref{#1}}         
\newcommand{\lbl}[1]      {\label{#1}}       
\newcommand{\be}          {\begin{eqnarray}}
\newcommand{\bel}[1]      {\begin{eqnarray} \label{#1}}
\newcommand{\ee}           {\end{eqnarray}}
\newcommand {\SECT}[2] {\section*{\centerline{\normalsize
{\bf #1}}} \setcounter{section}{#2}
\setcounter{theorem}{0}\setcounter{equation}{0}}
\begin{document}
\parindent 1em
\parskip 0mm
\medskip
\centerline{{\bf The Core of a 2-Dimensional Set-Valued Mapping.}}
\msk
\centerline{{\bf Existence Criteria and Efficient Algorithms }}
\msk
\centerline{{\bf for Lipschitz Selections of Low Dimensional Set-Valued Mappings}}
\vspace*{12mm}
\centerline{By~ {\sc Pavel Shvartsman}}\vspace*{3 mm}
\centerline {\it Department of Mathematics, Technion - Israel Institute of Technology,}\vspace*{1 mm}
\centerline{\it 32000 Haifa, Israel}\vspace*{1 mm}
\centerline{\it e-mail: pshv@technion.ac.il}
\bsk\bsk
\renewcommand{\thefootnote}{ }
\footnotetext[1]{{\it\hspace{-6mm}Math Subject
Classification:} 46E35\\
{\it Key Words and Phrases:} Set-valued mapping, Lipschitz selection, the Finiteness Principle, Helly's theorem, the core of a set-valued mapping, Hausdorff distance, balanced refinement.\smallskip
\par This research was supported by Grant No 2014055 from the United States-Israel Binational Science Foundation (BSF).}
\begin{abstract} Let $\MS=(\Mc,\rho)$ be a metric space and let $\X$ be a Banach space. Let $F$ be a set-valued mapping from $\Mc$ into the family $\Kc_m(\X)$ of all compact convex subsets of $\X$ of dimension at most $m$. The main result in our recent joint paper \cite{FS-2018}
with Charles Fefferman (which is referred to as a ``Finiteness Principle for Lipschitz selections'') provides efficient conditions for the existence of a Lipschitz selection of $F$, i.e., a Lipschitz mapping $f:\Mc\to\X$ such that $f(x)\in F(x)$ for every $x\in\Mc$. We give new alternative proofs of this result in two special cases. When $m=2$ we prove it for $X=\mathbf{R}^{2}$, and when $m=1$ we prove it for all choices of $X$. Both of these proofs make use of a simple reiteration formula for the ``core'' of a set-valued mapping $F$, i.e., for a mapping $G:\Mc\to\Kc_m(\X)$ which is Lipschitz with respect to the Hausdorff distance, and such that $G(x)\subset F(x)$ for all $x\in\Mc$.
\par We also present several constructive criteria for the existence of Lipschitz selections of set-valued mappings from $\Mc$ into the family $\HPL$ of all closed half-planes in $\RT$.
%
%
%
\end{abstract}
\renewcommand{\contentsname}{ }
\tableofcontents
\addtocontents{toc}{{\centerline{\sc{Contents}}}
\vspace*{10mm}\par}

\SECT{1. Introduction.}{1}
\addtocontents{toc}{\hspace*{3.2mm} 1. Introduction.\hfill \thepage\par\VST}
\indent\par Let  $\mfM=(\Mc,\rho)$ be a {\it pseudometric space}, i.e., suppose that the ``distance function'' $\rho:\Mc\times\Mc\to [0,+\infty]$ satisfies

$$\rho(x,x)=0, \rho(x,y)=\rho(y,x),~~~~
\text{and}~~~~
\rho(x,y)\le \rho(x,z)+\rho(z,y)$$

for all $x,y,z\in\Mc$. Note that $\rho(x,y)=0$ may hold with $x\ne y$, and $\rho(x,y)$ may be $+\infty$.
\par Let  $(\X,\|\cdot\|)$ be a real Banach space. Given a non-negative integer $m$ we let $\KM$ denote the family of all {\it non-empty compact convex subsets} $K\subset \X$ of dimension at most $m$. (We say that a convex subset of $\X$ has dimension at most $m$ if it is contained in an affine subspace of $\X$ of dimension at most $m$.) We let $$\Kc(\X)=\bigcup\{\KM: m=0,1,...\}$$ denote the family of all non-empty compact convex finite-dimensional subsets of $\X$.
\par By $\Lip(\Mc,\X)$ we denote the space of all Lipschitz mappings from $\Mc$ to $\X$ equipped with the Lipschitz seminorm
$$
\|f\|_{\Lip(\Mc,\X)}=\inf\{\,\lambda>0:\|f(x)-f(y)\|
\le\lambda\,\rho(x,y)~~~\text{for all} ~~~x,y\in\Mc\,\}.
$$
\par In this paper we study the following problem.
\begin{problem} {\em Suppose that we are given a set-valued mapping $F$ which to each point $x\in\Mc$ assigns a set $F(x)\in\KM$. A {\it selection} of $F$ is a map $f:\Mc\to \X$ such that $f(x)\in F(x)$ for all $x\in\Mc$.
\smsk
\par We want to know {\it whether there exists a selection $f$ of $F$ in the space $\Lip(\Mc,\X)$}. Such an $f$ is called a {\it Lipschitz selection} of the set-valued mapping $F:\Mc\to\KM$.
\smsk
\par If a Lipschitz selection $f$ exists, then we ask {\it how small we can take its Lipschitz seminorm}.
}
\end{problem}
\smsk

See Fig.1.
\smsk
\par The following result provides efficient conditions for the existence of a Lipschitz selection of an arbitrary set-valued mapping from a pseudometric space into the family $\KM$. We refer to it as a ``Finiteness Principle for Lipschitz selections'', or simply as a ``Finiteness Principle''.
\begin{theorem}\label{MAIN-FP} (Fefferman,Shvartsman \cite{FS-2018}) Fix $m\ge 1$. Let $(\Mc,\rho)$ be a pseudometric space, and let $F:\Mc\to\KM$ for a Banach space $\X$. Let
\bel{NMY-1}
N(m,\X)=2^{\ell(m,X)}~~~~~\text{where}~~~~~ \ell(m,X)=\min\{m+1,\dim\X\}.
\ee
\par Suppose that for every subset $\Mc'\subset\Mc$ consisting of at most $N=N(m,\X)$ points, the restriction $F|_{\Mc'}$ of $F$ to $\Mc'$ has a Lipschitz selection $f_{\Mc'}$ with Lipschitz  seminorm $\|f_{\Mc'}\|_{\Lip(\Mc',\X)}\le 1$.
\par Then $F$ has a Lipschitz selection $f$ with Lipschitz  seminorm
\begin{align}\lbl{IN-GM1}
\|f\|_{\Lip(\Mc,\X)}\le \gamma
\end{align}
where $\gamma=\gamma(m)$ is a positive constant depending only $m$.
\end{theorem}

\begin{figure}[h!]
\hspace{20mm}
\includegraphics[scale=0.6]{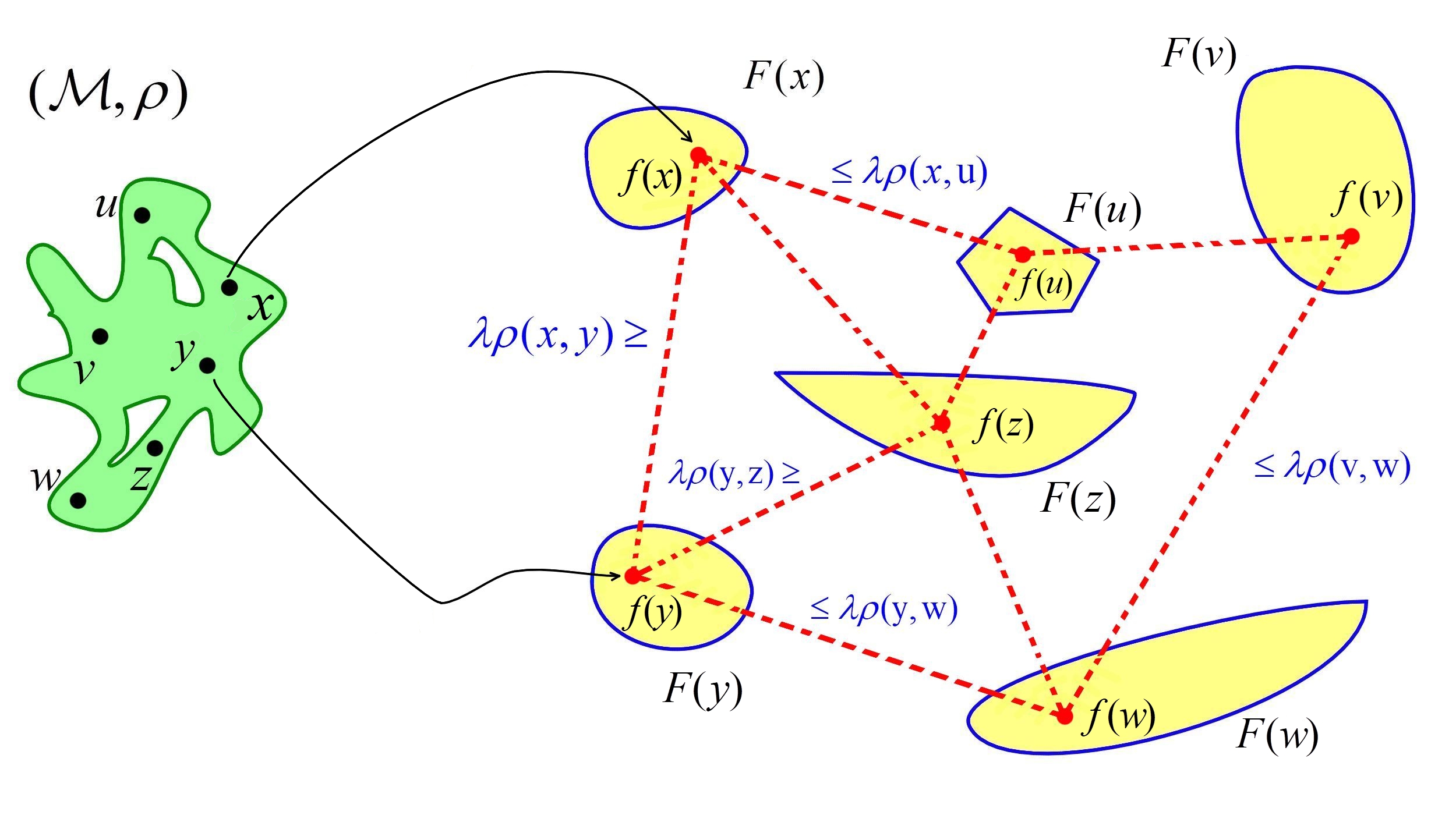}
\caption{$f:\Mc\to\RT$ is a Lipschitz selection of the set-valued mapping $F:\Mc\to\Kc(\RT)$.}
\end{figure}

\par There is an extensive literature devoted to the
Finiteness Principle for Lipschitz selection and related topics. We refer the reader to the papers \cite{Ar,AF,BL, FIL-2017,FP-2019,FS-2017,FS-2018,S-1992,S-2001,PR,PY1,PY2,
S-2002,S-2004,S-2008} and references therein for numerous results in this direction.
\par We note that the ``finiteness number'' $N(m,\X)$ in Theorem \reff{MAIN-FP} is optimal; see \cite{S-1992,S-2002}.
\par For the case of the trivial distance function $\rho\equiv 0$, Theorem \reff{MAIN-FP} agrees with the classical Helly's Theorem \cite{DGK}, except that the optimal finiteness constant for $\rho\equiv 0$ is
$$
n(m,X)=\ell(m,X)+1=\min\{m+2,\dim \X+1\}~~~~\text{in place of}~~~~N(m,X)=2^{\ell(m,X)}.
$$
Thus, Theorem \reff{MAIN-FP} may be regarded as a generalization of Helly's Theorem.
\par Our interest in Helly-type criteria for the existence of Lipschitz selections was initially motivated by some intriguing close connections of this problem with the
classical Whitney extension problem \cite{W1}, namely, the problem of characterizing those functions defined on a closed subset, say $E\subset\RN$, which are the restrictions to $E$ of $C^m$-smooth functions on $\RN$.
We refer the reader to the papers \cite{BS-1994,BS-1997,BS-2001,F-2005,F-2006,F-2009,
FI-2020,S-2008} and references therein for numerous results and techniques concerning this topic.
\par One of the main ingredients of the proof of Theorem \reff{MAIN-FP} is the construction of a special set-valued mapping $G:\Mc\to\KM$ introduced in \cite{FS-2018} which we call a {\it ``core''} of the set-valued mapping $F$. In fact each core is associated with a positive constant. Here are the relevant definitions.
\begin{definition}\lbl{D-CORE} {\em  Let $\gamma$ be a positive constant, and let $F:\Mc\to\KM$ be a set-valued mapping. A set-valued mapping $G:\Mc\to\KM$ is said to be a $\gamma$-core of $F$ if
\msk
\par (i). $G(x)\subset F(x)$ for all $x\in\Mc$;
\msk
\par (ii). $G$ is $\gamma$-Lipschitz with respect to the Hausdorff distance, i.e.,
$$
\dhf(G(x),G(y))\le \gamma\,\rho(x,y)~~~~~~\text{for all} ~~~x,y\in\Mc.
$$
}
\end{definition}
\smsk
\par We refer to a map $G$ as a {\it core of $F$} if $G$ is a $\gamma$-core of $F$ for some $\gamma>0$. See Fig. 2, 3, 4.

\bsk

\begin{figure}[h!]
\hspace{40mm}
\includegraphics[scale=1.30]{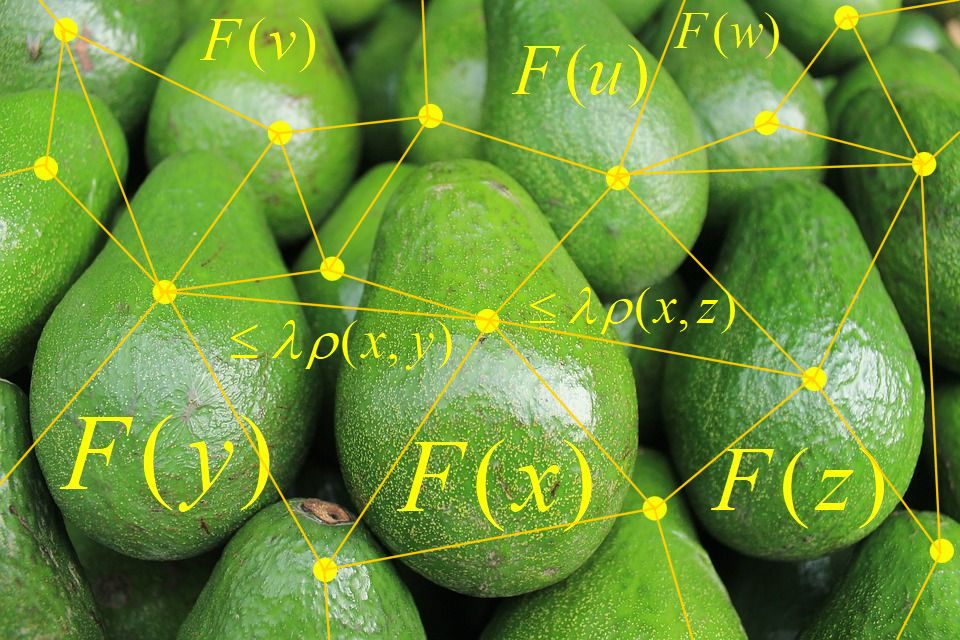}
\caption{A set-valued mapping $F$ into a family of avocados
and its Lipschitz selection
\\
\hspace*{11.5mm}
with Lipschitz seminorm at most $\lambda$.}
\end{figure}

\bsk
\begin{figure}[h!]
\hspace{40mm}
\includegraphics[scale=1]{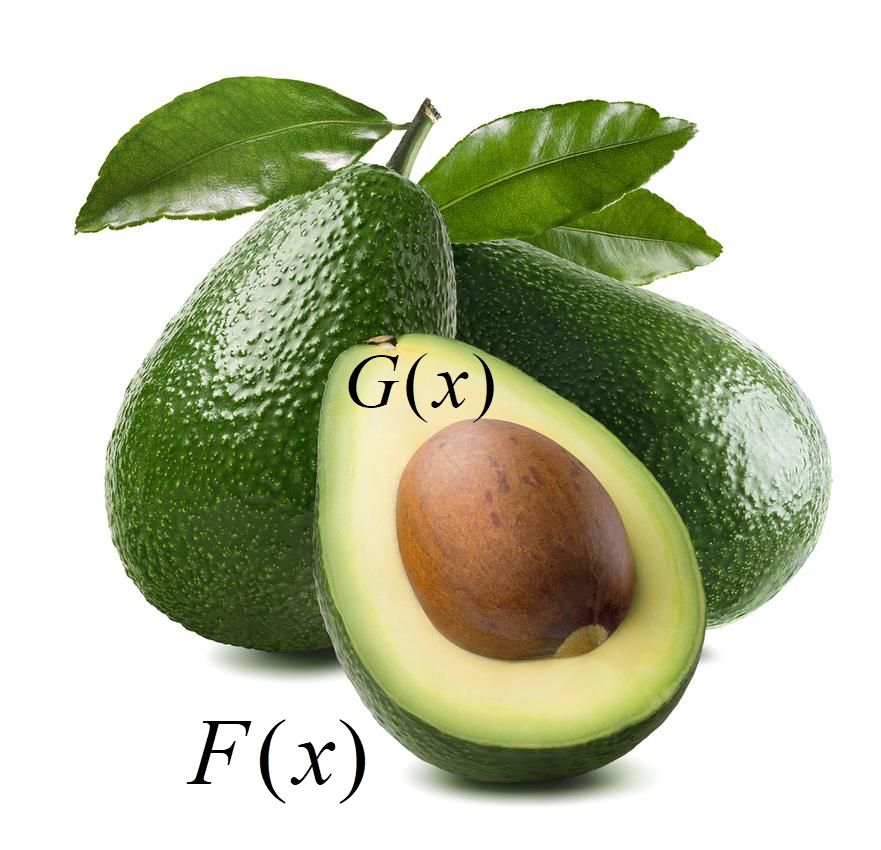}
\vspace*{-2mm}
\caption{The core $G(x)$ is a convex closed subset of $F(x)$.}
\end{figure}

\begin{figure}[h!]
\hspace{40mm}
\includegraphics[scale=0.35]{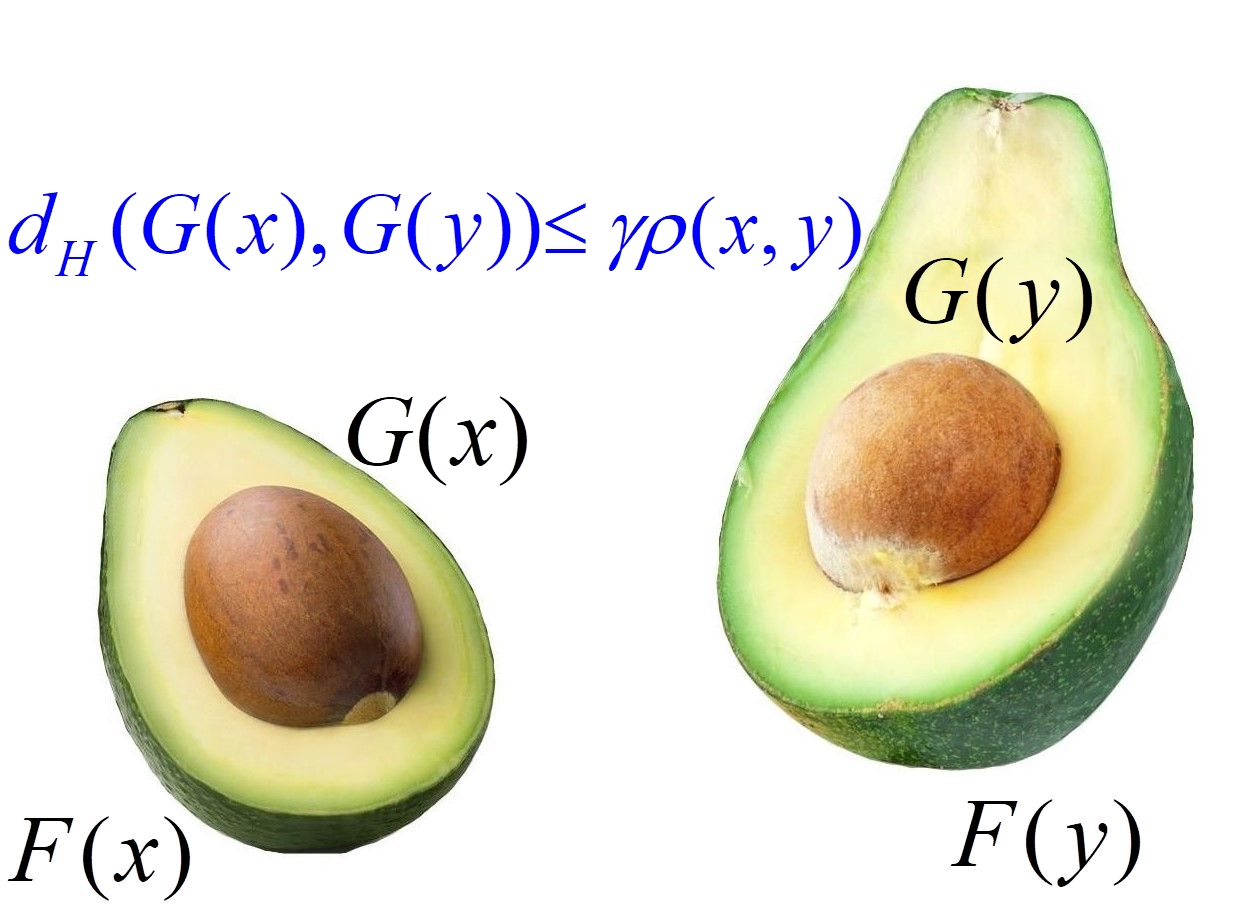}
\caption{The $\gamma$-core $G$ is $\gamma$-Lipschitz with respect to the Hausdorff distance.}
\end{figure}

\par Recall that the Hausdorff distance $\dhf(A,B)$ between two non-empty bounded sets $A,B\subset \X$ is defined as the least $r\ge 0$ such that for each $x\in A$ there exists $y\in B$ such that $\|x-y\|\le r$, and for each $x\in B$ there exists $y\in A$ such that $\|x-y\|\le r$. Thus,
\bel{HD-R}
\dhf(A,B)=\inf\{r>0:~A+\BX(0,r)\supset B~~\text{and}~~
B+\BX(0,r)\supset A\}.
\ee
\par Here and throughout this paper, for
each $x\in X$ and $r>0$, we use the standard notation $\BX(x,r)$ for the closed ball in $X$ with center $x$ and radius $r$. We also let $\BX=\BX(0,1)$ denote the unit ball in $\X$, and we write $r\BX$ to denote the ball $\BX(0,r)$.
\par In Definition \reff{D-CORE} $m$ can be any non-negative integer not exceeding the dimension of the Banach space $X$. It can happen that a core $G:\mathcal{M}\to\mathcal{K}_{m}(X)$ of a given
set-valued mapping $F:\mathcal{M}\to\mathcal{K}_{m}(X)$ in fact maps $\mathcal{M}$ into the smaller collection $\mathcal{K}_{m'}(X)$ for some integer $m'\in[0,m$). The next claim shows that the existence of some core $G:\mathcal{M}\to\mathcal{K}_{m}(X)$
for $F$ implies the existence of a (possibly different) core which maps $\mathcal{M}$ into $\mathcal{K}_{0}(X)$. Since $\mathcal{K}_{0}(X)$ is identified with $X$, that core is simply a Lipschitz selection of $F$.
\begin{claim}\label{CLS} (\cite[Section 5]{FS-2018}) Let $\gamma$ be a positive constant, let $m$ be a non-negative
integer, and let $G:\Mc\to\KM$ be a $\gamma$-core of a set-valued mapping $F:\Mc\to\KM$ for some Banach space $X$. Then $F$ has a Lipschitz selection $f:\Mc\to \X$ with $\|f\|_{\Lip(\Mc,\X)}\le C\,\gamma$ where $C=C(m)$ is a constant depending only on $m$.
\end{claim}
\par In \cite{FS-2018} we showed that this claim follows  from Definition \reff{D-CORE} and the existence of the so-called {\it ``Steiner-type point''} map $\ST:\KM\to \X$ \cite{S-2004}. See Section 2 for more detail.
\msk
\par In \cite{FS-2018} given a set-valued mapping $F:\Mc\to\KM$ satisfying the hypothesis of Theorem \reff{MAIN-FP}, we constructed a $\gamma$-core $G$ of $F$ with a positive constant $\gamma$ depending only on $m$. We produced the core $G$ using a rather delicate and complicated procedure whose main ingredients are families  of {\it Basic Convex Sets} associated with $F$, metric spaces with bounded {\it Nagata dimension}, ideas and methods of work \cite{FIL-2017} related to the case $\Mc=\RN$, and Lipschitz selections on finite metric trees. See \cite{FS-2018} for more details.
\msk
\par In the present paper we suggest and discuss a different new geometrical method for pro\-du\-cing a core of a set-valued mapping. Its main ingredient is the so-called {\it balanced refinement} of a set-valued mapping which we define as follows.
\begin{definition}\lbl{BREF} {\em  Let $\lambda\ge 0$, let $\MR$ be a pseudometric space, let $X$ be a Banach space, and let $F:\Mc\to\KM$ be a set-valued mapping for some non-negative integer $m$. For each $x\in\Mc$ we consider the subset of $F(x)$ defined by
$$
\BAL{F}{\lambda}{\rho}(x)=
\bigcap_{z\in\Mc}\,
\left[F(z)+\lambda\,\rho(x,z)\,\BX\right].
$$
\par We refer to the set-valued mapping
$\BAL{F}{\lambda}{\rho}:\Mc\to\KM\cup\{\emp\}$ as
the {\it $\lambda$-balanced refinement} of the mapping $F$.}
\end{definition}
\par We note that any Lipschitz selection $f$ of a set-valued mapping $F:\Mc\to\KM$ with Lipschitz seminorm $\|f\|_{\Lip(\Mc,\X)}\le \lambda$ is also a Lipschitz selection of the $\lambda$-balanced refinement of $F$, i.e.,
$$
f(x)\in\BAL{F}{\lambda}{\rho}(x)~~~~\text{for all}~~~~ x\in\Mc.
$$
\par Various geometrical parameters of the set $\BAL{F}{\lambda}{\rho}(x)$ (such as diameter and width, etc.) may turn out to be smaller than the same parameters for the set $F(x)$ which contains it. When attempting
to find Lipschitz selections of $F$ it may turn out to be convenient for our purposes to search for them in the more ``concentrated'' setting provided by the sets $\BAL{F}{\lambda}{\rho}(x)$. One can take this approach still further by searching in even smaller sets which can be obtained from consecutive iterations of balanced refinements of $F$, i.e. from the set functions which we describe in the following definition.
\begin{definition}\label{F-IT} {\em  Let $\ell$ be a positive integer, and let
$\vec{\lambda}=\left\{\lambda_k:1\le k\le\ell\right\}$ be a finite sequence of $\ell$ non-negative numbers $\lambda_k$. We set $F^{[0]}=F$, and, for every $x\in\Mc$ and integer $k\in[0,\ell-1]$, we define
\begin{align}\lbl{IT-F}
F^{[k+1]}(x)=\BAL{F^{[k]}}{\lambda_{k+1}}{\rho}(x)=
\bigcap_{z\in \Mc}\,
\left[F^{[k]}(z)+\lambda_{k+1}\,\rho(x,z)\,\BX\right].
\end{align}
\par We refer to the set-valued mapping $F^{[k]}:\Mc\to\KM\,\cupbig\,\{\emp\}$, $k\in[1,\ell]$, as
{\it the $k$-th order $(\vl,\rho)$-balanced refinement of $F$}.}
\end{definition}
\smsk
\par Clearly,
\begin{align}\lbl{CH-N}
F^{[k+1]}(x)\subset F^{[k]}(x)~~~~\text{on}~~~\Mc~~~ \text{for every}~~~k\in[0,\ell-1].
\end{align}
(Put $z=x$ in the right hand side of \rf{IT-F}.)

\begin{remark} {\em Of course, for each integer $k\in[1,\ell]$ the set $F^{[k]}(x)$ also depends on the sequence $\vec{\lambda}=\{\lambda_k:1\le k\le\ell\}$, on the pseudometric space $\mfM=\MR$ and the Banach space $\X$. However, in all places where we use $F^{[k]}$'s,  these objects, i.e., $\vec{\lambda}$, $\mfM$ and $\X$, are clear from the context. Therefore, in these cases, we omit any mention of $\vec{\lambda}$, $\mfM$ and $\X$ in the notation of $F^{[k]}$'s.\rbx}
\end{remark}
\smsk
\par We formulate the following
\begin{conjecture}\label{BR-IT} Let $\MR$ be a pseudometric space, and let $\X$ be a Banach space. Let $m$ be a fixed positive integer and (as in the formula \rf{NMY-1}
of Theorem \reff{MAIN-FP}) let $N(m,\X)$ denote the ``finiteness number'' $N(m,\X)=2^{\ell}$ where $\ell=\ell(m,X)=\min\{m+1,\dim\X\}$.
\par There exist a constant $\gamma\ge 1$ and a sequence
$\vec{\lambda}=\{\lambda_k:1\le k\le\ell\}$ of $\ell$ numbers $\lambda_{k}$ all satisfying $\lambda_{k}\ge1$ such that the following holds:
\smsk
\par Let $F:\Mc\to\KM$ be a set-valued mapping such that, for every $\Mc'\subset\Mc$ with $\#\Mc'\le N(m,\X)$, the restriction $F|_{\Mc'}$ of $F$ to $\Mc'$ has a Lipschitz selection $f_{\Mc'}:\Mc'\to \X$ with Lipschitz seminorm $\|f_{\Mc'}\|_{\Lip(\Mc',\X)}\le 1$.
\par Then the $\ell-th$ order balanced refinement
of the mapping $F$, namely the set-valued mapping $F^{[\ell]}:\Mc\to\KM$ is a $\gamma$-core of $F$.
\smsk
\par Here $F^{[\ell]}$ is defined as in Definition \reff{F-IT} using the particular sequence $\vec{\lambda}$.
\end{conjecture}
\par Our main results, Theorem \reff{MAIN-RT} and Theorem \reff{X-LSGM} below, state that Conjecture \reff{BR-IT} holds in two special cases, when either (i) $m=2$ and $\dim X=2$, or (ii) $m=1$ and $X$ is an arbitrary Banach space. Note that in both of these cases the above mentioned  finiteness number $N(m,\X)$ equals $4$.
\begin{theorem}\label{MAIN-RT} Let $\mfM=\MR$ be a pseudometric space, and let $X$ be a two dimensional Banach space. Let $m=2$ so that the number $\ell(m,X)=2$. In this case Conjecture \reff{BR-IT} holds for every $\lambda_1,\lambda_2$ and $\gamma$ such that
\begin{align}\lbl{GM-FN}
\lambda_1\ge e(\mfM,X), ~~~~~~\lambda_2\ge 3\lambda_1,
~~~~~~
\gamma\ge \lambda_2\,(3\lambda_2+\lambda_1)^2/
(\lambda_2-\lambda_1)^2.
\end{align}
Here $e(\mfM,X)$ denotes the Lipschitz extension constant of $X$ with respect to $\mfM$. (See Definition \reff{LIP-C}.)
\par Thus, the following statement is true: Let $F:\Mc\to\Kc(\X)$ be a set-valued mapping from a pseudometric space $(\Mc,\rho)$ into the family $\Kc(X)$ of all non-empty convex compact subsets of $\X$. Given $x\in\Mc$ let
\begin{align}\lbl{F12-LTI}
F^{[1]}(x)=
\bigcap_{z\in\Mc}\,
\left[F(z)+\lambda_1\,\rho(x,z)\BX\right],~~~~
F^{[2]}(x)=\bigcap_{z\in\Mc}\,
\left[F^{[1]}(z)+\lambda_2\,\rho(x,z)\BX\right].
\end{align}
\par Suppose that for every subset $\Mc'\subset\Mc$ with $\#\Mc'\le 4$, the restriction $F|_{\Mc'}$ of $F$ to $\Mc'$ has a Lipschitz selection with Lipschitz seminorm at most $1$.
\par Then for every $\lambda_1,\lambda_2$ and $\gamma$ satisfying \rf{GM-FN} the set
\begin{align}\lbl{F2-NEM}
F^{[2]}(x)\ne\emp~~~~\text{for every}~~~~x\in\Mc.
\end{align}
Furthermore,
\begin{align}\lbl{HD-RT}
\dhf(F^{[2]}(x),F^{[2]}(y))\le \gamma\,\rho(x,y)
~~~~\text{for every}~~~~x,y\in\Mc.
\end{align}
\par If $X$ is a Euclidean two dimensional space, \rf{F2-NEM} and  \rf{HD-RT} hold when \rf{GM-FN} is replaced by the weaker requirements that
\begin{align}\lbl{X-HSP}
\lambda_1\ge e(\mfM,X), ~~~~~~\lambda_2\ge 3\lambda_1,
~~~~~~
\gamma\ge \,\lambda_2\left\{1+2\lambda_2\,\slbig
\left(\lambda_2^2-\lambda_1^2\right)^{\frac12}\right\}^2.
\end{align}
\end{theorem}
\smsk
\smsk
\par In particular, in Section 3 we show that the mapping $F^{[2]}$ satisfies \rf{F2-NEM} and \rf{HD-RT} whenever $X$ is an {\it arbitrary} two dimensional Banach space and $\lambda_1=4/3$, $\lambda_2=4$, $\gamma=100$. If $X$ is also Euclidean, then one can set $\lambda_1=4/\pi$, $\lambda_2=12/\pi$ and $\gamma=38$. Furthermore, we prove that if $\Mc$ is a subset of a Euclidean space $E$, $\rho$ is the Euclidean metric in $E$, and $\X$ is a two dimensional Euclidean space, then properties \rf{F2-NEM} and \rf{HD-RT} hold for $\lambda_1=1$, $\lambda_2=3$, and $\gamma=25$.
\msk
\par In Section 6 we prove Theorem \reff{LTI-M} which improves the result of Theorem \reff{MAIN-RT} for the space $X=\LTI$, i.e., for $\RT$ equipped with the uniform norm
$$
\|x\|=\max\{|x_1|,|x_2|\},~~~~x=(x_1,x_2).
$$
More specifically, we show that in this case properties \rf{F2-NEM} and \rf{HD-RT} hold provided
$$
\lambda_1\ge 1,~~~\lambda_2\ge 3\lambda_1,~~~~\text{and}~~~~ \gamma\ge \lambda_2\, (3\lambda_2+\lambda_1)/(\lambda_2-\lambda_1).
$$
In particular, these properties hold for $\lambda_1=1$, $\lambda_2=3$ and $\gamma=15$.
\smsk
\par Let us now explicitly formulate the above mentioned
second main result of the paper. We prove it in Section 4. It deals with set-valued mappings from a pseudometric space into the family $\Kc_1(\X)$ of all {\it bounded closed line segments} of an arbitrary Banach space $\X$.
\begin{theorem}\label{X-LSGM} Let $\MR$ be a pseudometric space. Let $m=1$ and let $X$ be a Banach space with $\dim X>1$; thus, $\ell(m,X)=2$, see \rf{NMY-1}. In this case Conjecture \reff{BR-IT} holds for every  $\lambda_1,\lambda_2$ and $\gamma$ such that
\begin{align}\lbl{GM-FN-1}
\lambda_1\ge 1, ~~~~~~\lambda_2\ge 3\lambda_1,
~~~~~~
\gamma\ge \lambda_2\,(3\lambda_2+\lambda_1)/
(\lambda_2-\lambda_1).
\end{align}
\par Thus, the following statement is true: Let $F:\Mc\to\Kc_1(\X)$ be a set-valued mapping such that for every subset $\Mc'\subset\Mc$ with $\#\Mc'\le 4$, the restriction $F|_{\Mc'}$ of $F$ to $\Mc'$ has a Lipschitz selection with Lipschitz seminorm at most $1$.
\par Let $F^{[2]}$ be the mapping defined by \rf{F12-LTI}. Then properties \rf{F2-NEM} and \rf{HD-RT} hold whenever $\lambda_1$, $\lambda_2$ and $\gamma$ satisfy \rf{GM-FN-1}. In particular, one can set $\lambda_1=1$, $\lambda_2=3$ and $\gamma=15$.
\par If $X$ is a Euclidean space, the same statement is also true whenever, instead of \rf{GM-FN-1}, $\lambda_1,\lambda_2$ and $\gamma$ satisfy the weaker  condition
\begin{align}\lbl{X1-H}
\lambda_1\ge 1, ~~~~~~\lambda_2\ge 3\lambda_1,
~~~~~~
\gamma\ge \,\lambda_2+2\lambda_2^2\,\slbig
\left(\lambda_2^2-\lambda_1^2\right)^{\frac12}.
\end{align}
\par In particular, in this case, \rf{F2-NEM} and \rf{HD-RT} hold whenever $\lambda_1=1$, $\lambda_2=3$ and $\gamma=10$.
\end{theorem}
\smsk

\par In Section 5 we note that Conjecture \reff{BR-IT} also holds for a one dimensional space $X$ and $m=1$. In this case the statement of the conjecture is true for every $\lambda_1\ge 1$ and $\gamma\ge 1$. See Proposition \reff{X-1DIM}.

\par Note that Theorem \reff{MAIN-RT} tells us that for every set-valued mapping $F$ satisfying the hypothesis of this theorem, the mapping $F^{[2]}$ determined by \rf{F12-LTI} with $\lambda_1=4/3$ and $\lambda_2=4$ provides a $\gamma$-core of $F$ with $\gamma=100$. (See Definition \reff{D-CORE}.) In turn, Theorem \reff{X-LSGM} states that the mapping $F^{[2]}$ corresponding to the parameters $\lambda_1=1$ and $\lambda_2=3$ is a $15$-core of any $F$ satisfying the conditions of this theorem.
\smsk
\par We note that the proofs of Theorem
\reff{MAIN-RT} and Theorem \reff{X-LSGM} rely on Helly's Intersection Theorem and a series of auxiliary results about neighborhoods of intersections of convex sets. See Section 2.

\begin{remark} {\em Let us compare Conjecture \reff{BR-IT} (and Theorems \reff{MAIN-RT} and \reff{X-LSGM}) with
the Finiteness Principle (FP) formulated in Theorem \reff{MAIN-FP}. First we note that FP is invariant with respect to the transition to an equivalent norm on $X$, while the statement of Conjecture \reff{BR-IT} is not.
\par To express this more precisely, let $\|\cdot\|_1$ and $\|\cdot\|_2$ be two equivalent norms on $X$, i.e., suppose that for some $\alpha\ge 1$  the following inequality

$$
(1/\alpha)\,\|\cdot\|_1\le\|\cdot\|_2\le \alpha\,\|\cdot\|_1
$$

holds. Clearly, if FP holds for  $(X,\|\cdot\|_1)$ then it immediately holds also for $(X,\|\cdot\|_2)$ (with the constant $\alpha^2\gamma$ in \rf{IN-GM1} instead of $\gamma$). However the validity of Conjecture \reff{BR-IT} for the norm $\|\cdot\|_1$ does not imply its validity for an equivalent norm $\|\cdot\|_2$ on $X$ (at least we do not see any obvious way for obtaining such an implication). For example, the validity of Conjecture \reff{BR-IT} in $\ell^n_\infty$ (i.e., $\RN$ equipped with the uniform norm) does not automatically imply its validity in the space $\ell^n_2$ (i.e., $\RN$ with the Euclidean norm).
\smsk
\par  We also note the following: in a certain sense,
the result of Theorem \reff{MAIN-RT} is ``stronger'' than Theorem \reff{MAIN-FP} (i.e., FP for the case of a two dimensional Banach space $X$). Indeed, in this case, the hypotheses of FP and Theorem \reff{MAIN-RT} coincide. Moreover, Theorem \reff{MAIN-RT} ensures that the set-valued mapping $F^{[2]}$ is a core of $F$. This
property of $F^{[2]}$ implies, via arguments in \cite{FS-2018} that the function
$$
f(x)=\ST\,(F^{[2]})\,(x),~~~~~~x\in\Mc,
$$
is a Lipschitz selection of $F$. Here $\ST:\KM\to\X$  is the Steiner-type point map \cite{S-2004}.
\par Thus, FP (in the two dimensional case) follows immediately from Theorem \reff{MAIN-RT}. However, it is absolutely unclear how the statement of Theorem \reff{MAIN-RT} can be deduced from FP. I would like to thank Charles Fefferman who kindly drew my attention to this interesting fact.\rbx}
\end{remark}
\par Let us reformulate Conjecture \reff{BR-IT} in a way  which {\it does not require the use of the notion of a core of a set-valued mapping}. We recall that the mapping $F^{[\ell]}:\Mc\to \KM$ which appears in Conjecture \reff{BR-IT} is a $\gamma$-core of $F$ if
$$
\dhf(F^{[\ell]}(x),F^{[\ell]}(y))\le \gamma\, \rho(x,y)~~~\text{for all}~~~x,y\in\Mc.
$$
See part (ii) of Definition \reff{D-CORE}. Hence, given $x\in\Mc$,
\begin{align}\lbl{IM-L}
F^{[\ell]}(x)\subset F^{[\ell]}(y)+
\gamma\,\rho(x,y) \BX~~~\text{for every}~~~y\in\Mc.
\end{align}
\par We also recall that
$$
F^{[\ell+1]}(x)=\BAL{F^{[\ell]}}{\gamma}{\rho}(x)=
\bigcap_{y\in \Mc}\,
\left[F^{[\ell]}(y)+\gamma\,\rho(x,y)\,\BX\right].
$$
See \rf{IT-F}. This and \rf{IM-L} imply the inclusion
$F^{[\ell+1]}(x)\supset F^{[\ell]}(x)$, $x\in\Mc$. On the other hand, \rf{CH-N} tells us that

$$
F^{[\ell+1]}(x)\subset F^{[\ell]}(x)~~~~~\text{proving that}~~~~F^{[\ell+1]}=F^{[\ell]}~~~\text{on}~~~\Mc.
$$
\smsk
\par These observations enable us to reformulate Conjecture \reff{BR-IT} as follows.
\begin{conjecture}\label{ST-IT} Let $\MR$ be a pseudometric space, and let $\X$ be a Banach space. Let $m$ be a fixed positive integer and let $\ell=\ell(m,X)$, see\rf{NMY-1}.
\par There exists a sequence $\vec{\lambda}=\left\{\lambda_{k}:1\le k\le\ell+1\right\}$ of $\ell+1$ numbers $\lambda_{k}$ all satisfying $\lambda_{k}\ge1$ such that, for every set-valued mapping $F:\Mc\to\KM$ satisfying the hypothesis of the Finiteness Principle (Theorem \reff{MAIN-FP}), the family $\{F^{[k]}: k=1,...,\ell+1\}$ of  set-valued mappings constructed by formula \rf{IT-F} has the following property:
\begin{align}\lbl{ST-PR}
F^{[\ell]}(x)\ne\emp~~~~~\text{and}
~~~~~F^{[\ell+1]}(x)=F^{[\ell]}(x)~~~~ \text{for all} ~~~
x\in\Mc.
\end{align}
\end{conjecture}
\par We refer to \rf{ST-PR} as a {\it Stabilization Property} of balanced refinements.
\smsk
\par Thus, Theorem \reff{MAIN-RT} and Theorem \reff{X-LSGM} tell us that a Stabilization Property of ba\-lanced refinements holds whenever $\dim\X=2$ or $m=1$ (and $\X$ is an arbitrary). More specifically, Theorem \reff{MAIN-RT} shows that if $m=2$ and $\dim \X=2$, Conjecture \reff{ST-IT} holds with $\ell=2$ and $\vec{\lambda}=\{4/3,4,10^2\}$.
\par In other words, in this case, $F^{[2]}(x)\ne\emp$ for each $x\in\Mc$ and $F^{[3]}=F^{[2]}$ on $\Mc$. In turn, Theorem \reff{X-LSGM} states that the same property holds whenever $X$ is an arbitrary Banach space, $m=1$, and $\vec{\lambda}=\{1,3,15\}$.
\msk

\par In Sections 7 and 8 we present several explicit criteria for the existence of Lipschitz selections of set-valued mappings from a pseudometric space $\mfM=\MR$ into the family $\CRT$ of all convex closed subsets of $\RT$. These criteria develop the ideas and methods of a constructive criterion for Lipschitz selections in $\RT$ given in the paper \cite{S-2002}. Let us recall this result.
\par Let
$$
\X=\LTI~~~~~\text{and let}~~~~Q_0=[-1,1]\times[-1,1]
$$
be the unit ball of $\X$. Given a set-valued mapping $F:\Mc\to\Kc(\RT)$, a positive constant $\lambda$ and elements $x,x'\in\Mc$, we introduce a set
$$
\RL[x,x':\lambda]=
\HR[F(x)\cap\{F(x')+\lambda\,\rho(x,x') Q_0\}].~~~~~~~~~\text{(See Fig. 5.)}
$$
\smsk
\begin{figure}[h!]
\hspace{20mm}
\includegraphics[scale=0.3]{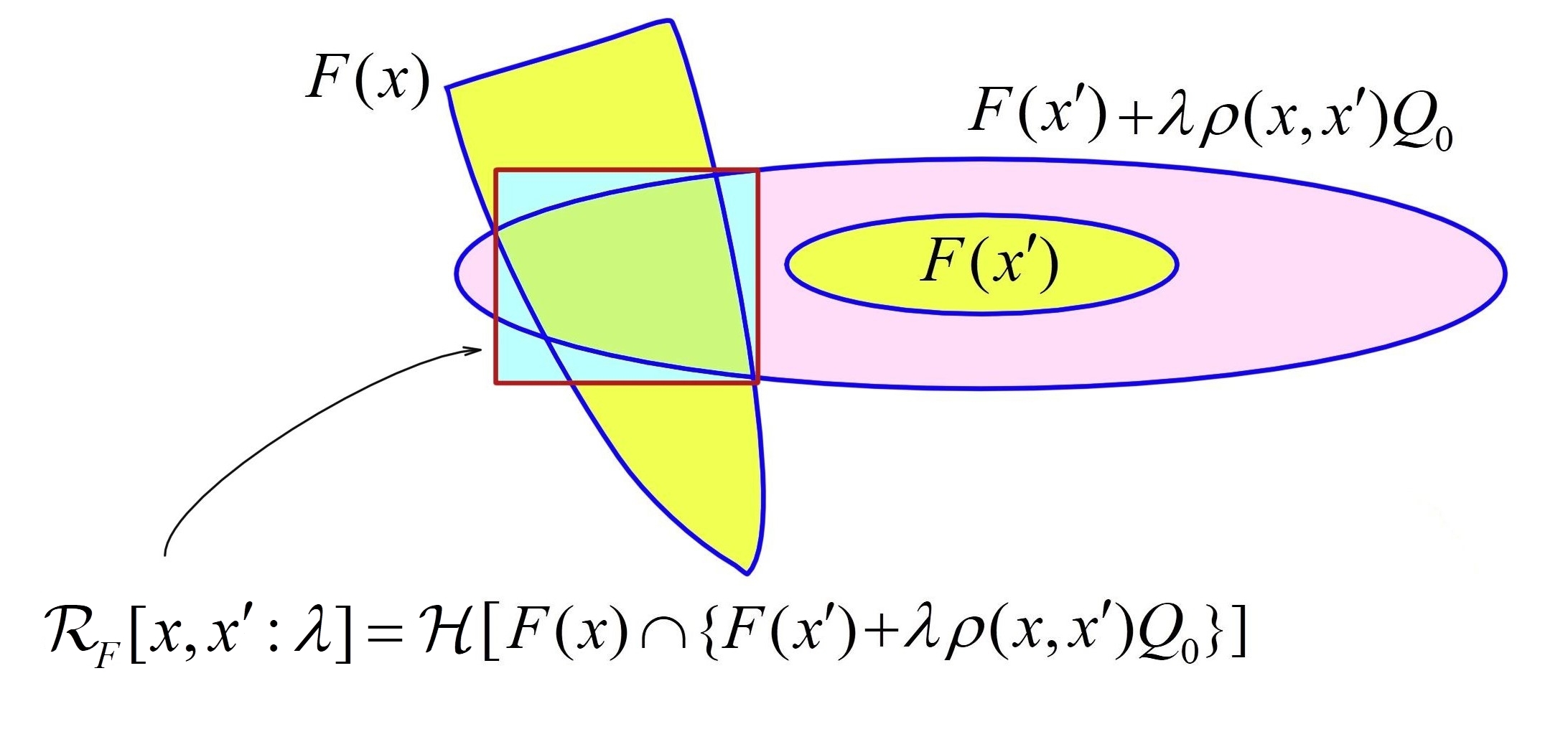}
\caption{The rectangular $\RL[x,x':\lambda]$ for $x,x'\in\Mc$ and $\lambda>0$.}
\end{figure}
\par Here, given a set $S\subset\RT$, by $\HR[S]$ we denote the smallest rectangular with sides parallel to the coordinate axes containing $S$.
We also set
\bel{OM-L}
|F|_{\mfM,\X}=\inf\{\,\|f\|_{\Lip(\Mc,\X)}: f~~\text{is a Lipschitz selection of}~~F\}.
\ee
\begin{theorem}\lbl{CR-LS1} A set-valued mapping $F:\Mc\to \Kc(\RT)$ has a Lipschitz selection if and only if there exists a constant $\lambda>0$ such that  the following conditions are satisfied:
\par (i) $\dist(F(x),F(y))\le\lambda\,\rho(x,y)$ for all $x,y\in\Mc$;
\par (ii) for every $x,x',y,y'\in\Mc$, we have
\bel{DRXY}
\dist\left(\,\RL[x,x':\lambda],\RL[y,y':\lambda]\,\right)
\le\lambda\,\rho(x,y).~~~~~~~\text{(See Fig. 6.)}
\ee
\par Furthermore, in this case
\bel{OP-NF}
\inf\lambda\le \,|F|_{\mfM,\LTI}\le 8\inf\lambda.
\ee
\end{theorem}

\begin{figure}[h!]
\hspace{5mm}
\includegraphics[scale=0.64]{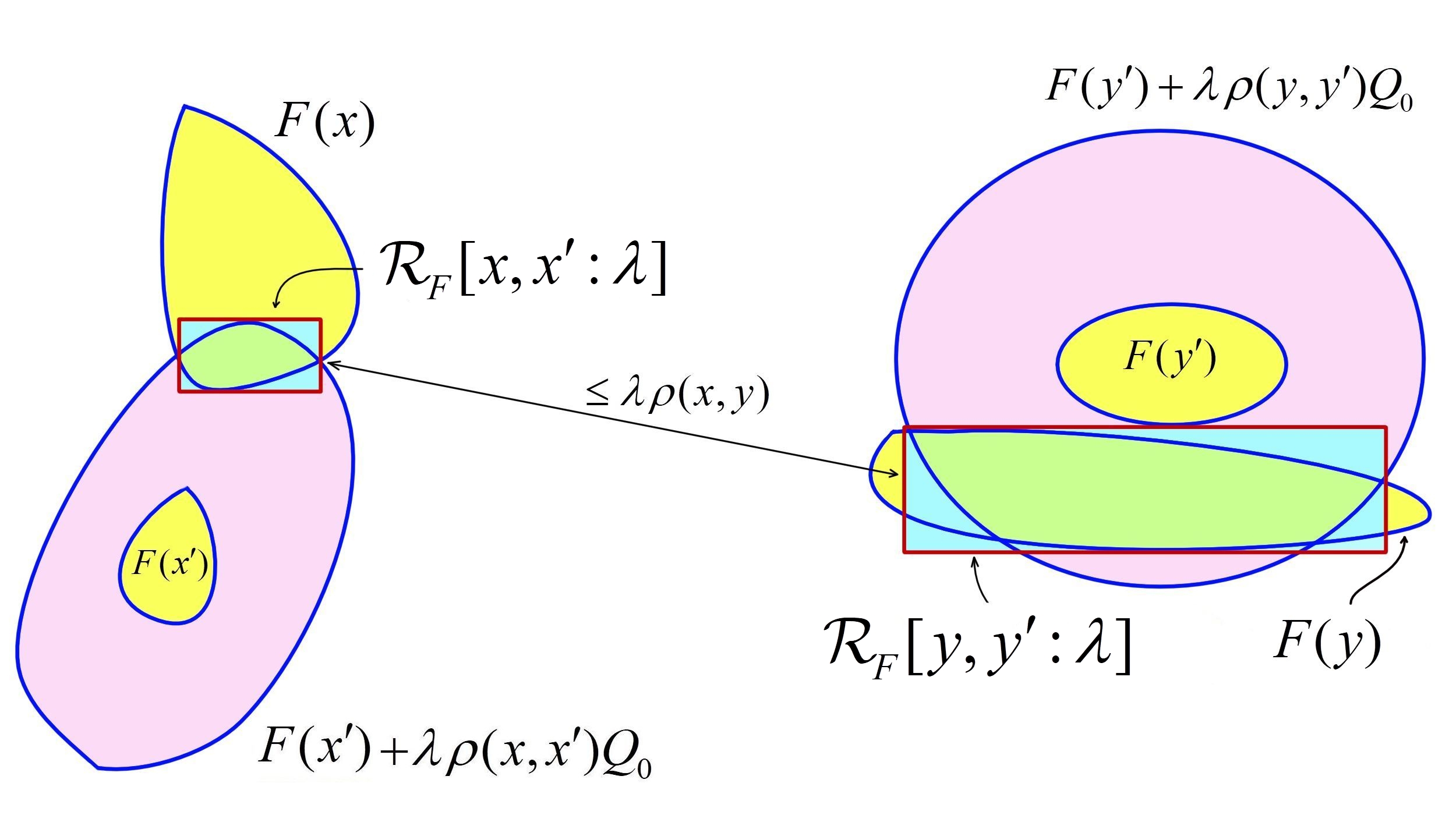}
\hspace*{-20mm}
\caption{The Lipschitz selection criterion in $\RT$.}
\end{figure}

\par For the convenience of the reader, in Section 7 we give a refined version of the proof of this theorem.
\msk
\par In Sections 8 and 9 we study Lipschitz selections of set-valued mappings which take their values in the family
$\HPL$ of all closed half-planes in $\RT$.
\par Let $\SO$ be the unit circle in $\RT$, and let $n:\Mc\to\SO$ and $\al:\Mc\to\R$ be two mappings defined on $\Mc$. These mappings determine a set-valued mapping $F:\Mc\to \HPL$ defined by
\bel{F-NAL-I}
F(x)=\{a\in\RT:\ip{a,n(x)}+\al(x)\le 0\},~~~~x\in\Mc.
\ee
(Here, given $a=(a_1,a_2)$, $n(x)=(n_1(x),n_2(x))\in\RT$, by $\ip{a,n(x)}=a_1 n_1(x)+a_2 n_2(x)$ we denote the standard inner product in $\RT$.) Thus, for each $x\in\Mc$, the set $F(x)$ is a half-plane in $\RT$ whose boundary is a straight line
$\ell_F(x)=\{a\in\RT:\ip{a,n(x)}+\al(x)=0\}$. The unit vector $n(x)$ is directed outside of the half-plane $F(x)$ and orthogonal to the line  $\ell_F(x)$.
\par Given $x,y\in\Mc$ such that $n(x)\nparallel n(y)$ we set $w(x,y:F)=\ell_F(x)\cap\ell_F(y)$. (In Section 8 we
give explicit formulae for the coordinates of the point
$w(x,y:F)=(w_1(x,y:F),w_2(x,y:F))$. See \rf{FW-XY}.)
Finally, by $\DXY$ we denote the determinant
$$
\DXY=\det\left(
\begin{array}{ll}
n_1(x)& n_1(y)\vspace{2mm}\\
n_2(x)& n_2(y)
\end{array}
\right)=n_1(x)\,n_2(y)-n_2(x)\,n_1(y).
$$
\begin{theorem}\lbl{CR-1L} Let $F:\Mc\to\HPL$ be a set-valued mapping defined by \rf{F-NAL-I}. Assume that either $\Mc$ is finite or there exist elements $x_1,...x_m\in\Mc$ such that the interior of convex hull of points $n(x_1),...,n(x_m)$ contains $0$.
\par The set-valued mapping $F$ has a Lipschitz selection if and only if there exists a constant $\lambda>0$ such that the following two conditions hold:
\msk
\par ($\bigstar 1$) $\al(x)+\al(y)\le \lambda\,\rho(x,y)$ for every $x,y\in\Mc$ such that $n(y)=-n(x)$;
\smsk
\par ($\bigstar 2$) For every $x,x',y,y'\in\Mc$ such that $n(x)\nparallel n(x')$, $n(y)\nparallel n(y')$, we have
\be
&&w_1(x,x':F)-w_1(y,y':F)\nn\\
&\le& \lambda\left\{
\frac{\rho(x,x')}{|\DXP|}\min\{|n_2(x)|,|n_2(x')|\}+
\frac{\rho(y,y')}{|\DYP|}\min\{|n_2(y)|,|n_2(y')|\}+
\rho(x,y)\right\}
\lbl{W-1Q}
\ee
provided
$$
n_2(x)\,n_2(x')\le 0,~n_1(x)+n_1(x')\le 0~~~~~\text{and}~~~~~n_2(y)\,n_2(y')\le 0, ~n_1(y)+n_1(y')\ge 0,
$$
and
\be
&&w_2(x,x':F)-w_2(y,y':F)\nn\\
&\le& \lambda\left\{
\frac{\rho(x,x')}{|\DXP|}\min\{|n_1(x)|,|n_1(x')|\}+
\frac{\rho(y,y')}{|\DYP|}\min\{|n_1(y)|,|n_1(y')|\}+
\rho(x,y)\right\}
\lbl{W-2Q}
\ee
provided
$$
n_1(x)\,n_1(x')\le 0,~ n_2(x)+n_2(x')\le 0,~~~~~\text{and}~~~~~n_1(y)\,n_1(y')\le 0,~ n_2(y)+n_2(y')\ge 0.
$$
\par Furthermore,
\bel{L-OP-8}
\tfrac{1}{\sqrt{2}}\,\inf\lambda
\le \,|F|_{\mfM,\LTI} \le 8\inf\lambda.
\ee
\end{theorem}
\smsk
\par Necessary and sufficient conditions for the existence of a Lipschitz selection given in Theorem \reff{CR-1L} involve Cartesian coordinates of certain geometric objects determined by the set-valued mapping $F$. Theorem \reff{CF-CR} below presents another explicit criterion for Lipschitz selections of $F$. This criterion formulates in terms of geometrical objects which depend only on $F$ and independent of the coordinate system in $\RT$. We refer to this criterion as a ``coordinate-free'' Lipschitz selection criterion.
\par Let us prepare the ingredients that are needed to formulate Theorem \reff{CF-CR}. Let $F:\Mc\to\HPL$ be a set-valued mapping defined by formula \rf{F-NAL-I}. Given $x,y\in\Mc$, we let $\AF(x,y)\in[0,\pi/2]$ denote the angle between the boundaries of $F(x)$ and $F(y)$, i.e., between the straight lines $\ell_F(x)$ and $\ell_F(y)$. Given a set $\Mc'\subset\Mc$, by $\diam_\rho(\Mc')$ we denote the diameter of $\Mc$ in $\MR$. Finally, we set
$0/0=0$, $a/0=+\infty$ for every $a>0$, and $\dist(\hspace{0.2mm}\emp,A)=0$ provided $A\subset\RT$.
\begin{theorem}\lbl{CF-CR} Let $\mfM=\MR$ be a pseudometric space, and let $F:\Mc\to\HPL$ be a set-valued mapping defined by \rf{F-NAL-I}. Assume that either $\Mc$ is finite or there exist elements $x_1,...x_m\in\Mc$ such that the interior of convex hull of points $n(x_1),...,n(x_m)$ contains $0$.
\par The mapping $F$ has a Lipschitz selection $f:\Mc\to\LTI$ if and only if there exists a constant
$\lambda>0$ such that for every four elements $x,x',y,y'\in\Mc$ the following inequality
\bel{MC-2}
\dist(F(x)\cap F(x'),F(y)\cap F(y'))\le
\lambda\left\{
\frac{\rho(x,x')}{\SXP}+\frac{\rho(y,y')}{\SYP}+
\diam_\rho\{x,x',y,y'\}
\right\}
\ee
holds. Furthermore,
$$
\tfrac{1}{\sqrt{2}}\,\inf\lambda
\le \,|F|_{\mfM,\LTI}
\le \gamma\,\inf\lambda.
$$
Here $\gamma>0$ is an absolute constant, $\gamma\le 5\cdot 10^5$.
\end{theorem}
\par In the next version of the paper, we will present several results related to the following problem formulated by C. Fefferman \cite{F-2019}:
\begin{problem}  {\em Let $(\Mc,\rho)$ be an $N$-point metric space. For each $x \in \Mc$, let $F(x) \subset \R^D$ be a convex polytope.
\par {\it How can one compute a map $f: \Mc\to\R^D$ such that $f(x)\in F(x)$ for all $x\in \Mc$, with Lipschitz norm as small as possible up to a factor $C(D)?$
\smsk
\par This is a big ill-conditioned linear programming problem. Can we do better than just applying general-purpose linear programming? How does the work of an optimal algorithm scale with the number of points $N?$}}
\end{problem}
\par Let $\mfM=\MR$ be an $N$-point pseudometric space (i.e., $N=\#\Mc$). We will present several efficient algorithms for Lipschitz selections of set-valued mappings from $\Mc$ into the family $\HPL$  of all closed half-planes in $\RT$. These algorithms rely on the methods of proofs of the constructive criteria for Lipschitz selections given in Sections 7-9.
\par In particular, we will exhibit an algorithm which, given a set-valued mapping $F:\Mc\to\HPL$ computes the order of magnitude of the quantity $|F|_{\mfM,\RT}$ (i.e., the Lipschitz seminorm of an optimal Lipschitz selection of $F$, see \rf{OM-L}), and a nearly optimal Lipschitz selection $f$ of $F$ {\it using work at most $C N^3$ and storage at most $CN$}. Here $C$ is an absolute constant.
\par Also, we will present a result related to a set-valued mapping $F$ from $\Mc$ into the family $\Kc_1(\R^M)$ of all bounded closed line segments of $\R^M$. (Here $M$ is a positive integer). In this case, we will exhibit an algorithm which computes the order of magnitude of $|F|_{\mfM,\R^M}$ and a nearly optimal Lipschitz selection $f$ of $F$ {\it using work at most $C(M+N^3)$ and storage at most $C(M+N)$}.
\smsk
\par The main ingredients of the proofs of these results are linear-time algorithms for linear programming in $\R^3$ due to Megiddo \cite{M}, and Lipschitz selection criteria of Theorems \reff{CR-LS-RI} and \reff{BS-LS-E}.

\bsk
\par {\bf Acknowledgements.} I am very thankful to Michael Cwikel for useful suggestions and remarks. I am also very grateful to Charles Fefferman for stimulating discussions and valuable advice.
\par The results obtained in Sections 2-6 of this paper were presented at the 12th Whitney Problems Workshop, August 2019, the University of Texas at Austin, TX. I am very thankful to all participants of that workshop for valuable conversations and useful remarks.

%
\SECT{2. Neighborhoods of intersections of convex sets in a Banach space.}{2}
\addtocontents{toc}{2. Neighborhoods of intersections of convex sets in a Banach space.\hfill \thepage\par\VST}
\indent

\par We first need to fix some notation. Let $(\X,\|\cdot\|)$ be a Banach space. We write
$$
\diam A=\sup\{\|a-b\|:a,b\in A\}~~~~~\text{and}~~~~~
\dist(A',A'')=\inf\{\|a'-a''\|:a'\in A',~a''\in A''\}
$$
to denote the diameter of a set $A\subset \X$ and the distance between sets $A',A''\subset\X$ respectively. For $x\in\X$ we also set $\dist(x,A)=\dist(\{x\},A)$, and put $\dist(\emp,A)=0$ provided $A$ is an arbitrary (possibly empty) subset of $X$. If $A\subset \X$ is finite, by $\#A$ we denote the number of elements of $A$.
\par Given non-empty sets $A,B\subset\X$ we let $A + B=\{a + b: a\in A, b\in B\}$ denote the Minkowski sum of these sets. Given a non-negative real number $\lambda$ by $\lambda A$ we denote the set $\lambda A =\{\lambda a: a\in A\}$.
\par Given $a,b\in\X$, $a\ne b$, by $[a,b]$ we denote a closed interval (a line segment) with ends in $a$ and $b$:
$$
[a,b]=\{x\in\X:x=(1-t)\,a+t\,b, 0\le t\le 1\}.
$$
We also write $[a,a]=\{a\}$ and consider $[a,a]$ as a closed ``interval'' in $\X$. By $\CX$ we denote the family of all bounded convex closed non-empty subsets of $\X$.
\par Given a set $A\subset\R$ we put
$\min A=\{\min x:x\in A\}$ and $\max A=\{\max x:x\in A\}$ provided $A$ is a closed subset of $\R$ bounded from above or below respectively. We let
\bel{IR-DF1}
\Ic(\R)=\{[a,b]:a,b\in\R, a\le b\}\cup
\{[a,+\infty):a\in\R\}\cup
\{(-\infty,b]:b\in\R\}\cup\{\R\}
\ee
denote the family of all closed intervals in $\R$ (bounded or unbounded).
We write $[x]_+$ for the positive part of the real $x$, i.e., $[x]_+=\max\{x,0\}$. We set $\frac{0}{0}=0$ and $\frac{a}{0}=+\infty$ for $a>0$.
\msk
\par Sometimes, given a set $\Mc$, we will be looking simultaneously at two distinct pseudometrics on $\Mc$, say $\rho$ and $\delta$. In this case we will speak of a $\rho$-Lipschitz selection and $\rho$-Lipschitz seminorm, or a $\delta$-Lipschitz selection and $\delta$-Lipschitz seminorm to make clear which pseudometric we are using. Furthermore, given a mapping $f:\Mc\to X$ we will write
$\|f\|_{\Lip((\Mc;\rho),\X)}$ or  $\|f\|_{\Lip((\Mc;\delta),\X)}$ to denote the Lipschitz seminorm of $f$ with respect to the pseudometric $\rho$ or $\delta$ respectively.
\smsk
\par We let $\LNI$ denote the space $\RN$ equipped with the uniform norm $\|x\|_\infty=\max\{|x_i|:i=1,...,n\}$ for $x=(x_1,...,x_n)\in\RN$. By $\LNT$ we denote $\RN$ equipped with the Euclidean norm
$\|x\|_2=\left(\sum_{i=1}^n x_i^2\right)^{1/2}$.
\smsk
\par  By $Ox_1=\{x=(t,0):t\in\R\}$ and  $Ox_2=\{x=(0,t):t\in\R\}$ we denote coordinate axes in $\RT$. Let
$$
B_0=\{a\in\RT:\|a\|_{\ell^2_2}\le 1\}
~~~~\text{and}~~~~
\SO=\{a=(a_1,a_2)\in\RT:
\|a\|_{\ell^2_2}=(a_1^2+a_2^2)^{\frac12}=1\}
$$
be the closed unit disk and the unit circle in $\RT$ respectively. Given non-zero vectors $u,v\in\RT$ we write $u\parallel v$ if $u$ and $v$ are collinear, and we write $u\nparallel v$ whenever these vectors are non-collinear. By $\theta(u,v)\in[0,2\pi)$ we denote
\bel{ANGLE}
\text{the angle of rotation from}~~~u/\|u\|_{\LTT}~~~\text{to}~~~  v/\|v\|_{\LTT}~~~\text{in the counterclockwise direction}.
\ee
(Thus, $\theta(v,u)=2\pi-\theta(u,v)$.) We refer to $\theta(u,v)$ as the angle between the vectors $u$ and $v$.
\smsk
\par Let us $\ell_1$, $\ell_2$ be two non-parallel straight lines in $\RT$ (we write $\ell_1\nparallel\ell_2$), and let  $A=\ell_1\cap \ell_2$. These two lines naturally form two angles $\vf_1,\vf_2\in[0,\pi)$, $\vf_1+\vf_2=\pi$, with  vertex at the point A. Let
\bel{ANG-L}
\vf(\ell_1,\ell_2)=\min\{\vf_1,\vf_2\};~~~~~
\text{clearly,}~~~~ \vf(\ell_1,\ell_2)\in[0,\pi/2].
\ee
\par Everywhere in this paper we refer to $\vf(\ell_1,\ell_2)$ as {\it ``the angle between straight lines $\ell_1$ and $\ell_2$''.} In other words, the angle between two non-parallel lines in $\RT$ means the {\it smallest} angle between these lines. If $\ell_1\parallel \ell_2$ (i.e., $\ell_1$ and $\ell_2$ are parallel), we set $\vf(\ell_1,\ell_2)=0$.
\par We let
$$
\RCT=\{\Pi=I_1\times I_2: I_1,I_2\in\Ic(\R)\}
$$
denote the family of all closed rectangles in $\RT$ with sides parallel to the coordinate axes. Finally, by $\HPL$ we denote the family of all closed half-planes in $\RT$, and by $\CRT$ the family of all closed convex subsets of $\RT$.
\msk

\par Given a Banach space $\X$, Przes{\l}awski and Yost \cite{PY2} have introduced an important geometrical characteristic of $\X$, the so-called {\it modulus of squareness} of $\X$. Let us recall its definition.
\par We observe that for any $x,y\in\X$ with $\|y\|<1<\|x\|$ there exists a unique $z=z(x,y)$ with $\|z\|=1$ which belongs to the line segment $[x,y]$. We set
\bel{D-OM-XY}
\omega(x,y)=\frac{\|x-z(x,y)\|}{\|x\|-1}
\ee
and define a function $\xi:[0,1)\to [1,\infty)$ by
\bel{D-XI-B}
\xi_{\X}(\beta)=\sup\,\{\omega(x,y): x,y\in X,\,\|y\|\le \beta<1<\|x\|\,\}.
\ee
We also put
\bel{FI-PS}
\vf(\beta)=(1+\beta)/(1-\beta)~~~~~\text{and}~~~~~ \psi(\beta)=(1-\beta^2)^{-\frac12},~~~~\beta\in[0,1).
\ee
\par It is shown in \cite{PY2}, that for any Banach space $\X$
\bel{F-B}
\xi_X(\beta)\le \vf(\beta)~~~~\text{for every}~~~\beta\in[0,1),
\ee
and
\bel{PS-H}
\xi_X(\beta)=\psi(\beta)~~~~\text{for every}~~~\beta\in[0,1),
\ee
provided $X$ is a Euclidean space.
\par Theorem \reff{PY-M} below recalls an important result from this paper. Recall that $\CX$ denotes the family of all non-empty bounded convex closed subsets of $\X$; let us equip this family with the Hausdorff distance.
\begin{theorem}(\cite[Theorem 4]{PY2})\lbl{PY-M} Let $(S,\delta)$ be a metric space, let $\X$ be a Banach space, and let $f:S\to\X$ and $g:S\to[0,\infty)$ be Lipschitz mappings. Let $F:S\to\CX$ be a Lipschitz (with respect to the Hausdorff distance) set-valued mapping.
\msk
\par Suppose that there exists a constant $\gamma>1$ such that $g(x)\ge\gamma\dist(f(x),F(x))$ for every $x\in S$. Then the intersection mapping $G:S\to\CX$ defined by
$$
G(x)=F(x)\cap \BX(f(x),g(x))
$$
is Lipschitz continuous on $S$ (with respect to $\dhf$) with Lipschitz seminorm
$$
\|G\|_{\Lip(S,\CX)}\le \|F\|_{\Lip(S,\CX)}+
(\,\|F\|_{\Lip(S,\CX)}+\|f\|_{\Lip(S,\X)}+\|g\|_{\Lip(S,\R)})
\,\xi(1/\gamma).
$$
\end{theorem}
\par This theorem enables to prove the following
\begin{proposition}\lbl{HD-IM} Let $X$ be a Banach space, $a\in \X$, $r\ge 0$, and let $C\subset \X$ be a convex set. 
\par Suppose that $C\cap \BX(a,r)\ne\emp$. Then for every $s>0$ and every $L>1$ the following inequality
$$
\dhf\left(C\cap\BX(a,Lr),(C+s\BX)\cap\BX(a,Lr+s)\right)\le
\left(1+2\,\xi_{\X}\left(\tfrac{1}{L}\right)\right)\,s
$$
holds.
\end{proposition}
\par {\it Proof.} Let $S=\{x,y\}\subset\R$ where $x=0$ and $y=s$, and let $\delta(x,y)=s$.
\par We define a mapping $f:S\to X$ and a function $g:S\to\R$ by letting $f(x)=f(y)=a$ and $g(x)=Lr, g(y)=Lr+s$. Clearly, $\|f\|_{\Lip(S,\X)}=0$, and $\|g\|_{\Lip(S,\R)}=1$.
\par We put $\gamma=L$. We note that $C\cap \BX(a,r)\ne\emp$ so that
$$
\dist(f(x),F(x))=\dist(a,C)\le r.
$$
Hence, $g(x)=Lr=\gamma\,r\ge \gamma\,\dist(f(x),F(x))$.
\smsk
\par Then we define a mapping $F:S\to\CX$ by setting  $F(x)=C$ and $F(y)=C+s\BX$. Clearly,
$$
\dhf(F(x),F(y))\le s=\delta(x,y)~~~~\text{so that}~~~~~
\|F\|_{\Lip(S,\CX)}\le 1.
$$
\par Thus, the conditions of Theorem \reff{PY-M} are satisfied for the metric space $(S,\delta)$ and the mappings $f$, $g$ and $F$. This theorem tells us that
the mapping $G:S\to\CX$ defined by
$$
G(u)=F(u)\cap \BX(f(u),g(u)),~~~~~~u\in S,
$$
is Lipschitz on $S$ with respect to the Hausdorff distance. Furthermore,
$$
\|G\|_{\Lip(S,\CX)}\le \|F\|_{\Lip(S,\CX)}+
(\,\|F\|_{\Lip(S,\CX)}+\|f\|_{\Lip(S,\X)}+\|g\|_{\Lip(S,\R)})
\,\xi(1/\gamma)\le 1+2\,\xi_{\X}\left(\tfrac{1}{L}\right).
$$
Hence,
\be
\dhf\left(C\cap\BX(a,Lr),(C+s\BX)\cap\BX(a,Lr+s)\right)
&=&
\dhf(G(x),G(y))\nn\\
&\le& \|G\|_{\Lip(S,\CX)}\,\delta(x,y)
\le
\left(1+2\,\xi_{\X}\left(\tfrac{1}{L}\right)\right)\,s
\nn
\ee
proving the proposition.\bx
\msk
\par Proposition \reff{HD-IM} implies the following important
\begin{proposition}\lbl{N-S} Let $X$ be a Banach space, and let $C\subset \X$ be a convex set. Let $a\in \X$ and let $r\ge 0$. Suppose that 
\bel{C-PR}
C\cap \BX(a,r)\ne\emp.
\ee
\par Then for every $s>0$ and $L>1$
\bel{IN-MN}
C\cap \BX(a,Lr)+\theta(L)\,s\,\BX\supset (C+s\BX)\cap (\BX(a,Lr)+s\BX)
\ee
where
\bel{TH-L}
\theta(L)=(3L+1)/(L-1).
\ee
\par If $X$ is a Euclidean space then \rf{IN-MN} holds with
\bel{TH-L-H}
\theta(L)=1+\frac{2L}{\sqrt{L^2-1}}\,.
\ee
\end{proposition}
\par {\it Proof.} Let
$$
G=C\cap \BX(a,Lr)~~~~\text{and}~~~~
\tG=(C+s\BX)\cap (\BX(a,Lr+s).
$$
\par Definition \rf{HD-R} tells us that $\tG\subset G+\dhf(G,\tG)\,\BX$. In turn, Proposition \reff{HD-IM} states that
$$
\dhf(G,\tG)\le
\left(1+2\,\xi_{\X}\left(\tfrac{1}{L}\right)\right)\,s.
$$
Hence,
$$
\tG\subset G+\Theta(L)\,s\,\BX
~~~~~\text{where}~~~~~
\Theta(L)=1+2\,\xi_{\X}\left(\tfrac{1}{L}\right).
$$
\par Now, let $X$ be an arbitrary Banach space. In this case, thanks to \rf{FI-PS} and \rf{F-B}, we have
$$
\Theta(L)\le 1+2\,\xi_{\X}\left(\tfrac{1}{L}\right)\le 1+
2\frac{1+1/L}{1-1/L}=\frac{3L+1}{L-1}.
$$
This inequality and \rf{TH-L} imply inclusion \rf{IN-MN} in the case under consideration.
\msk
\par Finally, let $X$ be a Euclidean space. In this case, from \rf{FI-PS}, \rf{PS-H} and \rf{TH-L-H}, we have
$$
\Theta(L)=1+2\,\xi_{\X}\left(\tfrac{1}{L}\right)=
1+2\,(1-(1/L)^2)^{-\frac12}=1+\frac{2L}{\sqrt{L^2-1}}
=\theta(L).
$$
\par The proof of the proposition is complete.\bx
\msk
\par For the case  of a Banach space, Proposition \reff{N-S} was proved in \cite[p. 279]{PR}. For similar results we refer the reader to \cite{Ar}, \cite[p. 369]{AF} and \cite[p. 26]{BL}.
\msk
\par For the sake of completeness, and for the reader's convenience, below we give
\msk
\par {\it A direct proof of Proposition \reff{N-S}.} We follow  the proof of Lemma 5.3 from \cite[p. 279]{PR}. If $r=0$ then \rf{IN-MN} holds trivially, so we assume that $r>0$. Without loss of generality, we may also assume that $a=0$. Thus we should prove that
$$
C\cap (Lr\BX)+\theta s\,\BX\supset
(C+s\BX)\cap(Lr\BX+s\BX)
$$
provided $r>0$, $s>0$, $L>1$. Let
\bel{IN-Z}
z\in(C+s\BX)\cap(Lr\BX+s\BX)=(C+s\BX)\cap[(Lr+s)\BX].
\ee
Prove that
\bel{F-Z}
z\in C\cap(Lr\BX)+\theta s\,\BX.
\ee
\par Thanks to \rf{IN-Z}, $z\in(C+s\BX)$ so that there exists an element $v\in C$ such that
\bel{VV-Z}
\|v-z\|\le s.
\ee
\par If $\|v\|\le Lr$, then $v\in C\cap(Lr\BX)$ proving \rf{F-Z}.
\msk
\par Suppose that
\bel{V-Z}
\|v\|>Lr.
\ee
Property \rf{IN-Z} tells us that $\|z\|\le Lr+s$ so that
\bel{V-NR}
\|v\|\le \|z\|+s\le Lr+2s.
\ee
In turn, assumption \rf{C-PR} tells us that there exists an element $v'\in C$ such that
\bel{VP-NR}
\|v'\|\le r<Lr.
\ee
\par Choose $\lambda\in(0,1)$ such that the element
$$
\tv=\lambda v'+(1-\lambda)v
$$
has the norm $\|\tv\|=Lr$. We know that $C$ is convex so that $[v',v]\subset C$ proving that
\bel{VP-D}
\tv\in C\cap (Lr\BX).
\ee
\par Thanks to \rf{VP-NR}, \rf{V-NR} and the triangle inequality,
$$
Lr=\|\tv\|=\|\lambda v'+(1-\lambda)v\|\le \lambda r+(1-\lambda)(Lr+2s)
$$
proving that
$$
\lambda\le \frac{2s}{(L-1)r+2s}.
$$
Consequently, thanks to this inequality, \rf{VP-NR} and \rf{V-NR}
$$
\|v-\tv\|=\lambda\|v-v'\|\le  \lambda(\|v\|+\|v'\|)
\le \frac{2s}{((L-1)r+2s)}\cdot(Lr+2s+r)\le 2s(L+1)/(L-1).
$$
From this inequality and \rf{VV-Z} we have
$$
\|z-\tv\|\le \|z-v\|+\|v-\tv\|\le s+2s(L+1)/(L-1)
=\theta(L)\, s
$$
which together with \rf{VP-D} implies \rf{F-Z}.
\msk
\par Let now $X$ be a Euclidean space. We modify the above proof after \rf{VP-D} as follows.
\par We put $\beta=1/L$, and
\bel{XYW}
x=\mathlarger{\tfrac{1}{L r}}\,v,
~~~~~~y=\mathlarger{\tfrac{1}{L r}}\,v',~~~~~~
w=\mathlarger{\tfrac{1}{L r}}\,\tv.
\ee
\par Then, thanks to \rf{V-Z} and \rf{VP-NR},
\bel{XYB}
\|y\|\le \beta<1<\|x\|.
\ee
\par We note that for any $u,\tu\in X$ such that $\|\tu\|<1<\|u\|$, there exists a unique $w=w(u,\tu)\in[u,\tu]$ with $\|w\|=1$. Hence, thanks to
\rf{D-OM-XY},
\bel{D-OM}
\omega(u,\tu)=\frac{\|u-w(u,\tu)\|}{\|u\|-1}.
\ee
We also recall the definition of the function $\xi_X$, see \rf{D-XI-B}:
\bel{D-XI}
\xi_X(\beta)=\sup\,\{\omega(u,\tu): u,\tu\in X,\,\|\tu\|\le \beta<1<\|u\|\,\}.
\ee
\par Prove that
\bel{XI}
\xi_X(\beta)=(1-\beta^2)^{-\frac12}.
\ee
\par In fact, fix $u$ with $\|u\|>1$. One can easily see that $\sup\,\{\omega(u,\tu): \|\tu\|\le \beta\}$ is attained for some $\tu$ with $\|\tu\|=\beta$, and the line segment $[\tu,u]$ is contained in a line tangent to the sphere with center at the origin and radius $\beta$. Thus, $u-\tu$ is perpendicular to $\tu$. Hence,
$$
\sup\{\omega(u,\tu): \|\tu\|\le \beta\}=f(\|u\|)
$$
where
$$
f(t)=
\frac{\sqrt{t^2-\beta^2}-\sqrt{1-\beta^2}}{t-1},~~~~
t>1.
$$
\par The function $f$ is decreasing on $(1,+\infty)$ so that
$$
\sup_{t>1} f(t)=\lim_{t\to 1} f(t)=(1-\beta^2)^{-\frac12}
$$
proving \rf{XI}.
\msk
\par We apply formula \rf{XI} to the points $x,y,w$ defined
by \rf{XYW}, and to $\beta=1/L$. We have
$$
\frac{\|v-\tv\|}{\|v\|-Lr}=\frac{\|x-w\|}{\|y\|-1}=
\omega(x,y).
$$
See \rf{D-OM}. Thanks to \rf{XYB} and \rf{D-XI},
$$
\frac{\|v-\tv\|}{\|v\|-Lr}\le
\sup\,\{\omega(u,\tu): u,\tu\in X,\,\|\tu\|\le \beta<1<\|u\|\,\}=\xi_X(\beta)=\xi_X(1/L)
$$
so that, thanks to \rf{XI},
$$
\frac{\|v-\tv\|}{\|v\|-Lr}\le
\xi_X(1/L)=\frac{L}{\sqrt{L^2-1}}.
$$
In turn, thanks to \rf{V-NR}, $\|v\|-Lr\le 2s$, so that
$$
\|v-\tv\|\le \frac{L}{\sqrt{L^2-1}}\,(\|v\|-Lr)
\le \frac{2sL}{\sqrt{L^2-1}}.
$$
This inequality and \rf{VV-Z} imply the following:
$$
\|z-\tv\|\le \|z-v\|+\|v-\tv\|\le s+\frac{2sL}{\sqrt{L^2-1}}=
\left(1+2L/\sqrt{L^2-1}\,\right)\,s.
$$
\par This and \rf{VP-D} imply \rf{F-Z} with $\theta=\theta(L)$ defined by \rf{TH-L-H} proving the proposition for a Euclidean space $X$.
\smsk
\par The proof of the proposition is complete.\bx
\msk
\par Proposition \reff{P-F3} below is one of the main ingredients in the proofs of Theorems \reff{MAIN-RT} and \reff{X-LSGM}. The proof of this proposition relies on Proposition \reff{N-S} and Helly's Intersection Theorem for {\it two dimensional} Banach spaces. We recall this theorem below.
\begin{theorem}\lbl{H-TH} Let $\Kc$ be a collection of convex closed subsets of a two dimensional Banach space $\X$. Suppose that $\Kc$ is finite or at least one member of the family $\Kc$ is bounded.
\par If every subfamily of $\Kc$ consisting of at most three elements has a common point then there exists a point common to all of the family $\Kc$.
\end{theorem}
\begin{proposition}\label{P-F3} Let $\X$ be a two dimensional Banach space. Let $C,C_1,C_2\subset \X$ be convex subsets, and let $r>0$. Suppose that
\bel{A-PT}
C_1\cap C_2\cap(C+r\BX)\ne\emp.
\ee
\par Then for every $L>1$ and every $\ve>0$ the following inclusion
\be
&&(C_1\cap C_2+Lr\BX)\cap C+\theta(L)\,\ve\BX
\supset\nn\\
&&[C_1\cap C_2+(Lr+\ve) \BX]
\cap
[(C_1+r\BX)\cap C+\ve \BX]
\cap
[(C_2+r\BX)\cap C+\ve \BX]\nn
\ee
holds. Here $\theta$ is the function from Proposition \reff{N-S}. (Thus, $\theta(L)=(3L+1)/(L-1)$ for an arbitrary $X$, and $\theta(L)=1+2L/\sqrt{L^2-1}$ for a  Euclidean $\X$.)
\end{proposition}
\par {\it Proof.} Suppose that
\bel{A-IN}
a\in
[C_1\cap C_2+(Lr+\ve)\BX]
\cap
[(C_1+r\BX)\cap C+\ve \BX]
\cap
[(C_2+r\BX)\cap C+\ve\BX]
\ee
and prove that
\bel{A-FN1}
a\in
(C_1\cap C_2+Lr\BX)\cap C+\theta(L)\,\ve\BX.
\ee
\par First, let us show that
\bel{N-EM1}
C_1\cap C_2\cap(C+r\BX)\cap \BX(a,Lr+\ve)\ne\emp.
\ee
\par Helly's Theorem \reff{H-TH} tells us that this statement holds provided any three sets in the left hand size of \rf{N-EM1} have a common point.
\par Note that $C_1$, $C_2$ and $C+r\BX$ have a common point. See \rf{A-PT}. We also know that
$$
a\in C_1\cap C_2+(Lr+\ve)\BX,
$$
see \rf{A-IN}, so that
$$
C_1\cap C_2\cap \BX(a,Lr+\ve)\ne\emp\,.
$$
\par Let us prove that
\bel{N-5}
C_1\cap (C+r\BX)\cap \BX(a,Lr+\ve)\ne\emp\,.
\ee
Property \rf{A-IN} tells us that
$$
a\in (C_1+r\BX)\cap C +\ve\BX.
$$
Therefore, there exist elements $b_1\in C_1$ and $b\in C$ such that
$$
\|b_1-b\|\le r~~~\text{and}~~~\|a-b\|\le \ve\,.
$$
In particular, $b_1\in C_1\cap(C+r\BX)$. Furthermore,
$$
\|a-b_1\|\le \|a-b\|+\|b-b_1\|\le \ve +r\le \ve+Lr.
$$
so that $b_1\in \BX(a,Lr+\ve)$.
\par Hence,
$$
b_1\in C_1\cap(C+r\BX)\cap \BX(a,Lr+\ve)
$$
proving \rf{N-5}. In a similar way we show that
$$
C_2\cap (C+r\BX)\cap \BX(a,Lr+\ve)\ne\emp\,.
$$
\par Thus \rf{N-EM1} holds proving the existence of a point $x\in \X$ such that
\bel{N-6}
x\in C_1\cap C_2\cap (C+r\BX)\cap \BX(a,Lr+\ve)\,.
\ee
\par In particular, $x\in C+r\BX$ so that
$\BX(x,r)\cap C\ne\emp$ proving that condition \rf{C-PR} of Proposition \reff{N-S} holds. This proposition tells us that
$$
C\cap \BX(x,Lr)+\theta(L)\,\ve\BX
\supset
(C+\ve \BX)
\cap
(\BX(x,Lr)+\ve \BX)=
(C+\ve \BX)\cap \BX(x,Lr+\ve)\,.
$$
\par From \rf{N-6} and \rf{A-IN} we learn that $a\in \BX(x,Lr+\ve)$ and $a\in C+\ve\BX$. Hence,
$$
(C+\ve\BX)\cap \BX(x,Lr+\ve)\ni a
$$
proving that
$$
C\cap \BX(x,Lr)+\theta(L)\,\ve\BX\ni a.
$$
\par Finally, property \rf{N-6} tells us that $x\in C_1\cap C_2$ proving the required inclusion \rf{A-FN1}.
\par The proof of the proposition is complete.\bx
\bsk
\par We finish the section with the proof of Claim \reff{CLS}, see \cite[Section 5]{FS-2018}. For completeness, we give this simple proof here.
\msk
\par {\it Proof of Claim \reff{CLS}.} The proof relies on the following selection theorem which is a special case of \cite[Theorem 1.2]{S-2004}.
\begin{theorem} Let $\X$ be a Banach space, and let $m\ge 1$. Then there exists a map $\ST:\KM\to \X$ such that
\smallskip
\par $(\alpha)$~~ $\ST(K)\in K$ for all $K\in\KM$
\smallskip
\par\noindent and
\par $(\beta)$~~ $\|\ST(K)-\ST(K')\|\le C(m)\cdot \dhf(K,K')$~~ for all $K,K'\in\KM$.
\smallskip
\par Here $C(m)$ depends only on $m$.
\end{theorem}
\smallskip
\par We refer to $\ST(K)$ as the {\it ``Steiner-type point'' of $K$}, and we call the mapping $\ST:\KM\to \X$ the {\it ``Steiner-type selector''}. In the special case $\X=\R^m$, we can take $\ST(K)$ to be the {\it Steiner point} of $K$, see, e.g., \cite{BL}.
\msk
\par To construct the Lipschitz selection $f$ and establish the claim, we just set
$$
f(x)=\ST(G(x))~~~\text{for}~~~x\in\Mc.
$$
\par Since $G(x)\in\KM$ for each $x\in\Mc$, the function $f$ is well defined on $\Mc$. By part (i) of Definition \reff{D-CORE} and by property $(\alpha)$ of the Steiner-type point,
$$
f(x)=\ST(G(x))\in G(x)\subset F(x)~~~\text{for every}~~~x\in\Mc.
$$
On the other hand, part (ii) of Definition \reff{D-CORE} and property $(\beta)$  of the Steiner-type point imply that
$$
\|f(x)-f(y)\|=\|\ST(G(x))-\ST(G(y))\|\le C(m)\cdot\dhf(G(x),G(y))
\le C(m)\cdot \gamma\,\rho(x,y)
$$
for all $x,y\in\Mc$ proving that $f$ is a Lipschitz selection of the set-valued mapping $F$ with Lipschitz seminorm at most $C(m)\cdot \gamma$.\bx

\SECT{3. Main Theorem for two dimensional Banach spaces.}{3}
\addtocontents{toc}{3. Main Theorem for two dimensional Banach spaces.\hfill \thepage\par\VST}
\indent

\par In this section we prove Theorem \reff{MAIN-RT}
\par First, let us recall the notion of the Lipschitz extension constant $e(\mfM,\X)$ which we use in the formulation of this theorem.
\begin{definition}\lbl{LIP-C} Let $\mfM=(\Mc,\rho)$ be a pseudometric space, and let $X$ be a Banach space.
We define Lipschitz extension constant $e(\mfM,\X)$ of $X$ with respect to $\mfM$ as the infimum of the constants $\lambda>0$ such that for every subset $\Mc'\subset \Mc$, and every Lipschitz mapping $f:\Mc'\to X$, there exists a Lipschitz extension $\tf:\Mc\to X$ of $f$ to all of $\Mc$ such that $\|\tf\|_{\Lip(\Mc,\X)}\le
\lambda \|f\|_{\Lip(\Mc',\X)}$.
\end{definition}
\begin{remark} {\em It is known that
\bel{LNI-E}
e(\mfM,\LNI)=1~~~~\text{for every pseudometric space}~~~ \mfM=\MR.
\ee
(It is immediate from the case $n=1$ which coincides with the McShane-Whitney extension theorem.)
\par It follows from \cite{Rif} and \cite{CL} that
$$
e(\mfM,X)\le 4/3
$$
provided $X$ is an {\it arbitrary two dimensional Banach space}. See also \cite{Bas}.
\msk
\par It is also known that
$e(\mfM,\X)\le n\,\Gamma(\frac{n}{2})/
(\sqrt{\pi}\,\Gamma(\frac{n+1}{2}))$
provided $X$ is an {\it $n$-dimensional Euclidean space}. See \cite{Rif} and \cite{Gr}. In particular,
\bel{E-LTT}
e(\mfM,X)=4/\pi
\ee
whenever $X$ is a {\it two dimensional Euclidean space}.
\par We also note that, by Kirszbraun's extension theorem \cite{Kir}, $e(\mfM,\X)=1$ provided {\it $\X$ is a Euclidean space, $\Mc$ is a subset of a Euclidean space $E$}, and $\rho$ is the metric in $E$.\rbx}
\end{remark}
\smsk
\par {\it Proof of Theorem \reff{MAIN-RT}.} Let $\mfM=(\Mc,\rho)$ be a pseudometric space, and let $\X$ be a two dimensional Banach space. Let $F:\Mc\to \Kc(X)$ be a set-valued mapping satisfying the hypothesis of Theorem \reff{MAIN-RT}. This enables us to make the following
\begin{assumption}\lbl{A-DS} For every subset $\Mc'\subset\Mc$ with $\#\Mc'\le 4$, the restriction $F|_{\Mc'}$ of $F$ to $\Mc'$ has a $\rho$-Lipschitz selection $f_{\Mc'}:\Mc'\to \X$ with $\rho$-Lipschitz seminorm $\|f_{\Mc'}\|_{\Lip((\Mc',\,\rho),\X)}\le 1$.
\end{assumption}
\par Fix constants
\bel{L-GE3}
L\ge 3
\ee
and
\bel{AL-E}
\al\ge e(\mfM,\X)~~~~\text{where}~~~~\mfM=(\Mc,\rho).
\ee
We introduce a new pseudometric on $\Mc$ defined by
\bel{D-AL}
\ds(x,y)=\al\rho(x,y),~~~~~x,y\in \Mc.
\ee
\par This definition, Definition \reff{LIP-C} and inequality \rf{AL-E} imply the following
\begin{claim}\lbl{CL-D} Let $\Mc'\subset \Mc$, and let  $f:\Mc'\to \X$ be a $\rho$-Lipschitz mapping on $\Mc'$. There exists a $\ds$-Lipschitz extension $\tf:\Mc\to X$ of $f$ to all of $\Mc$ with $\ds$-Lipschitz seminorm
$$
\|\tf\|_{\Lip((\Mc,\ds),\X)}\le
\|f\|_{\Lip((\Mc',\,\rho),\X)}.
$$
\end{claim}
\par We introduce set-valued mappings
\bel{F-1G}
F^{[1]}(x)=
\bigcap_{z\in\Mc}\,
\left[F(z)+\ds(x,z)\,\BX\right],~~~~~x\in\Mc,
\ee
and
\bel{F-2G}
F^{[2]}(x)=\bigcap_{z\in\Mc}\,
\left[F^{[1]}(z)+L\ds(x,z)\,\BX\right],~~~x\in\Mc.
\ee
\par Thus, $F^{[1]}$ and $F^{[2]}$ are the first and the second order $(\{1,L\},\ds)$-balanced refinements of $F$ respectively. See Definition \reff{F-IT}.
\msk
\par Formulae \rf{F-1G} and \rf{F-2G} imply an explicit formula for the mapping $F^{[2]}$:
$$
F^{[2]}(x)=\bigcap_{z\in\Mc}\,
\left\{\left(\,\bigcap_{y\in\Mc}\,
\left[F(y)+\ds(z,y)\,\BX\right]\right)
+L\ds(x,z)\,\BX\right\},
~~~x\in\Mc.
$$
\par We will prove that if $L$ and $\al$ are the constants satisfying \rf{L-GE3} and \rf{AL-E} respectively, then the following two facts holds:
\bel{G-NR}
F^{[2]}(x)\ne\emp~~~\text{for every}~~~x\in\Mc,
\ee
and
\bel{H-NR}
\dhf(F^{[2]}(x),F^{[2]}(y))\le \gamma_0(L)\ds(x,y)
~~~\text{for all}~~~x,y\in\Mc.
\ee
Here
\bel{G0}
\gamma_0(L)=L\cdot\theta(L)^2
\ee
where $\theta=\theta(L)$ is the function from Proposition \reff{N-S}. See \rf{TH-L} and \rf{TH-L-H}.
\smsk
\par We prove property \rf{G-NR} and inequality \rf{H-NR} in Proposition \reff{N-EM} and Proposition \reff{HD-G1} respectively.
\par We begin with the proof of property \rf{G-NR}. This proof relies on a series of auxiliary lemmas.
\begin{lemma}\label{H-IN} Let $X$ be a two dimensional Banach space, and let $\Kc$ be a collection of convex closed subsets of $X$ with non-empty intersection. Let $\BS\subset X$ be a convex closed subset symmetric with respect to $0$. Supposed that either (i) at least one member of the family $\Kc$ is bounded or (ii) $\BS$ is bounded. Then
\begin{align}\lbl{H-SM}
\left(\,\bigcap_{K\in\,\Kc} K\right) +\BS
=\bigcap_{K,K'\in\,\Kc}\,
\left\{\,\left(\,K \,\cbig\, K'\right)+\BS\,\right\}.
\end{align}
\end{lemma}
\par {\it Proof.} Obviously, the right hand side of \rf{H-SM} contains its left hand side. Let us prove the converse statement. Fix a point
\begin{align}\lbl{X-IN}
x\in
\bigcap_{K,K'\in\,\Kc}\,
\left\{\,\left(\,K \,\cbig\, K'\right)+\BS\,\right\}
\end{align}
and prove that $x\in \cbg\{K:K\in\,\Kc\} +\BS$.
We know that $\BS$ is symmetric with respect $0$ so that   $-\BS=\BS$. Therefore, $x\in \cbg\{K:K\in\,\Kc\} +\BS$ if and only if
\begin{align}\lbl{BX-1}
B^{(x)}\bigcap\,\left(\,\bigcap_{K\in\,\Kc} K\right)\ne\emp~~~~\text{where}~~~~B^{(x)}=x+\BS.
\end{align}
\par Let $S=\Kc\cbg\,\{B^{(x)}\}$. Helly's intersection Theorem \reff{H-TH} tells us that property \rf{BX-1} holds provided $\cbg\{K:K\in\,S'\}\ne\emp$ for every subfamily $S'\subset S$ consisting of at most three elements.
Clearly, this is true if $B^{(x)}\notin S'$ because there exists a point common to all of the sets from $\Kc$.
\par Suppose that $B^{(x)}\in S'$. Then $S'=\{B^{(x)},K,K'\}$ for some $K,K'\in\Kc$. Then, thanks to \rf{X-IN}, $x\in\left(\,K \cbig K'\right)+B$ proving that
$B^{(x)}\,\cbig\,K\,\cbig\, K'\ne\emp$.
\par Thus, \rf{BX-1} holds, and the proof of the lemma is complete.\bx
\smsk
\begin{lemma}\lbl{G-NE1} For each $x\in\Mc$ the set $F^{[1]}\in\Kc(X)$. Furthermore, for every $x,z\in\Mc$ the following equality
\bel{R-F1}
F^{[1]}(z)+L\ds(x,z)\BX=\bigcap_{y',y''\in\Mc}
\left\{[F(y')+\ds(z,y')\BX]\cap [F(y'')+\ds(z,y'')\BX]+L\ds(x,z)\BX\right\}
\ee
holds.
\end{lemma}
\par {\it Proof.} Let $x\in\Mc$. Formula \rf{F-1G} and Helly's Theorem \reff{H-TH} tell us that $F^{[1]}(x)\ne\emp$ provided
\bel{N-3S}
[F(z_1)+\ds(x,z_1)\BX]\cap [F(z_2)+\ds(x,z_2)\BX]\cap[F(z_3)+\ds(x,z_3)\BX]\ne\emp
\ee
for every $z_1,z_2,z_3\in\Mc$.
\par This property easily follows from Assumption \reff{A-DS}. Indeed, let $\Mc'=\{x, z_1,z_2,z_3\}$. Then $\#\Mc'\le 4$. Assumption \reff{A-DS} tells us that there exists a $\rho$-Lipschitz selection $f_{\Mc'}:\Mc'\to\X$ of $F$ with $\rho$-Lipschitz seminorm $\|f_{\Mc'}\|_{\Lip((\Mc',\,\rho),X)}\le 1$. In particular,
$f_{\Mc'}(z_i)\in F(z_i)$ and
$$
\|f_{\Mc'}(x)-f_{\Mc'}(z_i)\|\le \rho(x,z_i)\le \al\rho(x,z_i)=\ds(x,z_i)~~~\text{for every}~~~i=1,2,3.
$$
See \rf{D-AL}. These properties of $f_{\Mc'}$ and $F$ tell us that the point $f_{\Mc'}(x)$ belongs to the left hand side of \rf{N-3S}. Thus \rf{N-3S} holds for arbitrary $z_i\in\Mc$, $i=1,2,3$, proving that $F^{[1]}(x)\ne\emp$.
\par Finally, this property, Lemma \reff{H-IN} and formula \rf{F-1G} imply \rf{R-F1} proving the lemma.\bx
\begin{lemma}\lbl{REP-G} For every $x\in\Mc$ the following representation
$$
F^{[2]}(x)=\bigcap_{u,u',u''\in\Mc}
\left\{[F(u')+\ds(u',u)\BX]\cap [F(u'')+\ds(u'',u)\BX]+L\ds(u,x)\BX\right\}
$$
holds.
\end{lemma}
\par {\it Proof.} The lemma is immediate from \rf{F-2G} and \rf{R-F1}.\bx
\msk
\msk
\par Given $x,u,u',u''\in\Mc$ we set
\bel{H-D}
T_x(u,u',u'')=
[F(u')+\ds(u',u)\BX]\cap [F(u'')+\ds(u'',u)\BX]+L\ds(u,x)\,\BX
\,.
\ee
In these settings Lemma \reff{REP-G} reformulates as follows:
\bel{G-XP}
F^{[2]}(x)=\bigcap_{u,u',u''\in\Mc} T_x(u,u',u'').
\ee
\begin{proposition}\lbl{N-EM} For every $x\in\Mc$ the set $F^{[2]}(x)\ne\emp$.
\end{proposition}
\par {\it Proof.} Formula \rf{G-XP} and Helly's Theorem \reff{H-TH} tell us that $F^{[2]}(x)\ne\emp$ provided for every choice of points $u_i,u'_i,u''_i\in\Mc$, $i=1,2,3$, we have
\bel{H-INT}
T_x(u_1,u'_1,u''_1)\cap T_x(u_2,u'_2,u''_2)\cap
T_x(u_3,u'_3,u''_3)\ne\emp.
\ee
\par We set
\bel{RI-DF}
r_i=\ds(x,u_i),~~~~i=1,2,3.
\ee
Without loss of generality, we may assume that
\bel{R-123}
r_1\le r_2\le r_3.
\ee
\par For each $i\in\{1,2,3\}$ we also set
\bel{G-UI}
G(u'_i)=F(u'_i)+\ds(u_i',u_i)\BX~~~~\text{and}~~~~
G(u''_i)=F(u''_i)+\ds(u_i'',u_i)\BX.
\ee
\par We will prove that there exist points $y_i\in\X$, $i=1,2,3$, such that
\bel{YK-S}
y_i\in G(u'_i)\cap G(u''_i)~~~~\text{for every}~~~~
i=1,2,3,
\ee
and
\bel{Y-GA}
\|y_1-y_2\|\le r_1+r_2~~~~\text{and}~~~~
\|y_1-y_3\|\le r_1+2r_2+r_3.
\ee
\par Let us see that the existence of the points $y_i$ with these properties implies \rf{H-INT}. In fact, let us set
$$
z=y_1+\tau(y_3-y_1)=y_3+(1-\tau)(y_1-y_3)
$$
with  $\tau=r_1/(r_1+2r_2+r_3)$. Then, thanks to \rf{Y-GA},
$$
\|y_1-z\|=\tau\|y_3-y_1\|\le \frac{r_1}{r_1+2r_2+r_3}\cdot
(r_1+2r_2+r_3)=r_1,
$$
and
$$
\|y_3-z\|=(1-\tau)\|y_3-y_1\|\le \frac{2r_2+r_3}{r_1+2r_2+r_3}\cdot
(r_1+2r_2+r_3)=2r_2+r_3.
$$
Hence,
$$
\|y_2-z\|\le \|y_2-y_1\|+\|y_1-z\|\le r_1+r_2+r_1=2r_1+r_2.
$$
\par From these inequalities, \rf{R-123} and \rf{RI-DF} we have
\bel{Z-YI}
\|z-y_i\|\le 3r_i=3\ds(x,u_i),~~~~i=1,2,3.
\ee
\par Let us prove that $z\in T_x(u_i,u'_i,u''_i)$ for each
$i\in\{1,2,3\}$. In fact, we know that $L\ge 3$, see \rf{L-GE3}. Furthermore, we know that $y_i\in G(u'_i)\cap G(u''_i)$, see \rf{YK-S}, so that, thanks to \rf{Z-YI}, \rf{G-UI} and \rf{H-D},
$$
z\in G(u'_i)\cap G(u''_i)+3\ds(x,u_i)\BX\subset
G(u'_i)\cap G(u''_i)+L\ds(x,u_i)\BX=T_x(u_i,u'_i,u''_i)
$$
proving \rf{H-INT}.
\msk
\par Thus, our aim is to prove the existence of points $y_i$ satisfying \rf{YK-S} and \rf{Y-GA}. We will do this in three steps.
\msk
\par {\it STEP 1}. We introduce sets $W_i\subset\X$, $i=1,...,4$, defined by
\bel{W-123-D}
W_1=G(u'_1),~~~~~W_2=G(u''_1),~~~~W_3=G(u'_2)\cap G(u''_2)+(r_1+r_2)\BX,
\ee
and
\bel{W-4}
W_4=G(u'_3)\cap G(u''_3)+(r_1+2r_2+r_3)\BX.
\ee
\par Obviously, there exist the points $y_i$ satisfying \rf{YK-S} and \rf{Y-GA} whenever
\bel{G-S1}
W_1\cap W_2\cap W_3\cap W_4\ne \emp.
\ee
\par Thus, it suffices to prove property \rf{G-S1}. Helly's Theorem \reff{H-TH} tells us that \rf{G-S1} holds provided
the intersection of any three elements of the family of sets $\{W_1,W_2,W_3,W_4\}$ is non-empty.
\msk
\par {\it STEP 2.} Prove that
\bel{W-1234-P}
W_1\cap W_3\cap W_4\ne\emp.
\ee
This statement is equivalent to the following one:
\bel{W-134}
G(u'_1)\cap [G(u'_2)\cap G(u''_2)+(r_1+r_2)\BX]\cap
[G(u'_3)\cap G(u''_3)+(r_1+2r_2+r_3)\BX]\ne\emp.
\ee
\par Let
\bel{V-123}
V_1=G(u'_1)+(r_1+r_2)\BX,~~~~~V_2=G(u'_2),~~~~V_3=G(u''_2),
\ee
and let
\bel{V-4}
V_4=G(u'_3)\cap G(u''_3)+(r_2+r_3)\BX.
\ee
\par Let us see that if
\bel{G-S2}
V_1\cap V_2\cap V_3\cap V_4\ne \emp
\ee
then \rf{W-1234-P} and \rf{W-134} hold.
\smsk
\par Indeed, definitions \rf{V-123} and \rf{V-4}, and property \rf{G-S2} imply the existence of points $z_1\in G(u'_1)$, $z_2\in G(u'_2)\cap G(u''_2)$ and $z_3\in G(u'_3)\cap G(u''_3)$ such that $\|z_1-z_2\|\le r_1+r_2$ and $\|z_2-z_3\|\le r_2+r_3$. Hence,
$$
\|z_1-z_3\|\le\|z_1-z_2\|+\|z_2-z_3\|\le (r_1+r_2)+(r_2+r_3)=r_1+2r_2+r_3.
$$
\par Thus, thanks to \rf{W-123-D} and \rf{W-4}, the point
$z_1\in W_1\cap W_3\cap W_4$ proving \rf{W-134}.
\msk
\par Let us prove \rf{G-S2}. We will again make use of Helly's Theorem \reff{H-TH} which tells us that \rf{G-S2} holds provided every three elements of the family $\{V_1,V_2,V_3,V_4\}$ have a common point.
\smsk
\par First, let us prove that
\bel{V-124}
V_1\cap V_2\cap V_4=[G(u'_1)+(r_1+r_2)\BX]\cap G(u'_2)\cap[G(u'_3)\cap G(u''_3)+(r_2+r_3)\BX]\ne \emp.
\ee
\par Let $\Mc_1=\{u'_1,u'_2,u'_3,u''_3\}$. Clearly, $\#\Mc_1\le 4$, so that, thanks to Assumption \reff{A-DS}, there exists a $\rho$-Lipschitz mapping $f_{\Mc_1}:\Mc_1\to\X$ with $\|f_{\Mc_1}\|_{\Lip((\Mc_1,\,\rho),\X)}\le 1$ such that
$$
f_{\Mc_1}(u'_1)\in F(u'_1),~~~~
f_{\Mc_1}(u'_2)\in F(u'_2),~~~~
f_{\Mc_1}(u'_3)\in F(u'_3),~~~~\text{and}~~~~
f_{\Mc_1}(u''_3)\in F(u''_3).
$$
\par Claim \reff{CL-D} tells us that there exists a $\ds$-Lipschitz mapping $\tf_1:\Mc\to\X$ with $\ds$-Lipschitz seminorm $\|\tf_1\|_{\Lip((\Mc,\,\ds),\X)}\le
\|f_{\Mc_1}\|_{\Lip((\Mc_1,\,\rho),\X)}\le 1$ such that $\tf_1|_{\Mc_1}=f_{\Mc_1}$.
\smsk
\par Prove that
$$
\tf_1(u_2)\in V_1\cap V_2\cap V_4.
$$
\par We know that
$$
\tf_1(u'_2)=f_{\Mc_1}(u'_2)\in F(u'_2)~~~~\text{and}~~~~
\|\tf_1(u'_2)-\tf_1(u_2)\|\le \ds(u'_2,u_2).
$$
Hence,
$$
\tf_1(u_2)\in F(u'_2)+\ds(u'_2,u_2)\BX=G(u'_2)=V_2.
$$
\par In the same way we prove that $\tf_1(u_1)\in G(u'_1)$. We also know that
$$
\|\tf_1(u_1)-\tf_1(u_2)\|\le\ds(u_1,u_2)
$$
so that $\tf_1(u_2)\in G(u'_1)+\ds(u_1,u_2)\BX.$
By the triangle inequality,
$$
\ds(u_1,u_2)\le \ds(u_1,x)+\ds(x,u_2)=r_1+r_2
$$
proving that $\tf_1(u_2)\in G(u'_1)+(r_1+r_2)\BX=V_1$.
\par We also know that
$$
\tf_1(u'_3)=f_{\Mc_1}(u'_3)\in F(u'_3),~~~~
\tf_1(u''_3)=f_{\Mc_1}(u''_3)\in F(u''_3)
$$
and
$$
\|\tf_1(u_3)-\tf_1(u'_3)\|\le\ds(u_3,u'_3),~~~~
\|\tf_1(u_3)-\tf_1(u''_3)\|\le\ds(u_3,u''_3).
$$
Hence,
$$
\tf_1(u_3)\in [F(u'_3)+\ds(u'_3,u_3)\BX]\cap
[F(u''_3)+\ds(u''_3,u_3)\BX]=G(u'_3)\cap G(u''_3).
$$
Furthermore, $\|\tf_1(u_2)-\tf_1(u_3)\|\le\ds(u_2,u_3)$. These properties of $\tf_1(u_3)$ and the triangle inequality
$\ds(u_2,u_3)\le \ds(u_2,x)+\ds(x,u_3)=r_2+r_3$ imply the following:
$$
\tf_1(u_2)\in G(u'_3)\cap G(u''_3)+\ds(u_2,u_3)\BX
\subset G(u'_3)\cap G(u''_3)+(r_2+r_3)\BX=V_4.
$$
\par Thus, $\tf_1(u_2)\in V_1\cap V_2\cap V_4$ proving \rf{V-124}.
\smsk
\par In the same fashion we show that $V_1\cap V_3\cap V_4\ne\emp$.
\msk
\par Prove that
\bel{V-234}
V_2\cap V_3\cap V_4=G(u'_2)\cap G(u''_2)
\cap[G(u'_3)\cap G(u''_3)+(r_2+r_3)\BX]\ne \emp.
\ee
\par Following the same scheme as in the proof of \rf{V-124}, we introduce a set $\Mc_2=\{u'_2,u''_2,u'_3,u''_3\}$. Assumption \reff{A-DS} provides the existence of a $\rho$-Lipschitz selection $f_{\Mc_2}:\Mc_2\to\X$ of the restriction $F|_{\Mc_2}$ with $\|f_{\Mc_2}\|_{\Lip((\Mc_2,\,\rho),\X)}\le 1$. In turn,
Claim \reff{CL-D} tells us that there exists a $\ds$-Lipschitz mapping $\tf_2:\Mc\to\X$ with $\ds$-Lipschitz seminorm $\|\tf_2\|_{\Lip((\Mc,\,\ds),\X)}\le
\|f_{\Mc_2}\|_{\Lip((\Mc_2,\,\rho),\X)}\le 1$ such that $\tf_2|_{\Mc_1}=f_{\Mc_2}$.
\par Considerations, similar to those in the proof of \rf{V-124}, enable us to show that
$$
\tf_2(u_2)\in V_2\cap V_3\cap V_4
$$
proving \rf{V-234}.
\smsk
\par Finally, following the same ideas, we prove that
\bel{V-F}
V_1\cap V_2\cap V_3=[G(u'_1)+(r_1+r_2)\BX]\cap G(u'_2)\cap G(u''_2)\ne \emp.
\ee
\par We introduce a set $\Mc_3=\{u'_2,u''_2,u'_3,u''_3\}$.
Assumption \reff{A-DS} guarantees the existence of a $\rho$-Lipschitz selection $f_{\Mc_3}:\Mc_3\to\X$ of the restriction $F|_{\Mc_3}$ with $\|f_{\Mc_3}\|_{\Lip((\Mc_3,\,\rho),\X)}\le 1$. Claim \reff{CL-D} provides the existence of a $\ds$-Lipschitz mapping $\tf_3:\Mc\to\X$ with $\ds$-Lipschitz seminorm $\|\tf_3\|_{\Lip((\Mc,\,\ds),\X)}\le
\|f_{\Mc_3}\|_{\Lip((\Mc_3,\,\rho),\X)}\le 1$ such that $\tf_3|_{\Mc_3}=f_{\Mc_3}$.
\par Then we show that $\tf_3(u_2)\in V_1\cap V_2\cap V_3$ proving \rf{V-F}.
\smsk
\par We leave the details of the proofs of properties \rf{V-234} and \rf{V-F} to the interested reader.
\msk
\par Thus, we have proved \rf{G-S2}. As we have seen above, \rf{G-S2} implies \rf{W-1234-P}. See {\it STEP 2.}
\msk
\par {\it STEP 3}. We return to the proof of the property \rf{G-S1}.
\par We prove that
$$
W_2\cap W_3\cap W_4=
G(u''_1)\cap [G(u'_2)\cap G(u''_2)+(r_1+r_2)\BX]\cap
[G(u'_3)\cap G(u''_3)+(r_1+2r_2+r_3)\BX]
\ne\emp
$$
in the same fashion as property \rf{W-1234-P}.
\smsk
\par Let us show that
\bel{W-124}
W_1\cap W_2\cap W_4=
G(u'_1)\cap G(u''_1)\cap
[G(u'_3)\cap G(u''_3)+(r_1+2r_2+r_3)\BX]
\ne\emp.
\ee
\par We put $\Mc_4=\{u'_1,u''_1,u'_3,u''_3\}$. This set contains at most $4$ points. In this case Assumption \reff{A-DS} guarantees us the existence of a $\rho$-Lipschitz mapping $f_{\Mc_4}:\Mc_4\to\X$ with $\|f_{\Mc_4}\|_{\Lip((\Mc_4,\,\rho),\X)}\le 1$ such that $f_{\Mc_4}(u)\in F(u)$ for every $u\in\Mc_4$.
\par Claim \reff{CL-D} enables us to extend $f_{\Mc_4}$ to a $\ds$-Lipschitz mapping $\tf_4:\Mc\to\X$ with $\ds$-Lipschitz seminorm $\|\tf_4\|_{\Lip((\Mc,\,\ds),\X)}\le
\|f_{\Mc_4}\|_{\Lip((\Mc_4,\,\rho),\X)}\le 1$.
Then we show that $\tf_4(u_1)\in W_1\cap W_2\cap W_4$ proving \rf{W-124}.
\msk
\par In a similar way we prove that
\bel{W-123}
W_1\cap W_2\cap W_3=
G(u'_1)\cap G(u''_1)\cap
[G(u'_2)\cap G(u''_2)+(r_1+r_2)\BX]
\ne\emp.
\ee
\par More specifically, we set $\Mc_5=\{u'_1,u''_1,u'_2,u''_2\}$. In this case, Assumption \reff{A-DS} provides the existence of a $\rho$-Lipschitz mapping $f_{\Mc_5}:\Mc_5\to\X$ with $\|f_{\Mc_5}\|_{\Lip((\Mc_5,\,\rho),\X)}\le 1$ such that $f_{\Mc_5}(u)\in F(u)$ for every $u\in\Mc_5$.
\par We again use Claim \reff{CL-D} to extend $f_{\Mc_5}$ to a $\ds$-Lipschitz mapping $\tf_5:\Mc\to\X$ with $\ds$-Lipschitz seminorm $\|\tf_5\|_{\Lip((\Mc,\,\ds),\X)}\le
\|f_{\Mc_5}\|_{\Lip((\Mc_5,\,\rho),\X)}\le 1$.
Finally, we show that
$\tf_5(u_1)\in W_1\cap W_2\cap W_3$ proving \rf{W-123}.
\smsk
\par We leave the details of the proofs of \rf{W-124} and \rf{W-123} to the interested reader.
\msk
\par The proof of the proposition is complete.\bx
\msk
\par We turn to the proof of inequality \rf{H-NR}.
\begin{proposition}\lbl{HD-G1} For every $x,y\in\Mc$ the following inequality
\bel{HD-KJ}
\dhf(F^{[2]}(x),F^{[2]}(y))\le \gamma_0(L)\ds(x,y)
\ee
holds. (We recall that $\gamma_0(L)=L\,\theta(L)^2$ is defined by \rf{G0}.)
\end{proposition}
\par {\it Proof.} Let $x,y\in\Mc$. Formula \rf{G-XP} tells us that
\bel{G-Y}
F^{[2]}(x)=\bigcap_{u,u',u''\in\Mc} T_x(u,u',u'')~~~\text{and}~~~
F^{[2]}(y)=\bigcap_{u,u',u''\in\Mc} T_y(u,u',u'').
\ee
\par Let
$$
\tau=\gamma_0(L)\ds(x,y).
$$
Representation \rf{G-Y}, Lemma \reff{H-IN} and Proposition \reff{N-EM} imply the following:
\bel{GX-1}
F^{[2]}(x)+\tau\,\BX=
\bigcap\,\,
\{\,T_x(u,u',u'')\cap T_x(v,v',v'')
+\tau\,\BX\,\}.
\ee
Here the first intersection in the right hand side of this equality is taken over all
$$
u,u',u'',v,v',v''\in\Mc.
$$
\par Fix $u,u',u'',v,v',v''\in\Mc$ and prove that
\bel{A-1}
A=T_x(u,u',u'')\cap T_x(v,v',v'')
+\tau\,\BX\,\supset F^{[2]}(y)\,.
\ee
\par We introduce the following sets:
\bel{CI-S}
C_1=F(u')+\ds(u',u)\BX,~~~C_2=F(u'')+\ds(u'',u)\BX,~~~
C=T_x(v,v',v'')\,.
\ee
Let
\bel{DL}
\ve=L\,\theta(L)\ds(x,y)~~~~~\text{and}~~~~~r=\ds(x,u)\,.
\ee
\par Then $\tau=\gamma_0(L)\ds(x,y)=\theta(L)\,\ve$,
and
$$
A=T_x(u,u',u'')\cap T_x(v,v',v'')+\tau\,\BX
=(C_1\cap C_2+Lr\BX)\cap C+\theta(L)\,\ve\,\BX.
$$
\par We want to apply Proposition \reff{P-F3} to the set $A$. To do this we have to verify condition \rf{A-PT} of this proposition, i.e., to show that
\bel{C-A}
C_1\cap C_2\cap(C+r\BX)\ne\emp.
\ee
\par Let $\Mc'=\{u',u'',v',v''\}$. This set contains at most four elements so that, thanks to Assumption \reff{A-DS}, there exists a $\rho$-Lipschitz selection $f_{\Mc'}$ of the restriction $F|_{\Mc'}$ with $\|f_{\Mc'}\|_{\Lip((\Mc',\,\rho),\X)}\le 1$.
Claim \reff{CL-D} enables us to extend $f_{\Mc'}$ to a $\ds$-Lipschitz mapping $\tf_1:\Mc\to\X$ defined
on all of $\Mc$ with $\ds$-Lipschitz seminorm
$$
\|\tf_{1}\|_{\Lip((\Mc,\ds),\X)}\le
\|f_{\Mc'}\|_{\Lip((\Mc',\,\rho),\X)}\le 1.
$$
\par In particular, $\tf_1(u')=f_{\Mc'}(u')\in F(u')$, $\tf_1(u'')=f_{\Mc'}(u'')\in F(u'')$,
$$
\|\tf_1(u')-\tf_1(u)\|\le \ds(u',u),~~~
\|\tf_1(u'')-\tf_1(u)\|\le \ds(u'',u)
$$
and
\bel{U-FC}
\|\tf_1(x)-\tf_1(u)\|\le \ds(x,u)=r.
\ee
Hence, $\tf_1(u)\in C_1\cap C_2$.
\par In a similar way we show that $\tf_1(x)\in T_x(v,v',v'')=C$. From this and \rf{U-FC} we have $\tf_1(u)\in C+r\BX$. Hence,
$$
C_1\cap C_2\cap(C+r\BX) \ni \tf_1(u)
$$
proving \rf{C-A}.
\par Thus, property \rf{A-PT} of Proposition \reff{P-F3} holds. This proposition tells us that
\begin{align}
A&=(C_1\cap C_2+Lr\BX)\cap C+\theta(L)\,\ve\,\BX
\nn\\
&\supset [C_1\cap C_2+(Lr+\ve)\BX]
\cap [(C_1+r\BX)\cap C +\ve \BX]
\cap [(C_2+r\BX)\cap C +\ve \BX]
\nn\\
&=S_1\cap S_2\cap S_3.
\nn
\end{align}
\par Prove that
\bel{S-AI}
S_i\supset F^{[2]}(y)~~~\text{for every}~~~i=1,2,3.
\ee
\par We begin with the set
$S_1=C_1\cap C_2+(Lr+\ve)\BX$. Thus,
$$
S_1=
\{F(u')+\ds(u',u)\BX\}\cap \{F(u'')+\ds(u'',u)\BX\}
+(L\ds(u,x)+L\,\theta(L)\ds(x,y))\BX\,.
$$
See \rf{CI-S}. By the triangle inequality,
$$
\ds(u,x)+\theta(L) \ds(x,y)\ge \ds(u,x)+\ds(x,y)\ge \ds(u,y)
$$
so that
$$
S_1\supset
\{F(u')+\ds(u',u)\BX\}\cap \{F(u'')+\ds(u'',u)\BX\}
+L\ds(u,y)\BX=T_y(u,u',u'')\,.
$$
But $T_y(u,u',u'')\supset F^{[2]}(y)$, see \rf{G-Y}, which implies the required inclusion $S_1\supset F^{[2]}(y)$.
\msk
\par We turn to the proof of the inclusion
$S_2\supset F^{[2]}(y)$. Note that $S_2$ is defined by
\bel{DF-S2}
S_2=(C_1+r\BX)\cap C +\ve \BX.
\ee
\par By the triangle inequality,
\bel{N-10}
C_1+r\BX=F(u')+\ds(u',u)\BX+\ds(u,x)\BX\supset
F(u')+\ds(u',x)\BX\,.
\ee
\par Let
\bel{CI-DF}
\tC=F(u')+\ds(u',x)\BX,~~~\tC_1=F(v')+\ds(v',v)\BX,~~~
\tC_2=F(v'')+\ds(v'',v)\BX,
\ee
and let
\bel{TR}
\tr=\ds(v,x).
\ee
\par In these settings
$$
C=T_x(v,v',v'')=\tC_1\cap\tC_2+L\tr \BX.
$$
\par Let
\bel{TC-1}
\tA=(\tC_1\cap\tC_2+L\tr \BX)\cap \tC+\ve\BX.
\ee
Then, thanks to \rf{DF-S2} and \rf{N-10},
\bel{S-TA}
S_2\supset
\{F(u')+\ds(u',x)\BX\}\cap C+\ve \BX=
(\tC_1\cap\tC_2+L\tr \BX)\cap \tC+\ve \BX=\tA\,.
\ee
\par Prove that
\bel{A-GR}
\tA\supset F^{[2]}(y).
\ee
\par As in the previous case, we will do this by applying Proposition \reff{P-F3} to the set $\tA$. But first we have to show that the hypothesis of this proposition holds for $\tA$, i.e.,
\bel{PR-1E}
\tC_1\cap\tC_2\cap(\tC+\tr\BX)\ne\emp\,.
\ee
\par Let $\hM=\{u',v',v''\}$. Assumption \reff{A-DS} tells us that the restriction $F|_{\hM}$ of $F$ to $\hM$ has a $\rho$-Lipschitz selection $f_{\hM}:\hM\to\X$ with $\|f_{\hM}\|_{\Lip((\hM,\,\rho),\X)}\le 1$. In turn,
Claim \reff{CL-D} tells us that there exists a $\ds$-Lipschitz mapping $\tf_2:\Mc\to\X$ with $\ds$-Lipschitz seminorm $\|\tf_2\|_{\Lip((\Mc,\,\ds),\X)}\le
\|f_{\hM}\|_{\Lip((\hM,\,\rho),\X)}\le 1$ such that $\tf_2|_{\hM}=f_{\hM}$.
\par In particular,
$$
\tf_{2}(u')=f_{\hM}(u')\in F(u'),~~\tf_{2}(v')=f_{\hM}(v')\in F(v'), ~~\tf_{2}(v'')=f_{\hM}(v'')\in F(v'').
$$
In addition, $\|\tf_{2}(x)-\tf_{2}(u')\|\le \ds(x,u')$,
$$
\|\tf_{2}(v')-\tf_{2}(v)\|\le \ds(v',v),~~~
\|\tf_{2}(v'')-\tf_{2}(v)\|\le \ds(v'',v)~~~\text{and}~~~
\|\tf_{2}(x)-\tf_{2}(v)\|\le \ds(x,v).
$$
Combining these properties of $\tf_{2}$ with definitions \rf{CI-DF} and \rf{TR}, we conclude that
$$
\tC_1\cap\tC_2\cap(\tC+\tr\BX)\ni \tf_{2}(v)
$$
proving \rf{PR-1E}.
\smsk
\par We recall that $\ve=L\,\theta(L)\ds(x,y)$, see \rf{DL}, so that
$$
\tA=
(\tC_1\cap \tC_2+L\tr \BX)\cap \tC+L\,\theta(L)\ds(x,y)\BX, ~~~~\text{(see \rf{TC-1}).}
$$
We apply Proposition \reff{P-F3} to $\tA$ and obtain the following:
\begin{align}
\tA&\supset \{\tC_1\cap \tC_2+(L\tr+L\ds(x,y))\BX\}\nn\\
&\cap \{(\tC_1+\tr \BX)\cap \tC +L\ds(x,y) \BX\}
\cap \{(\tC_2+\tr \BX)\cap \tC +L\ds(x,y) \BX\}\nn\\
&=\tS_1\cap\tS_2\cap\tS_3.
\nn
\end{align}
\par Prove that
\bel{S-I13}
\tS_i\supset F^{[2]}(y))~~~\text{for every}~~~i=1,2,3.
\ee
\par First, let us show that
\bel{S1-ING}
\tS_1=\tC_1\cap \tC_2+(L\tr+L\,\ds(x,y))\BX\supset F^{[2]}(y).
\ee
By \rf{TR} and the triangle inequality,
$$
\tr+\ds(x,y)=\ds(v,x)+\ds(x,y)\ge \ds(v,y)
$$
so that
\begin{align}
\tS_1&\supset
\tC_1\cap \tC_2+L\ds(v,y)\BX\nn\\
&=
\{F(v')+\ds(v',v)\BX\}\cap
\{F(v'')+\ds(v'',v)\BX\}+L\ds(v,y)\BX\nn\\
&=T_y(v,v',v'').\nn
\end{align}
See \rf{CI-DF} and \rf{H-D}. This inclusion and \rf{G-Y} imply \rf{S1-ING}.
\smsk
\par Prove that
\bel{S2-IN}
\tS_2=(\tC_1+\tr \BX)\cap \tC +L\ds(x,y)\BX\supset F^{[2]}(y)).
\ee
\par Thanks to \rf{CI-DF}, \rf{TR} and the triangle inequality,
$$
\tC_1+\tr \BX=F(v')+\ds(v',v)\BX+\ds(v,x)\BX\supset
F(v')+\ds(v',x)\BX
$$
so that
$$
\tS_2\supset
\{F(v')+\ds(v',x)\BX\}\cap \{F(u')+\ds(u',x)\BX\}
+L\,\ds(x,y)\BX=T_y(x,u',v').
$$
See \rf{H-D}. From this inclusion and \rf{G-Y} it follows that $\tS_2\supset T_y(x,u',v')\supset F^{[2]}(y)$ proving \rf{S2-IN}.
\msk
\par In the same way we prove that
$$
\tS_3=(\tC_2+\tr \BX)\cap \tC +L\ds(x,y) \BX\supset
T_y(x,u',v'')\supset F^{[2]}(y)\,.
$$
\par This inclusion together with \rf{S1-ING} and \rf{S2-IN} imply \rf{S-I13}. Hence,
$$
\tA\supset \tS_1\cap\tS_2\cap\tS_3\supset F^{[2]}(y)
$$
proving \rf{A-GR}.
\par We know that $S_2\supset\tA$, see \rf{S-TA}, so that $S_2\supset F^{[2]}(y)$. In the same fashion we show that
$$
S_3=(C_2+r\BX)\cap C +L\,\ve\,\BX\supset F^{[2]}(y)
$$
proving \rf{S-AI}. Hence,
$$
A\supset S_1\cap S_2\cap S_3\supset F^{[2]}(y)
$$
proving \rf{A-1}.
\smsk
\par Combining \rf{A-1} with \rf{GX-1} we prove that 
$$
F^{[2]}(x)+\gamma_0(L)\ds(x,y)\BX=F^{[2]}(x)+\tau\BX
\supset F^{[2]}(y)\,.
$$
\par By interchanging the roles of $x$ and $y$ we obtain also
$$
F^{[2]}(y)+\gamma_0(L)\ds(x,y)\BX\supset F^{[2]}(x)\,.
$$
These two inclusions imply inequality \rf{HD-KJ} proving the proposition.\bx
\msk
\par We are in a position to finish the proof of Theorem \reff{MAIN-RT}.
\par Let $\lambda_1,\lambda_2$ and $\gamma$ be parameters satisfying \rf{GM-FN}. Thus, $\lambda_1\ge e(\mfM,X)$, $\lambda_2\ge 3\lambda_1$ and
\bel{GM-LL}
\gamma\ge \lambda_2\,(3\lambda_2+\lambda_1)^2/
(\lambda_2-\lambda_1)^2.
\ee
\par We set $\alpha=\lambda_1$, $L=\lambda_2/\lambda_1$. Then $L$ and $\alpha$ satisfies \rf{L-GE3} and \rf{AL-E} respectively, i.e., $L\ge 3$ and $\alpha\ge e(\mfM,X)$.
We also recall that
\bel{DS-L2}
\ds=\alpha\rho=\lambda_1\,\rho,~~~~~~\text{see}~~~~ \rf{D-AL}.
\ee
\par In these settings, the set values mappings $F^{[1]}$ and $F^{[2]}$ defined by formulae \rf{F-1G} and \rf{F-2G} has the following representations:
$$
F^{[1]}(x)=
\bigcap_{z\in\Mc}\,
\left[F(z)+\lambda_1\rho(x,z)\,\BX\right],~~~~~x\in\Mc,
$$
and
$$
F^{[2]}(x)=\bigcap_{z\in\Mc}\,
\left[F^{[1]}(z)+\lambda_2\,\rho(x,z)\,\BX\right],~~~x\in\Mc.
$$
\par In other words, $F^{[1]}$ and $F^{[2]}$ are the first and the second order $(\{\lambda_1,\lambda_2\},\rho)$-balanced refinements of $F$ respectively. See Definition \reff{F-IT}.
\msk
\par Proposition \reff{N-EM} tells us that, under these conditions $F^{[2]}(x)\ne\emp$ for every $x\in\Mc$. In turn, Proposition \reff{HD-G1} states that for every $x,y\in\Mc$ the following inequality
\bel{FF-DH}
\dhf(F^{[2]}(x),F^{[2]}(y))\le \gamma_0(L)\ds(x,y)
\ee
holds. Recall that $\gamma_0(L)=L\cdot\theta(L)^2$
where $\theta=\theta(L)=(3L+1)/(L-1)$, see \rf{TH-L}. Hence,
$$
\theta(L)=\frac{3L+1}{L-1}=\frac{3(\lambda_2/\lambda_1)+1}
{\lambda_2/\lambda_1-1}=\frac{3\lambda_2+\lambda_1}
{\lambda_2-\lambda_1}.
$$
From this,\rf{FF-DH} and \rf{DS-L2}, we have
$$
\dhf(F^{[2]}(x),F^{[2]}(y))\le L\cdot\theta(L)^2\,\ds(x,y))
=\frac{\lambda_2}{\lambda_1}\cdot \frac{(3\lambda_2+\lambda_1)^2}
{(\lambda_2-\lambda_1)^2}\,(\lambda_1\rho(x,y))
=\lambda_2\,\frac{(3\lambda_2+\lambda_1)^2}
{(\lambda_2-\lambda_1)^2}
\,\rho(x,y).
$$
\par This inequality together with \rf{GM-LL} implies the required inequality $\dhf(F^{[2]}(x),F^{[2]}(y))\le \gamma\,\rho(x,y)$ proving Theorem \reff{MAIN-RT}
for $\lambda_1,\lambda_2$ and $\gamma$ satisfying \rf{GM-FN}.
\msk
\par Prove that $F^{[2]}$ satisfies property \rf{F2-NEM} and inequality \rf{HD-RT} for $\lambda_1=4/3$, $\lambda_2=3\lambda_1=4$, and $\gamma=100$. Indeed, we know that $e(\mfM,\X)\le 4/3$ provided $\mfM=\MR$ is an arbitrary pseudometric space, and $\X$ is a {\it two dimensional} Banach space. Therefore, we can set  $\lambda_1=4/3$, $\lambda_2=3\lambda_1=4$, and
$$
\gamma=\lambda_2\,(3\lambda_2+\lambda_1)^2/
(\lambda_2-\lambda_1)^2=100.
$$
\par In these settings, inequalities \rf{GM-FN} for
$\lambda_1,\lambda_2$ and $\gamma$ hold proving \rf{F2-NEM} and \rf{HD-RT} in the case under considerations.
\msk

\par Next, let $X$ be a two dimensional Euclidean space, and let $\lambda_1,\lambda_2$ and $\gamma$ satisfy \rf{X-HSP}. We replace in the above proof the function $\theta=\theta(L)$ defined by \rf{TH-L} with the function $\theta(L)=1+2\,L/\sqrt{L^2-1}$ defined by \rf{TH-L-H}.
\par Following this scheme of the proof, we set $\alpha=\lambda_1$ and $L=\lambda_2/\lambda_1$.
Again, Proposition \reff{N-EM} tells us that $F^{[2]}(x)\ne\emp$ for each $x\in\Mc$. Then we show that inequality \rf{FF-DH} holds for all $x,y\in\Mc$ with $\gamma_0(L)=L\cdot\theta(L)^2$ and $\theta(L)=1+2\,L/\sqrt{L^2-1}$ provided $\lambda_1\ge e(\mfM,X)$ and $\lambda_2\ge 3\lambda_1$.
\par In these settings,
$$
\gamma_0(L)=L\,\theta(L)^2=
(\lambda_2/\lambda_1)\,\left\{1+2\lambda_2\,\slbig
\left(\lambda_2^2-\lambda_1^2\right)^{\frac12}\right\}^2.
$$
From this equality and \rf{FF-DH}, we have
$$
\dhf(F^{[2]}(x),F^{[2]}(y))\le \gamma_0(L)\,(\lambda_1\rho(x,y))
=\lambda_2\left\{1+2\lambda_2\,\slbig
\left(\lambda_2^2-\lambda_1^2\right)^{\frac12}\right\}^2
\rho(x,y)\le \gamma\,\rho(x,y)
$$
provided $\lambda_1, \lambda_2$ and $\gamma$ satisfy \rf{X-HSP}. This proves that inequalities \rf{X-HSP} imply \rf{F2-NEM} and \rf{HD-RT}.
\smsk
\par We know that $e(\mfM,X)=4/\pi$, see \rf{E-LTT}. This and \rf{X-HSP} enable us to set $\lambda_1=4/\pi$, $\lambda_2=12/\pi$. In this case, all three inequalities in \rf{X-HSP} will be satisfied provided
$$
\gamma\ge \lambda_2\left\{1+2\lambda_2\,\slbig
\left(\lambda_2^2-\lambda_1^2\right)^{\frac12}\right\}^2
=3\,(4/\pi)\,\left(1+6/\sqrt{8}\right)^2\approx 37.16\,.
$$
This shows that \rf{F2-NEM} and \rf{HD-RT} hold with $\lambda_1=4/\pi$, $\lambda_2=12/\pi$ and $\gamma=38$.
\smsk
\par Finally, let us assume that $\X$ is a Euclidean space, $\Mc$ is a subset of a Euclidean space $E$, and $\rho$ is the metric in $E$. We know that in this case  $e(\mfM,\X)=1$ (Kirszbraun's extension theorem \cite{Kir}). This enables us to set $\lambda_1=1$ and $\lambda_2=3$. In view of \rf{X-HSP}, for this choice of $\lambda_1$ and $\lambda_2$ one can set
$$
\gamma\ge \lambda_2\left\{1+2\lambda_2\,\slbig
\left(\lambda_2^2-\lambda_1^2\right)^{\frac12}\right\}^2
=3\left(1+6/\sqrt{8}\right)^2\approx 24.99\,.
$$
This proves that \rf{F2-NEM} and \rf{HD-RT} hold with $\lambda_1=1$, $\lambda_2=3$, and $\gamma=25$ provided $X$ is a Euclidean space and $\Mc$ is a subset of a Euclidean space.
\smsk
\par The proof of Theorem \reff{MAIN-RT} is complete.\bx

\SECT{4. Balanced refinements of line segments in a Banach space.}{4}
\addtocontents{toc}{4. Balanced refinements of line segments in a Banach space.\hfill \thepage\par\VST}
\indent

\par In this section we prove Theorem \reff{X-LSGM}. Let $(\Mc,\rho)$ be a pseudometric space, and let $(\X,\|\cdot\|)$ be a Banach space with $\dim X>1$. We recall that $\Kc_1(\X)$ is the family of all non-empty compact convex subsets of $\X$ of dimension at most $1$.
In other words, the family $\Kc_1(\X)=\{[a,b]\subset\X:a,b\in\X\}$ consists of all points and all closed bounded intervals in $\X$.
\par In this section we need the following version of Helly's Theorem.
\begin{theorem}\lbl{HT-IX} Let $\Kc$ be a collection of closed convex subsets of $\X$ containing a set $K_0\in\Kc_1(\X)$. If the intersection of $K_0$ with any two sets from $\Kc$ is non-empty, then there exists a point common to all of the collection $\Kc$.
\end{theorem}
\par {\it Proof.} We introduce a family $\tKc=\{K\cap K_0:K\in\Kc\}$, and apply to $\tKc$ one dimensional Helly's Theorem. (See next section, Lemma \reff{H-R}, part (a).)\bx
\smsk
\par We will also need the following variant of Proposition \reff{P-F3} for the family $\Kc_1(\X)$.
\begin{proposition}\label{C123} Let $\X$ be a Banach space, and let $r\ge 0$. Let $C,C_1,C_2\subset \X$ be convex closed subsets, and let $C_1\in\Kc_1(\X)$. Suppose that
\bel{A-PT-1}
C_1\cap C_2\cap(C+r\BX)\ne\emp.
\ee
\par Then for every $L>1$ and every $\ve>0$ the following inclusion
$$
(C_1\cap C_2+Lr\BX)\cap C+\theta(L)\,\ve\BX
\supset[C_1\cap C_2+(Lr+\ve)\BX]
\cap[(C_1+r\BX)\cap C+\ve\BX]
$$
holds.
\par Here $\theta(L)=(3L+1)/(L-1)$; if $\X$ is a  Euclidean space, one can set $\theta(L)=1+2L/\sqrt{L^2-1}$.
\end{proposition}
\par {\it Proof.} Let
\bel{A-IN1}
a\in
[C_1\cap C_2+(Lr+\ve)\BX]
\cap
[(C_1+r\BX)\cap C+\ve\BX].
\ee
\par Prove that
\bel{A-7}
a\in (C_1\cap C_2+Lr\BX)\cap C+\theta(L)\,\ve\BX\,.
\ee
\par First, let us show that
\bel{N-EM1-1}
C_1\cap C_2\cap(C+r\BX)\cap \BX(a,Lr+\ve)\ne\emp.
\ee
\par Recall that $C_1\in\Kc_1(\X)$. Helly's Theorem \reff{HT-IX} tells us that it is suffices to show that any two sets in the left hand size of \rf{N-EM1-1} have a common point with $C_1$.
\par First we note that $C_1$, $C_2$ and $C+r\BX$ have a common point. See \rf{A-PT-1}. We also know that
$$
a\in C_1\cap C_2+(Lr+\ve)\BX,
$$
see \rf{A-IN1}, so that $C_1\cap C_2\cap \BX(a,Lr+\ve)\ne\emp$.
\par Let us prove that
\bel{N-51}
C_1\cap (C+r\BX)\cap \BX(a,2r+\ve)\ne\emp\,.
\ee
Property \rf{A-IN1} tells us that
$a\in (C_1+r\BX)\cap C +\ve\BX$.
Therefore, there exist  points $b_1\in C_1$ and $b\in C$ such that $\|b_1-b\|\le r$ and $\|a-b\|\le \ve$.
In particular, $b_1\in C_1\cap(C+r\BX)$. Furthermore,
$$
\|a-b_1\|\le \|a-b\|+\|b-b_1\|\le \ve +r\le \ve +2r,
$$
so that $b_1\in \BX(a,2r+\ve)$. Hence,
$$
b_1\in C_1\cap(C+r\BX)\cap \BX(a,2r+\ve)
$$
proving \rf{N-51}.
\smsk
\par Thus, \rf{N-EM1-1} holds proving the existence of a point $x\in \X$ such that
\bel{N-61}
x\in C_1\cap C_2\cap (C+r\BX)\cap \BX(a,2r+\ve)\,.
\ee
\par In particular, $x\in C+r\BX$ so that
$\BX(x,r)\cap C\ne\emp$ proving that condition \rf{C-PR} of Proposition \reff{N-S} is satisfied. We apply this proposition to $x$, $r$ and the set $C$ and get:
$$
C\cap \BX(x,Lr)+\theta(L)\,\ve\BX\supset
(C+\ve\BX)\cap
(\BX(x,Lr)+\ve \BX)=
(C+\ve\BX)\cap \BX(x,Lr+\ve)\,.\nn
$$
\par From \rf{N-61} we learn that $a\in \BX(x,Lr+\ve)$. In turn, \rf{A-IN1} tells us that
$$
a\in (C_1+r\BX)\cap C+\ve\BX\subset
C+\ve\BX.
$$
Hence, $(C+\ve\BX)\cap \BX(x,Lr+\ve)\ni a$ proving that
$C\cap \BX(x,Lr)+\theta(L)\,\ve\BX\ni a$.
\smsk
\par Finally, property \rf{N-61} tells us that $x\in C_1\cap C_2$ proving the required inclusion \rf{A-7}.\bx
\msk
\par We note that the finiteness number $N(1,\X)=\min\{2^{2},2^{\dim \X}\}=4$. (Recall that $\dim\X>1$.) Let $F:\Mc\to \Kc_1(\X)$ be a set-valued mapping. We suppose that $F$ satisfies the hypothesis of Theorem \reff{X-LSGM}, i.e., that the following assumption is true.
\smsk
\begin{assumption}\lbl{A-XK1} For every subset $\Mc'\subset\Mc$ with $\#\Mc'\le 4$ the restriction $F|_{\Mc'}$ of $F$ to $\Mc'$ has a Lipschitz selection $f_{\Mc'}:\Mc'\to \X$ with $\|f\|_{\Lip(\Mc',\X)}\le 1$.
\end{assumption}
\par Let $\vl=\{\lambda_1,\lambda_2\}$. We introduce balanced $(\vl,\rho)$-refinements of $F$ of the first and the second order, i.e., set-valued mappings
$$
F^{[1]}(x)=
\bigcap_{y\in\Mc}\,
\left[F(y)+\lambda_1\rho(x,y)\,\BX\right],~~~x\in\Mc,
$$
and
$$
F^{[2]}(x)=
\bigcap_{z\in\Mc}\,
\left[F^{[1]}(z)+\lambda_2\rho(x,z)\,\BX\right],~~~x\in\Mc.
$$
See Definition \reff{F-IT}.
\msk
\par Our aim is to prove that if
\bel{GM-FN-P}
\lambda_1\ge 1, ~~~~~~\lambda_2\ge 3\lambda_1,
~~~~~~
\gamma\ge \lambda_2\,(3\lambda_2+\lambda_1)/
(\lambda_2-\lambda_1),
\ee
then the set-valued mapping $F^{[2]}$ is a $\gamma$-core of $F$ (with respect to $\rho$), i.e.,
$$
F^{[2]}(x)\ne\emp~~\text{for every}~~x\in\Mc,~~~\text{and}~~~
\dhf(F^{[2]}(x),F^{[2]}(y))\le \gamma\rho(x,y)~~\text{for all}~~x,y\in\Mc.
$$
\par We set $L=\lambda_2/\lambda_1$. We also introduce a new pseudometric on $\Mc$ defined by
$$
\ds(x,y)=\lambda_1\rho(x,y),~~~~~x,y\in \Mc.
$$
Thus, thanks to \rf{GM-FN-P},
\bel{L-DS2}
L\ge 3~~~~~~\text{and}~~~~~~ \rho\le \ds ~~~~\text{on}~~~~\Mc.
\ee
\par In these settings,
\bel{F12-X1}
F^{[1]}(x)=
\bigcap_{z\in\Mc}\,
\left[F(z)+\ds(x,z)\BX\right]
~~~\text{and}~~~
F^{[2]}(x)=\bigcap_{z\in\Mc}\,
\left[F^{[1]}(z)+L\ds(x,z)\BX\right],
~~~x\in\Mc.
\ee
\par We need the following analog of Lemma \reff{H-IN}.
\begin{lemma}\lbl{H-1} Let $\Kc$ be a collection of convex closed subsets of $\X$ containing a set $K_0\in\Kc_1(\X)$. Suppose that $\cap\{K:K\in\,\Kc\}\ne\emp$. Then for every $r\ge 0$ the following equality
$$
\left(\,\bigcap_{K\in\,\Kc} K\right) +r\BX
=\bigcap_{K\in\,\Kc}\,
\left\{\,\left[\,K \,\cbig\, K_0\right]+r\BX\,\right\}
$$
holds.
\end{lemma}
\par {\it Proof.} Let $\tKc=\{K\cap K_0:K\in\Kc\}$. Clearly, $\tKc\subset \Kc_1(\X)$. It is also clear that the statement of the lemma is equivalent to the equality
$$
\left(\,\bigcap_{\tK\in\,\tKc} \tK\right) +r\BX
=\bigcap_{\tK\in\,\tKc}\,
\left\{\,\tK +r\BX\,\right\}
$$
provided $\cap\{\tK:\tK\in\Kc\}\ne\emp$. We prove this equality by a slight modification of the proof of Lemma \reff{H-IN}. In particular, in this proof we use Helly's Theorem \reff{HT-IX} rather than Theorem \reff{H-TH}.
We leave the details to the interested reader.\bx
\smsk
\par The next lemma is an analog of Lemma \reff{G-NE1}.
\begin{lemma}\lbl{FZ-E} For every $x\in\Mc$ the set $F^{[1]}(x)\in\Kc_1(X)$. Moreover, for every $x,z\in\Mc$ we have
\bel{FD-2}
F^{[1]}(z)+L\ds(x,z)\BX=\bigcap_{v\in\Mc}
\left\{[F(v)+\ds(z,v)\BX]\cap F(z)+L\ds(x,z)\BX\right\}.
\ee
\end{lemma}
\par {\it Proof.} Let $\Kc=\{F(z)+\ds(z,x)\BX:z\in\Mc\}$. We know that $\Kc$ is a family of bounded closed convex subsets of $\X$ containing the set $F(x)\in\Kc_1(\X)$. Helly's Theorem \reff{HT-IX} tells us that the set $F^{[1]}(x)=\cap\{K:K\in\Kc\}\ne \emp$ provided the set
\bel{W-1}
E=F(x)\cap [F(z')+\ds(z',x)\BX]\cap [F(z'')+\ds(z'',x)\BX]
\ne\emp
\ee
for every $z',z''\in\Mc$.
\par Let $\Mc'=\{x,z',z''\}$. Assumption \reff{A-XK1} tells us that there exists a function $f_{\Mc'}:\Mc'\to\X$ satisfying the following conditions: $f_{\Mc'}(x)\in F(x)$, $f_{\Mc'}(z')\in F(z')$, $f_{\Mc'}(z'')\in F(z'')$,
$$
\|f_{\Mc'}(z')-f_{\Mc'}(x)\|\le \rho(z',x)\le \ds(z',x),~~~\text{and}
~~~\|f_{\Mc'}(z'')-f_{\Mc'}(x)\|\le \rho(z'',x)\le \ds(z'',x).
$$
See \rf{L-DS2}. Hence, $f_{\Mc'}(x)\in E$ proving \rf{W-1}. Thus, $F^{[1]}(x)\ne \emp$.
\smsk
\par We also know that $F^{[1]}(x)\in\Kc_1(\X)$. These properties of $F^{[1]}(x)$, Lemma \reff{H-1} and \rf{F12-X1} imply the required equality \rf{FD-2} proving the lemma.\bx
\msk
\par Note that, thanks to \rf{F12-X1}, the following explicit representation of the mapping $F^{[2]}$
$$
F^{[2]}(x)=\bigcap_{z\in\Mc}\,
\left\{\left(\,\bigcap_{y\in\Mc}\,
\left[F(y)+\ds(z,y)\,\BX\right]\right)
+L\ds(x,z)\,\BX\right\},
~~~x\in\Mc,
$$
holds. This representation and Lemma \reff{FZ-E} imply the following analog of Lemma \reff{REP-G}.
\begin{lemma}\lbl{GW-2} For every $x\in\Mc$ the following equality
$$
F^{[2]}(x)=\bigcap_{u,u'\in\Mc}
\left\{[F(u')+\ds(u',u)\BX]\cap F(u)+L\ds(u,x)\BX\right\}
$$
holds.
\end{lemma}
\par Given $x,u,u'\in\Mc$ we put
\bel{TH-D}
\tT_x(u,u')=
[F(u')+\ds(u',u)\BX]\cap F(u)+L\ds(u,x)\,\BX.
\ee
\par Now, Lemma \reff{GW-2} reformulates as follows:
\bel{G-XP1}
F^{[2]}(x)=\bigcap_{u,u'\in\Mc} \tT_x(u,u').
\ee
\begin{proposition}\lbl{TN-EM} For every $x\in\Mc$ the set $F^{[2]}(x)\ne\emp$.
\end{proposition}
\par {\it Proof.} Clearly, $F(x)=\tT_x(x,x)$, see \rf{H-D}. We also know that $F(x)\in\Kc_1(\X)$. Formula \rf{G-XP1} and Helly's Theorem \reff{HT-IX} tell us that $F^{[2]}(x)\ne\emp$ provided for every choice of points $u_i,u'_i\in\Mc$, $i=1,2$, we have
\bel{TH-INT}
F(x)\cap\tT_x(u_1,u'_1)\cap \tT_x(u_2,u'_2)\ne\emp.
\ee
\par We recall that
\bel{H-DLP}
\tT_x(u_i,u'_i)=
[F(u_i')+\ds(u_i',u_i)\BX]\cap F(u_i)+L\ds(u_i,x)\,\BX,
~~~i=1,2.
\ee
See \rf{TH-D}. Without loss of generality, we may assume that
\bel{A-R}
\rho(u_1,x)\ge \rho(u_2,x).
\ee
\par We introduce the following sets:
\bel{G-S123}
G_1=F(u_2),~~~~
G_2=F(u_2')+\rho(u_2,u_2')\BX,~~~~
G_3=F(x)+\rho(u_2,x)\BX,
\ee
and
\bel{G-S4}
G_4=[F(u_1')+\rho(u_1',u_1)\BX]\cap F(u_1)+\rho(u_1,u_2,)\BX\,.
\ee
\par Prove that
\bel{TG-NE}
\bigcap_{i=1}^4\,G_i\ne\emp.
\ee
We know that $G_1=F(u_2)\in\Kc_1(X)$. In this case, Helly's Theorem \reff{HT-IX} tells us that \rf{TG-NE} holds provided $G_1\cap G_i\cap G_j\ne\emp$ for every $2\le i,j\le 4$, $i\ne j$.
\smsk
\par First prove that
\bel{I-S123}
G_1\cap G_2\cap G_3=F(u_2)\cap [F(u_2')+\rho(u_2,u_2')\BX]\cap
[F(x)+\rho(u_2,x)\BX]\ne\emp.
\ee
\par Let $\Mc_1=\{u_2',u_2,x\}$. Because $\#\Mc_1\le 4$, Assumption \reff{A-XK1} guarantees the existence of a mapping $f_1:\Mc_1\to\X$ with the following properties:  $f_1(x)\in F(x)$, $f_1(u_2)\in F(u_2)$, $f_1(u_2')\in F(u_2')$,
$$
\|f_1(u_2)-f_1(x)\|\le \rho(u_2,x)~~~~\text{and}~~~~
\|f_1(u_2)-f_1(u_2')\|\le \rho(u_2,u_2').
$$
These properties of $f_1$ and definition \rf{G-S123} tell us that $f_1(u_2)\in G_1\cap G_2\cap G_3$ proving \rf{I-S123}.
\smsk
\par Prove that
$$
G_1\cap G_2\cap G_4=F(u_2)\cap [F(u_2')+\rho(u_2,u_2')\BX]\cap
\{[F(u_1')+\rho(u_1',u_1)\BX]\cap F(u_1)+\rho(u_1,u_2)\BX\}\ne\emp.
$$
\par Let $\Mc_2=\{u_1',u_1,u_2',u_2\}$. Clearly, $\#\Mc_2\le 4$. Assumption \reff{A-XK1} tells us that there  exists a mapping $f_2:\Mc_2\to\X$ with the following properties:  $f_2(u_i)\in F(u_i)$, $f_2(u_i')\in F(u_i')$,~
$i=1,2$,
$$
\|f_2(u_1)-f_2(u_1')\|\le \rho(u_1,u_1'),~~~~
\|f_2(u_1)-f_2(u_2)\|\le \rho(u_1,u_2),~~~~
\text{and}~~~~
\|f_2(u_2)-f_2(u_2')\|\le \rho(u_2,u_2').
$$
From these properties of $f_2$ and definitions \rf{G-S123} and \rf{G-S4}, we have $f_2(u_2)\in G_1\cap G_2\cap G_4$ proving that $G_1\cap G_2\cap G_4\ne\emp$.
\smsk
\par Finally, prove that
$$
G_1\cap G_3\cap G_4=F(u_2)\cap [F(x)+\rho(u_2,x)\BX]
\cap
\{[F(u_1')+\rho(u_1',u_1)\BX]\cap F(u_1)+\rho(u_1,u_2)\BX\}\ne\emp.
$$
\par We introduce a set $\Mc_3=\{u_1',u_1,x,u_2\}$. Because $\#\Mc_2\le 4$, by Assumption \reff{A-XK1}, there exists a mapping $f_3:\Mc_3\to\X$ with the following properties:  $f_3(v)\in F(v)$ for each $v\in\Mc_3$,
$$
\|f_3(u_1)-f_3(u_1')\|\le \rho(u_1,u_1'),~~~~
\|f_3(u_1)-f_3(u_2)\|\le \rho(u_1,u_2),~~~~
\text{and}~~~~
\|f_3(u_2)-f_3(x)\|\le \rho(u_2,x).
$$
These properties of $f_2$, \rf{G-S123} and \rf{G-S4} tell us that $f_3(u_2)\in G_1\cap G_3\cap G_4$ proving the required property $G_1\cap G_3\cap G_4\ne\emp$.
\smsk
\par Thus, property \rf{TG-NE} is proven. Let $\tMc=\{u_1',u_1,x,u_2,u_2'\}$. Property \rf{TG-NE} and definitions \rf{G-S123}, \rf{G-S4} imply the existence of a mapping $g:\tMc\to X$ with the following properties: $g(v)\in F(v)$ for every $v\in\tMc$,
\bel{G-LP}
\|g(u_1)-g(u_1')\|\le \rho(u_1,u_1'),~~~~
\|g(u_1)-g(u_2)\|\le \rho(u_1,u_2),~~~~
\|g(u_2)-g(u_2')\|\le \rho(u_2,u_2'),
\ee
and
\bel{G-U2X}
\|g(u_2)-g(x)\|\le \rho(u_2,x).
\ee
\par Prove \rf{TH-INT} by showing that
\bel{GX-TH}
g(x)\in F(x)\cap\tT_x(u_1,u'_1)\cap \tT_x(u_2,u'_2).
\ee
\par Indeed, from properties of $g$ we know that $g(x)\in F(x)$.
\par We also know that $g(u_2)\in F(u_2)$, $g(u_2')\in F(u_2')$. Thanks to \rf{G-LP}, \rf{G-U2X} and \rf{L-DS2},
$$
\|g(u_2)-g(u_2')\|\le \rho(u_2,u_2')\le \ds(u_2,u_2')
~~~~\text{and}~~~~
\|g(u_2)-g(x)\|\le \rho(u_2,x)\le L\ds(u_2,x).
$$
\par From these properties of $g$ and definition \rf{H-DLP}, we have
$$
g(x)\in
[F(u_2')+\ds(u_2',u_2)\BX]\cap F(u_2)+L\ds(u_2,x)\,\BX=\tT_x(u_2,u'_2).
$$
\par It remains to show that $g(x)\in \tT_x(u_1,u'_1)$. Indeed, as we know,
\bel{GU-IN}
g(x)\in F(x),~~~~~g(u_1)\in F(u_1),~~~~~\text{and}~~~~~ g(u_1')\in F(u_1').
\ee
Furthermore, thanks to \rf{G-LP} and \rf{L-DS2},
\bel{GU1}
\|g(u_1)-g(u_1')\|\le \rho(u_1,u_1')\le \ds(u_1,u_1').
\ee
\par Let us estimate $\|g(u_1)-g(x)\|$. Thanks to \rf{G-LP}, \rf{G-U2X} and the triangle inequality,
\be
\|g(u_1)-g(x)\|&\le& \|g(u_1)-g(u_2)\|+\|g(u_2)-g(x)\|\le \rho(u_1,u_2)+\rho (u_2,x)\nn\\
&\le& (\rho(u_1,x)+\rho(x,u_2)+\rho (u_2,x)=
\rho(u_1,x)+2\rho(x,u_2)\nn
\ee
so that, thanks to \rf{A-R} and \rf{L-DS2},
$$
\|g(u_1)-g(x)\|\le 3\rho(u_1,x)\le L\ds(u_1,x).
$$
\par From this inequality, \rf{GU1}, property \rf{GU-IN} and definition \rf{H-DLP} we obtain the required property $g(x)\in \tT_x(u_1,u'_1)$ proving \rf{GX-TH}.
\smsk
\par The proof of the proposition is complete.\bx
\msk
\par As in Section 3, we again set
$$
\gamma_0=\gamma_0(L)=L\,\theta(L)^2
$$
where $\theta=\theta(L)$ is the function from Proposition \reff{C123}. (Thus $\theta(L)=(3L+1)/(L-1)$; if $\X$ is a  Euclidean space, one can set $\theta(L)=1+2L/\sqrt{L^2-1}$.) Cf. \rf{G0}.
\begin{proposition}\lbl{HD-1D} For every $x,y\in\Mc$ the following inequality
\bel{HD-1}
\dhf(F^{[2]}(x),F^{[2]}(y))\le \gamma_0(L)\ds(x,y)
\ee
holds.
\end{proposition}
\par {\it Proof.} Let $x,y\in\Mc$. Formula \rf{G-XP1} tells us that
\bel{GXY-2}
F^{[2]}(x)=\bigcap_{u,u'\in\Mc} \tT_x(u,u')~~~\text{and}~~~
F^{[2]}(y)=\bigcap_{u,u'\in\Mc} \tT_y(u,u').
\ee
Recall that
\bel{H-RM}
\tT_x(u,u')=
[F(u')+\ds(u',u)\BX]\cap F(u)+L\ds(u,x)\,\BX.
\ee
\par We also know that the set $F^{[2]}(x)\ne\emp$, see Proposition \reff{TN-EM}, and the set $\tT_x(x,x)=F(x)\in\Kc_1(\X)$. These properties, the above formula for $F^{[2]}(x)$ and Lemma \reff{H-1} tell us that
\bel{THXY}
F^{[2]}(x)+\gamma_0(L)\ds(x,y)\,\BX=
\bigcap_{u,u'\in\Mc}\,\,
\left\{\,\tT_x(u,u')\cap F(x)+\gamma_0(L)\ds(x,y)\,\BX\,
\right\}.
\ee
\par We fix $u,u'\in\Mc$ and introduce a set
$$
\tA=\tT_x(u,u')\cap F(x)+\gamma_0(L)\ds(x,y)\,\BX.
$$
We also introduce sets
\bel{C-K1}
C_1=F(u),~~~C_2=F(u')+\ds(u',u)\BX,~~~~\text{and}~~~~C=F(x).
\ee
Let
\bel{DL1}
\ve=L\ds(x,y)~~~\text{and}~~~r=\ds(x,u)\,.
\ee
\par In these settings, $\gamma_0(L)\ds(x,y)=\theta(L)\,\ve$
and
$$
\tA=\tT_x(u,u')\cap F(x)+\gamma_0(L)\ds(x,y)\,\BX
=(C_1\cap C_2+Lr\BX)\cap C+\theta(L)\,\ve\,\BX.
$$
\par Let us apply Proposition \reff{C123} to the set $\tA$. First, we have to verify condition \rf{A-PT-1} of this proposition, i.e., to show that
\bel{C-INT}
C_1\cap C_2\cap(C+r\BX)\ne\emp.
\ee
\par Let $\tMc=\{x,u,u'\}$. This set consists of at most three points so that, thanks to Assumption \reff{A-XK1}, there exists a $\rho$-Lipschitz selection $f_{\tMc}$ of the restriction $F|_{\tMc}$ with $\|f_{\tMc}\|_{\Lip((\tMc;\rho),\X)}\le 1$. Thus, $f_{\tMc}(u')\in F(u')$, $f_{\tMc}(u)\in F(u)$, $f_{\tMc}(x)\in F(x)$,
$$
\|f_{\tMc}(u')-f_{\tMc}(u)\|\le \rho(u',u)~~~
\text{and}~~~
\|f_{\tMc}(x)-f_{\tMc}(u)\|\le \rho(x,u).
$$
\par Let us see that
\bel{FM-U}
f_{\tMc}(u)\in C_1\cap C_2\cap(C+r\BX).
\ee
\par Indeed, we know that $f_{\tMc}(u)\in F(u)=C_1$, see \rf{C-K1}. Furthermore, $f_{\tMc}(u')\in F(u')$ and, thanks to \rf{L-DS2}, $\rho\le \ds$ on $\Mc$. Hence,
$$
\|f_{\tMc}(u')-f_{\tMc}(u)\|\le \rho(u',u)\le \ds(u',u)
$$
proving that $f_{\tMc}(u)\in C_2$, see \rf{C-K1}.
\par Finally, thanks to \rf{C-K1} and \rf{DL1}, $f_{\tMc}(x)\in F(x)=C$ and
$$
\|f_{\tMc}(x)-f_{\tMc}(u)\|\le \rho(x,u)\le \ds(x,u)=r,
$$
proving that $f_{\tMc}(u)\in C+r\BX$.
\par Thus, \rf{FM-U} is true, and property \rf{C-INT} holds. We also recall that the set $C_1=F(u)\in\Kc_1(\X)$. Now, Proposition \reff{C123} tells us that
\begin{align}
\tA&=(C_1\cap C_2+Lr\BX)\cap C+\theta(L)\,\ve\,\BX
\nn\\
&\supset [C_1\cap C_2+(Lr+\ve)\BX]
\cap [(C_1+r\BX)\cap C +\ve\BX]
\nn\\
&=\tS_1\cap \tS_2.
\nn
\end{align}
\par Prove that $\tS_i\supset F^{[2]}(y)$ for every $i=1,2$. We begin with the set $\tS_1=C_1\cap C_2+(Lr+\ve)\BX$. Thanks to \rf{C-K1} and \rf{DL1},
$$
\tS_1=
\{F(u')+\ds(u',u)\BX\}\cap F(u)
+(L\ds(u,x)+L\ds(x,y))\BX\,.
$$
By the triangle inequality,
$\rho(u,x)+\rho(x,y)\ge \rho(u,y)$ so that
$$
\tS_1\supset
[F(u')+\ds(u',u)\BX]\cap F(u)
+L\ds(u,y)\BX=\tT_y(u,u'),~~~~~\text{see \rf{H-RM}}.
$$
But, thanks to \rf{GXY-2},  $\tT_y(u,u')\supset F^{[2]}(y)$ which implies the required inclusion $\tS_1\supset F^{[2]}(y)$.
\smsk
\par We turn to the set $\tS_2=
(C_1+r\BX)\cap C +\ve\BX$. Definitions \rf{TH-D}, \rf{C-K1} and  \rf{DL1} tell us that
$$
\tS_2=[F(u)+\ds(u,x)\BX]\cap F(x)+L\ds(x,y)\BX=T_y(u,x).
$$
Thanks to \rf{GXY-2}, $\tT_y(u,x)\supset F^{[2]}(y)$ proving that $\tS_2\supset F^{[2]}(y)$.
\smsk
\par Thus,
$$
\tA=\tT_x(u,u')\cap F(x)+\gamma_0(L)\ds(x,y)\,\BX
\supset \tS_1\cap \tS_2\supset F^{[2]}(y)~~~~\text{for every}~~~~u,u'\in\Mc.
$$
From this and representation \rf{THXY}, we have
$$
F^{[2]}(x)+\gamma_0(L)\ds(x,y)\,\BX\supset F^{[2]}(y)\,.
$$
\par By interchanging the roles of $x$ and $y$ we obtain also
$$
F^{[2]}(y)+\gamma_0(L)\ds(x,y)\,\BX\supset F^{[2]}(x)\,.
$$
\par These two inclusions imply the required inequality \rf{HD-1} proving the proposition.\bx
\msk
\par We finish the proof of Theorem \reff{X-LSGM} as follows. We fix $\lambda_1$, $\lambda_2$ and $\gamma$ satisfying inequalities \rf{GM-FN-1}. Proposition \reff{TN-EM} tells us that for such choice of these parameters the set $F^{[2]}(x)\ne\emp$ for every $x\in\Mc$.
\par In turn, Proposition \reff{HD-1D} tells us that in these settings  $\dhf(F^{[2]}(x),F^{[2]}(y))\le \gamma_0(L)\ds(x,y)$ for all $x,y\in\Mc$. We recall that here $L=\lambda_2/\lambda_1$, $\ds=\lambda_1\rho$, $\gamma_0(L)=L\,\theta(L)$ and
$\theta(L)=(3L+1)/(L-1)$. Hence,
\be
\dhf(F^{[2]}(x),F^{[2]}(y))&\le&
\gamma_0(L)\ds(x,y)=L\left(\frac{3L+1}{L-1}\right)\ds(x,y)
\nn\\
&=&
(\lambda_2/\lambda_1)
\frac{3(\lambda_2/\lambda_1)+1}{(\lambda_2/\lambda_1)-1}
\cdot (\lambda_1\rho(x,y))
=\{\lambda_2\,(3\lambda_2+\lambda_1)/
(\lambda_2-\lambda_1)\}\,\rho(x,y).
\nn
\ee
We recall that $\gamma\ge \lambda_2\,(3\lambda_2+\lambda_1)/
(\lambda_2-\lambda_1)$, see \rf{GM-FN-1}, so that
$\dhf(F^{[2]}(x),F^{[2]}(y))\le \gamma\,\rho(x,y)$ for all $x,y\in\Mc$.
\par We have proved that property \rf{F2-NEM} and inequality \rf{HD-RT} hold for the mapping $F^{[2]}$ provided $\lambda_1$, $\lambda_2$ and $\gamma$ satisfy inequalities \rf{GM-FN-1}. In particular, we can set
$\lambda_1=1$, $\lambda_2=3$. For these parameters  \rf{F2-NEM} and \rf{HD-RT} hold provided
$$
\gamma=\lambda_2\,(3\lambda_2+\lambda_1)/
(\lambda_2-\lambda_1)=3\,(3\cdot 3+1)/(3-1)=15.
$$
\par Let now $\X$ be a Euclidean space, and let $\lambda_1$, $\lambda_2$ and $\gamma$ be parameters satisfying inequalities \rf{X1-H}.
In this case, replacing in the above calculations $\theta(L)=(3L+1)/(L-1)$ with
$\theta(L)=1+2L/\sqrt{L^2-1}$ we obtain the following:
\be
\dhf(F^{[2]}(x),F^{[2]}(y))&\le&
L\left(1+2L/\sqrt{L^2-1}\right)\ds(x,y)
\nn\\
&=&
(\lambda_2/\lambda_1)
(1+2(\lambda_2/\lambda_1)/\sqrt{(\lambda_2/\lambda_1)^2-1}
\cdot (\lambda_1\rho(x,y))\nn\\
&=&
\left\{\lambda_2+2\lambda_2^2\,\slbig
\left(\lambda_2^2-\lambda_1^2\right)^{\frac12}\right\}
\,\rho(x,y)\le \gamma\,\rho(x,y).
\nn
\ee
See \rf{X1-H}. This proves that \rf{F2-NEM} and \rf{HD-RT} hold provided $\lambda_1$, $\lambda_2$ and $\gamma$ satisfy inequalities \rf{X1-H}.
\par In particular, we can set  $\lambda_1=1$, $\lambda_2=3$ and $\gamma=10$. In fact, in this case
$$
\lambda_2+2\lambda_2^2\,\slbig
\left(\lambda_2^2-\lambda_1^2\right)^{\frac12}=
3+2\cdot 3^2/\sqrt{3^2-1}=3+18/\sqrt{8}\approx 9.36\le 10=\gamma.
$$
\par The proof of Theorem \reff{X-LSGM} is complete.\bx

\SECT{5. The case $X=\R$ and related results.}{5}
\addtocontents{toc}{5. The case $X=\R$ and related results.\hfill \thepage\par\VST}
\indent

\indent\par {\bf 5.1 Main conjecture in the one dimensional case.}
\addtocontents{toc}{~~~~5.1 Main conjecture in the one dimensional case.\hfill \thepage\par\VST}

\msk
\par In this section we prove Conjecture \reff{BR-IT} for a one dimensional Banach space $\X$. Clearly, we may assume that $X=\R$. Thus, in this case the unit ``ball'' of $\X$ is the interval $I_0=[-1,1]$. Given $a\in\R$ and $r\ge 0$, we set $rI_0=[-r,r]$ and $I(a,r)=[a-r,a+r]$.
\begin{proposition}\label{X-1DIM} Let $\MR$ be a pseudometric space. Let $m=1$ and let $X=\R$; thus, $\ell=\ell(m,X)=1$, see \rf{NMY-1}. In this case Conjecture \reff{BR-IT} holds for every $\lambda_1\ge 1$ and $\gamma\ge 1$.
\par Thus, the following statement is true: Let $F$ be a set-valued mapping from $\Mc$ into the family $\Kc(\R)$ of all closed bounded intervals in $\R$. Suppose that for every $x,y\in\Mc$ there exist points $g(x)\in F(x)$ and $g(y)\in F(y)$ such that $|g(x)-g(y)|\le \rho(x,y)$.
\smsk
\par Let $F^{[1]}(x)$, $x\in\Mc$, be the $\lambda_1$-balanced refinement of the mapping $F$, i.e., the set
\bel{F1-PR1}
F^{[1]}(x)=
\bigcap_{z\in\Mc}\,\left[F(z)+
\lambda_1\,\rho(x,z)\,\BXR\right]~~~~
\text{where}~~~~I_0=[-1,1].
\ee
\par Then $F^{[1]}(x)\ne\emp$ for every $x\in\Mc$, and
$$
\dhf(F^{[1]}(x),F^{[1]}(y))\le \gamma\,\rho(x,y)~~~~\text{for all}~~~~x,y\in\Mc.
$$
\end{proposition}
\msk
\par As in the previous sections, one of the main tools in the proof of Conjecture \reff{BR-IT} will be Helly's Theorem. Let us recall its statement in the one-dimensional case. We will also give a formula for a neighborhood of the intersection of intervals in $\R$.
\begin{lemma}\lbl{H-R} Let $\Kc\subset\Ic(\R)$ be a collection of closed intervals in $\R$. (See \rf{IR-DF1}.)
\smsk
\par (a) Suppose that either $\Kc$ is finite or at least one member of $\Kc$ is bounded.
\par If the intersection of every two intervals from $\Kc$ is non-empty, then there exists a point in $\R$ common to all of the family $\Kc$.
\par (b) Suppose that $\cap\{K:K\in\Kc\}\ne\emp$. Then for every $r\ge 0$ the following equality
$$
\left(\,\bigcap_{K\in\,\Kc} K\right) +r\BXR
=\bigcap_{K\in\,\Kc}\,
\left\{\,K+r\BXR\,\right\}
$$
holds.
\end{lemma}
\par {\it Proof.} In Lemma \reff{H-IN} we have proved an analog of property {\it (b)} for $\RT$. The proof of {\it (b)} is an obvious modification of that proof where we replace Helly's Theorem \reff{H-TH} in $\RT$ with Helly's Theorem in $\R$ formulated in part (a) of the present lemma. We leave the details to the interested reader.\bx
\begin{remark}\lbl{R-S4} {\em We can slightly weaken the hypothesis of one dimensional Helly's theorem given in part $(a)$ of Lemma \reff{H-R} as follows: we may assume that
$(i)$ either $\Kc$ is finite or $(ii)$ there exists a {\it finite} subfamily $\tKc\subset\Kc$ such that the intersection $\cap\{I:I\in\tK\}$ is non-empty and {\it bounded}.
\par Furthermore, $(ii)$ can be replaced with the following requirement: $(ii')$ there exist intervals $I,I'\in\Kc$ such that the intersection $I\cap I'$ is non-empty and {\it bounded}.\rbx}
\end{remark}
\smsk
\par Let $F:\Mc\to \Kc(\R)$ be a set-valued mapping which to every $x\in\Mc$ assigns a closed bounded interval $F(x)=[a(x),b(x)]$. (Thus, $a(x)=\min F(x)$,  $b(x)=\max F(x)$, so that $a(x)\le b(x)$, $x\in\Mc$.) Let
$$
r(x)=\frac{b(x)-a(x)}{2},~~~~c(x)=\frac{a(x)+b(x)}{2},~~~
x\in\Mc.
$$
Thus, $F(x)=I(c(x),r(x))=c(x)+r(x)\BXR$. Clearly,
\bel{DR-XY}
\dist(F(x),F(y))=[|c(x)-c(y)|-r(x)-r(y)]_+=
\max\{[a(x)-b(y)]_+,[a(y)-b(x)]_+\}.
\ee
This formula leads us to the following
\begin{claim}\lbl{TWO} Given $x,y\in\Mc$, and $\lambda\ge 0$, there exist points
$g(x)\in F(x)$, $g(y)\in F(y)$ such that $|g(x)-g(y)|\le\lambda\,\rho(x,y)$ if and only if the following inequality
$$
|c(x)-c(y)|\le r(x)+r(y)+\lambda\,\rho(x,y)
$$
holds. This inequality is equivalent to the inequality
$$
\max\{a(x)-b(y),a(y)-b(x)\}\le\lambda\,\rho(x,y).
$$
\end{claim}
\par {\it Proof.} The claim is immediate from formula \rf{DR-XY} and the following obvious fact: such points $g(x),g(y)$ exist iff $\dist(F(x),F(y))\le \lambda\,\rho(x,y)$.\bx
\smsk
\par Given a set-valued mapping $F(x)=[a(x),b(x)]$, $x\in\Mc$, we set
\bel{LM-F}
\lambda_F=\sup_{x,y\in\Mc}\frac{[a(x)-b(y)]_+}{\rho(x,y)}=
\sup_{x,y\in\Mc}\frac{[\min F(x)-\max F(y)]_+}{\rho(x,y)}.
\ee
\par Note that from inequality \rf{DR-XY}, we have
\bel{LM-F-D}
\lambda_F=\sup_{x,y\in\Mc}
\frac{\dist(F(x),F(y))}{\rho(x,y)}.
\ee
Clearly,
$$
\lambda_F=\sup_{x,y\in\Mc}
\frac{[\min F(x)+\min\{-F(y)\}]_+}{\rho(x,y)}.
$$
(Recall that we set $\frac{0}{0}=0$ and $\frac{A}{0}=+\infty$ for $A>0$.)
\par Given $\lambda\ge 0$, we also introduce the following functions on $\Mc$:
\bel{FP-D1}
f^+[\lambda;F](x)=\inf_{y\in\Mc}\,\{b(y)+\lambda\,\rho(x,y)\}
=\inf_{y\in\Mc}\,\{\max F(y)+\lambda\,\rho(x,y)\},
\ee
\bel{FM-D1}
f^-[\lambda;F](x)=\sup_{y\in\Mc}\,\{a(y)-\lambda\,\rho(x,y)\}
=\sup_{y\in\Mc}\,\{\min F(y)-\lambda\,\rho(x,y)\},
\ee
and
\bel{FS-D1}
f[\lambda;F](x)=\frac{f^+[\lambda;F](x)
+f^-[\lambda;F](x)}{2}.
\ee
\begin{lemma}\lbl{FP-R1} Let $\lambda\ge 0$, and let $F:\Mc\to\Kc(\R)$ be a set-valued mapping.
\smsk
\par (i) (The Finiteness Principle for Lipschitz selections in $\R$.) Suppose that for every $x,y\in\Mc$ the restriction $F|_{\{x,y\}}$ of $F$ to $\{x,y\}$ has a Lipschitz selection $f_{\{x,y\}}$ with $\|f_{\{x,y\}}\|_{\Lip(\{x,y\},\R)}\le \lambda$. Then $F$ has a Lipschitz selection $f:\Mc\to\R$ with Lipschitz seminorm $\|f\|_{\Lip(\Mc,\R)}\le \lambda$.
\par Furthermore, one can set
$$
f=f^+[\lambda;F],~~~f=f^-[\lambda;F]~~~\text{or}~~~
f=f[\lambda;F].
$$
\par (ii) There exists a Lipschitz selection of $F$ if and only if
\bel{LF}
\lambda_F=\sup_{x,y\in\Mc}
\frac{[\min F(x)-\max F(y)]_+}{\rho(x,y)}<\infty.
\ee
Moreover, if this inequality holds then
$$
\lambda_F=\min\{\|f\|_{\Lip(\Mc,\R)}:f~~\text{is a Lipschitz selection of}~~F\}.
$$
\par The above minimum is attained at each of the following functions: $f^+[\lambda;F]$, $f^-[\lambda;F]$ or $f[\lambda;F]$. In other words,
$$
\lambda_F=\|f^+[\lambda;F]\|_{\Lip(\Mc,\R)}=
\|f^-[\lambda;F]\|_{\Lip(\Mc,\R)}=
\|f[\lambda;F]\|_{\Lip(\Mc,\R)}.
$$
\end{lemma}
\par {\it Proof.} (i) Let $F(x)=[a(x),b(x)]$, $x\in\Mc$. Prove that the function
$$
f(x)=f^+[\lambda;F](x)=\inf_{y\in\Mc}
\,\{b(y)+\lambda\,\rho(x,y)\}
$$
is a Lipschitz selection of $F$ with $\|f\|_{\Lip(\Mc,\R)}\le\lambda$.
\par Clearly, $f(x)\le b(x)$ on $\Mc$. (Take $y=x$ in the definition of $f$.) The hypothesis of part (i) of the lemma tells us that for every $x,y\in\Mc$ there exist points $g(x)\in[a(x),b(x)]$, $g(y)\in[a(y),b(y)]$
such that $|g(x)-g(y)|\le \lambda\,\rho(x,y)$. Therefore, thanks to Claim \reff{TWO},
$$
a(x)\le b(y)+\lambda\,\rho(x,y).
$$
Hence, $a(x)\le f(x)$. Thus, $f(x)\in[a(x),b(x)]=F(x)$ proving that $f$ is a selection of $F$ on $\Mc$.
\par Furthermore, thanks to the triangle inequality,
$$
|f(x)-f(y)|=|\inf_{u\in\Mc}
\,\{b(u)+\lambda\,\rho(x,u)\}-\inf_{u\in\Mc}
\,\{b(u)+\lambda\,\rho(y,u)\}|\le
\sup_{u\in\Mc}
\,|\lambda\,\rho(x,u)-\lambda\,\rho(y,u)|\le\lambda\,\rho(x,y)
$$
proving the required inequality $\|f\|_{\Lip(\Mc,\R)}\le\lambda$.
\par In the same way we show that the function $f=f^-[\lambda;F]$ is a Lipschitz selection of $F$ with  Lipschitz seminorm at most $\lambda$. Clearly, the function
$f[\lambda;F]=(f^+[\lambda;F]+f^-[\lambda;F])/2$ has the same property.
\msk
\par (ii) Let $f:\Mc\to\R$ be a Lipschitz selection of $F$ with $\|f\|_{\Lip(\Mc,\R)}\le\lambda$. We know that for every $x,y\in\Mc$ we have $f(x)\in F(x)$,  $f(y)\in F(y)$ and $|f(x)-f(y)|\le \lambda\,\rho(x,y)$. In this case Claim \reff{TWO} tells us that
$$
[a(x)-b(y)]_+=[\min F(x)-\max F(y)]_+\le \lambda\,\rho(x,y).
$$
Hence, $\lambda_F\le \lambda<\infty$, see \rf{LF}.
\par Conversely, suppose that $\lambda_F<\infty$. Then, thanks to \rf{LM-F}, for every $x,y\in\Mc$ we have
$$
a(x)-b(y)\le \lambda_F\,\rho(x,y)~~~\text{and}~~~
a(y)-b(x)\le \lambda_F\,\rho(x,y).
$$
\par This inequality and Claim \reff{TWO} tell us that there exist points $g(x)\in F(x)$, $g(y)\in F(y)$
such that $|g(x)-g(y)|\le \lambda_F\,\rho(x,y)$. In other words, $g$ is a Lipschitz selection of the restriction $F$ to the set $\{x,y\}$ with Lipschitz  seminorm $\|g\|_{\Lip(\{x,y\},\R)}\le \lambda_F$. Therefore, thanks to part (i) of the present lemma, there exists a Lipschitz selection $f:\Mc\to\R$ of $F$ with Lipschitz seminorm $\|f\|_{\Lip(\Mc,\R)}\le\lambda_F$.
\par It remains to note that in the proof of part (i) of the present lemma we have shown that each of the functions
$f^+[\lambda_F;F]$, $f^-[\lambda_F;F]$ or $f[\lambda_F;F]$ provides a Lipschitz selection of $F$ with Lipschitz seminorm at most $\lambda_F$.
\par The proof of the lemma is complete.\bx
\msk
\par {\it Proof of Proposition \reff{X-1DIM}.} We have to prove that the set $F^{[1]}(x)$ is non-empty for each $x\in\Mc$, and for every $x,y\in\Mc$
\bel{TM-1}
\dhf(F^{[1]}(x),F^{[1]}(y))\le \rho(x,y)
\ee
provided the restriction $F|_{\Mc'}$ of $F$ to every two point subset $\Mc'\subset\Mc$ has a Lipschitz selection $f_{\Mc'}:\Mc'\to\R$ with $\|f_{\Mc'}\|_{\Lip(\Mc',\R)}\le \eta$. Thus, for every $z,z'\in\Mc$ there exist points
\bel{GZ-ZP}
g(z)\in F(z),~g(z')\in F(z')~~~\text{such that}~~~ |g(z)-g(z')|\le \rho(z,z').
\ee
\par We recall that the set-valued mapping $F^{[1]}$ is defined by formula \rf{F1-PR1}.
\par Prove that $F^{[1]}(x)\ne\emp$\, for every $x\in\Mc$.
Indeed, thanks to \rf{F1-PR1} and Helly's Theorem for intervals (part (i) of Lemma \reff{H-R}), $F^{[1]}(x)\ne\emp$ provided
\bel{FFP}
(F(z)+\rho(x,z)\,\BXR)\cap(F(z')+\rho(x,z')
\,\BXR)\ne\emp
\ee
for every $z,z'\in\Mc$.
\par We know that there exist points $g(z)$ and $g(z')$ satisfying \rf{GZ-ZP}. Let
$$
a=\min\{g(z)+\rho(z,x),g(z')+\rho(z',x)\}.
$$
Thanks to the inequality $|g(z)-g(z')|\le \rho(z,z')$, we have
$$
g(z)=\min\{g(z),g(z')+\rho(z',z)\}
$$
so that, by the triangle inequality,
$$
|a-g(z)|\le \max\{\rho(z,x),|\rho(z',x)-\rho(z',z)\}
=\rho(z,x).
$$
We also know that $g(z)\in F(z)$, see \rf{GZ-ZP}, so that $a\in F(z)+\rho(x,z)\,\BXR$.
\par In the same way we show  that
$a\in F(z')+\rho(x,z')\,\BXR$ proving the required property
\rf{FFP}.
\smsk
\par Prove that
\bel{FXY}
F^{[1]}(x)+\rho(x,y)\,\BXR\supset F^{[1]}(y)
\ee
for every $x,y\in\Mc$.
\par We know that $F^{[1]}(x)\ne\emp$ which enables us to apply part (b) of Lemma \reff{H-R} to the left hand side of
\rf{FXY}. This lemma and definition \rf{F1-PR1} tell us that
$$
F^{[1]}(x)+\rho(x,y)\,\BXR=
\bigcap_{z\in \Mc}\,
\left[F(z)+\rho(x,z)\,\BXR\right]+\rho(x,y)\,\BXR
=\bigcap_{z\in \Mc}\,
\left[F(z)+(\rho(x,z)+\rho(x,y))\,\BXR\right]
$$
so that, thanks to the triangle inequality,
$$
F^{[1]}(x)+\rho(x,y)\,\BXR\supset
\bigcap_{z\in \Mc}\,
\left[F(z)+\rho(y,z)\,\BXR\right]=F^{[1]}(y)
$$
proving \rf{FXY}. By interchanging the roles of $x$ and $y$ we obtain also
$$
F^{[1]}(y)+\rho(x,y)\,\BXR\supset F^{[1]}(x).
$$
These two inclusions prove the required inequality \rf{TM-1}.
\smsk
\par The proof of Proposition \reff{X-1DIM} is complete.\bx

\indent\par {\bf 5.2 Several useful formulae for the Hausdorff distance.}
\addtocontents{toc}{~~~~5.2 Several useful formulae for the Hausdorff distance.\hfill \thepage\par\VST}
\msk
\par Let $\X$ be a Banach space and let $A,B\subset \X$. We recall formula \rf{HD-R} for the Hausdorff distance between $A$ and $B$:
$$
\dhf(A,B)=\inf\{r>0:~A+\BX(0,r)\supset B~~\text{and}~~
B+\BX(0,r)\supset A\}.
$$
\par We also useful introduce a function
$$
\dq(A,B)=\inf\{r>0:~A+\BX(0,r)\supset B\}=
\sup\{\,\dist(a,B): a\in A\}.
$$
Then,
$$
\dhf(A,B)=\max\{\,\dq(A,B),\dq(B,A)\}.
$$
\par Let us note the following useful formula for the Hausdorff distance, see \cite[p. 144]{DZ}:
$$
\dhf(A,B)=\sup\,\{\,|\dist(x,A)-\dist(x,B)|:x\in\X\}.
$$
\par Next, we recall a well known expression for $\dhf$ in terms of support functions. Let $\X^*$ be the dual space of $\X$, and let $B_{\X^*}$ be unit ball of $\X^*$. We recall
that the support function $h_A:\X^*\to \R$ is defined by
$$
h_A(f)=\sup\,\{f(x):x\in A\}.
$$
\par One can easily see that for every $\alpha,\beta\ge 0$ and every $A,B\subset\X$ we have
\bel{M-SP}
h_{\alpha A+\beta B}=\alpha h_A+\beta h_B.
\ee
\par Furthermore, if $A$ is a convex closed bounded sets
which is symmetric with respect to $0$, then for every $f\in\X^*$ we have
\bel{S-HG}
h_{A}(f)=\sup\{f(x):x\in A\}=\sup\{-f(x):x\in A\}=
h_{A}(-f).
\ee
\par In these settings, for every convex closed bounded subsets $A,B\subset \X$ the following equality
\bel{SF-HD}
\dhf(A,B)=\sup\,\{\,|h_A(f)-h_B(f)|:f\in B_{X^*}\}
\ee
holds. See, e.g., \cite{CV} or \cite{Gru}.
\par Let us also note the following result proven in \cite{Wil}: If $A,B\subset \X$ are non-empty, bounded and convex then $\dhf(A,B)\le \dhf(\partial A,\partial B)$. If
$A,B\subset \X$ are bounded, convex and have non-empty interior, then
$$
\dhf(A,B)=\dhf(\partial A,\partial B).
$$
Here $\partial A$ denotes the boundary of the set $A$.
\begin{lemma}\lbl{SYM} (i) Let $I_k=[a_k,b_k]$, $k=1,2$ be two line segments in $\R$. Then
$$
\dhf(I_1,I_2)=\max\{|a_1-a_2|,|b_1-b_2|\}.
$$
(ii) Let $A_1,A_2\subset\X$ be convex closed bounded sets. Suppose that $A_1, A_2$ are centrally symmetric with respect to points $c_1$ and $c_2$ respectively. Then
$$
\|c_1-c_2\|\le \dhf(A_1,A_2).
$$
\end{lemma}
\par {\it Proof.} (i) Let $\TI=[-1,1]$. Suppose that $I_1\ne I_2$; otherwise the statement (i) is trivial. In this case
$$
\ve=\max\{|a_1-a_2|,|b_1-b_2|\}>0.
$$
Then $a_1-\ve\le a_2\le b_2\le b_1+\ve$ so that
$$
I_1+\ve \TI=[a_1-\ve,b_1+\ve]\supset I_2=[a_2,b_2].
$$
In the same way we show that and $I_2+\ve \TI\supset I_1$, proving that
$$
\dhf(I_1,I_2)\le\ve=\max\{|a_1-a_2|,|b_1-b_2|\}.
$$
\par Prove the converse inequality. Let $\ve$ be a positive number such that
\bel{VE-K}
I_1+\ve\TI=[a_1-\ve,b_1+\ve]\supset I_2=[a_2,b_2]~~~~\text{and}~~~~ I_2+\ve \TI=[a_2-\ve,b_2+\ve]\supset I_1=[a_1,b_1].
\ee
Then
$$
a_1-\ve\le a_2,~~b_2\le b_1+\ve,~~~~\text{and}~~~~
a_2-\ve\le a_1,~~ b_1\le b_2+\ve.
$$
Hence,
$$
\max\{|a_1-a_2|,|b_1-b_2|\}\le \ve.
$$
\par We take the infimum over all $\ve>0$ satisfying \rf{VE-K}, and obtain the required inequality
$$
\max\{|a_1-a_2|,|b_1-b_2|\}\le \dhf(I_1,I_2).
$$
\par (ii) Let $A\subset\X$ be a convex closed bounded set. We assume that $A$ is centrally symmetric with respect to a point $\bar{a}\in\X$. Thus, $A=\bar{A}+\bar{a}$ where $\bar{A}$ is a convex closed bounded set with center of symmetry at $0$. Therefore, thanks to \rf{S-HG}, for every $f\in B_{\X^*}$ we have $h_{\bar{A}}\,(f)=h_{\bar{A}}\,(-f)$.
\par Hence, thanks to this property and \rf{M-SP},
$$
h_A(f)=h_{\bar{a}+\bar{A}}(f)=f(\bar{a})+h_{\bar{A}}(f)
~~~~\text{and}~~~~
h_A(-f)=-f(\bar{a})+h_{\bar{A}}(-f)=-f(\bar{a})+h_{\bar{A}}(f)
$$
proving that
$$
f(\bar{a})=\tfrac12(h_A(f)-h_A(-f)).
$$
\par Applying this formula to the sets $A_1, A_2$ and their centers $c_1,c_2$, we have
\be
|f(c_1)-f(c_2)|&=&
\tfrac12|(h_{A_1}(f)-h_{A_1}(-f))-
(h_{A_2}(f)-h_{A_2}(-f))|\nn\\
&\le&
\tfrac12|h_{A_1}(f)-h_{A_2}(f)|+
\tfrac12|h_{A_1}(-f)-h_{A_2}(-f)|.
\nn
\ee
This inequality and \rf{SF-HD} imply the following:
$$
|f(c_1-c_2)|=|f(c_1)-f(c_2)|\le \dhf(A_1,A_2).
$$
\par Hence,
$$
\|c_1-c_2\|=\sup_{f\in B_{\X^*}}|f(c_1-c_2)|\le \dhf(A_1,A_2)
$$
proving the lemma.\bx
\msk

\indent\par {\bf 5.3 Three criteria for Lipschitz selections.}
\addtocontents{toc}{~~~~5.3 Three criteria for Lipschitz selections.\hfill \thepage\par\VST}

\msk

\par Lemma \reff{SYM} and Theorem \reff{X-LSGM} imply the following Lipschitz selection theorem.
\begin{theorem}\lbl{LSEG-B} Let $\MR$ be a pseudometric space, and let $\X$ be a Banach space. Let $\lambda>0$ and let $F:\Mc\to\Kc_1(\X)$ be a set-valued mapping from $\Mc$ into the family $\Kc_1(\X)$ of all bounded closed line segments in $\X$.
\par Suppose that for every subset $\Mc'\subset\Mc$ with $\#\Mc'\le 4$, the restriction $F|_{\Mc'}$ of $F$ to $\Mc'$ has a Lipschitz selection with Lipschitz seminorm at most $\lambda$.
\par Then $F$ has a Lipschitz selection $f$ with Lipschitz  seminorm $\|f\|_{\Lip(\Mc,\X)}\le 15\lambda$. If $X$ is a Euclidean space, there exists a Lipschitz selection $f$ of $F$ with $\|f\|_{\Lip(\Mc,\X)}\le 10\lambda$.
\end{theorem}
\par {\it Proof.} Let $\vl=(1,3)$ and let $F^{[1]}$ and $F^{[2]}$ be the first and the second order $(\vl,\lambda\rho)$-balanced refinements of $F$. See Definition \reff{BREF}. Thus,
$$
F^{[1]}(x)=
\bigcap_{z\in\Mc}\,
\left[F(z)+\lambda\,\rho(x,z)\BX\right]
~~~\text{and}~~~
F^{[2]}(x)=\bigcap_{z\in\Mc}\,
\left[F^{[1]}(z)+3\lambda\,\rho(x,z)\BX\right],
~~~x\in\Mc.
$$
\par Theorem \reff{X-LSGM} tells us that the set-valued mapping $F^{[2]}$ is a $\gamma$-core of $F$ with $\gamma=15$ provided $\X$ is an arbitrary Banach space, and
with $\gamma=10$ whenever $\X$ is a Euclidean space. In other words, $F^{[2]}(x)\ne\emp$ for every $x\in\Mc$, and
\bel{DH-GM}
\dhf(F^{[2]}(x),F^{[2]}(y))\le \gamma\,\lambda\,\rho(x,y)
~~~~\text{for all}~~~~x,y\in\Mc.
\ee
\par Clearly, $F^{[2]}(x)\in \Kc_1(\X)$, i.e., $F^{[2]}(x)$
is a closed bounded line segment in $\X$ for each $x\in\Mc$. In other words, $F^{[2]}(x)=[a_1(x),a_2(x)]$, $x\in\Mc$, where $a_i:\Mc\to\X$, $i=1,2$, are certain mappings on $\Mc$.
\par We define a mapping $f:\Mc\to\X$ by letting
$$
f(x)=\tfrac12(a_1(x)+a_2(x)), ~~~~x\in\Mc.
$$
\par Thus, $f(x)$ is the center of the line segment
$F^{[2]}(x)=[a_1(x),a_2(x)]$ so that $f(x)\in F^{[2]}(x)\subset F(x)$ proving that $f$ is a selection of $F$ on $\Mc$. Furthermore, Lemma \reff{SYM} and inequality \rf{DH-GM} tell us that
$$
\|f(x)-f(y)\|\le \dhf(F^{[2]}(x),F^{[2]}(y))\le
\gamma\,\lambda\,\rho(x,y)
~~~~\text{for all}~~~~x,y\in\Mc.
$$
Thus, $\|f\|_{\Lip(\Mc,\X)}\le \gamma\lambda$, and the proof of the theorem is complete.\bx
\msk
\par We finish the section with a useful criterion for Lipschitz selections in $\R$. To its formulation, given $\lambda>0$ we set
$$
F^{[1]}_{\lambda}(x)=
\bigcap_{z\in \Mc}\,
\left[F(z)+\lambda\,\rho(x,z)\,I_0\right],~~~~~x\in\Mc.
$$
\begin{proposition}\lbl{CR-FL-1} Let $\mfM=\MR$ be a pseudometric space, and let $F:\Mc\to\Ic(\R)$ be a set-valued mapping. Suppose that either $\Mc$ is finite or there exist $x,y\in\Mc$ and  $\al\ge 0$ such that the set $F(x)\cap [F(y)+\al\,I_0]$ is non-empty and bounded. Then the following criterion holds: given $\lambda>0$ the mapping $F$ has a Lipschitz selection $f:\Mc\to\R$ with $\|f\|_{\Lip(\Mc;\R)}\le\lambda$ if and only if the set $F^{[1]}_{\lambda}(x)\ne\emp$ for every $x\in\Mc$.
\par Furthermore,
$$
|F|_{\mfM,\R}=
\inf\{\lambda: F^{[1]}_{\lambda}(x)\ne\emp~~~\text{for all}~~~x\in\Mc\}.
$$
See \rf{OM-L}.
\end{proposition}
\par {\it Proof.} The proposition easily follows from Lemma \reff{FP-R1}. Indeed, suppose that $F$ has a Lipschitz selection $f:\Mc\to\R$ with $\|f\|_{\Lip(\Mc;\R)}\le\lambda$. Then, given $x\in\Mc$, we have  $|f(x)-f(y)|\le\lambda\,\rho(x,z)$ for every $z\in\Mc$. But $f(x)\in F(x)$ and $f(z)\in F(z)$ (because $f$ is a selection of $F$) so that $f(x)\in F(z)+\lambda\,\rho(x,z) I_0$ proving that $f(x)\in F^{[1]}_{\lambda}(x)$.
\par Now, suppose that $F^{[1]}_{\lambda}(x)\ne\emp$ for every $x\in\Mc$. Then, for every $x,z\in\Mc$, we have
$$
F(x)\cap[F(z)+\lambda\,\rho(x,z)\,I_0]\ne\emp.
$$
Hence, $\dist(F(x),F(z))\le \lambda\,\rho(x,z)$ so that there exist points $g(x)\in F(x)$, $g(z)\in F(z)$ such that $|g(x)-g(y)|\le\lambda\,\rho(x,z)$. Part $(i)$ of Lemma \reff{FP-R1} tells us that in these settings the mapping $F$ has a Lipschitz selection $f:\Mc\to\R$ with $\|f\|_{\Lip(\Mc;\R)}\le\lambda$.
\par The proof of the proposition is complete.\bx
\msk
\par The following proposition is immediate from Proposition \reff{CR-FL-1}.
\begin{proposition} Let $G:\Mc\to\Rc(\RT)$ be a set-valued mapping from a pseudometric space $\mfM=(\Mc,\rho)$ into the family $\RCT$  of all closed rectangles in $\RT$ with sides parallel to the coordinate axes. Suppose that either $\Mc$ is finite or there exist $x,y\in\Mc$ and  $\al\ge 0$ such that the set $G(x)\cap [G(y)+\al\,Q_0]$ is non-empty and bounded.
\par Then the following criterion holds: given $\lambda>0$ the mapping $G$ has a Lipschitz selection $g:\Mc\to\RT$ with $\|g\|_{\Lip(\Mc;\LTI)}\le\lambda$ if and only if the set
$$
G^{[1]}_{\lambda}(x)=
\bigcap_{z\in \Mc}\,
\left[G(z)+\lambda\,\rho(x,z)\,Q_0\right]
$$
is not empty for every $x\in\Mc$. (Recall that $Q_0=[-1,1]^2$.)
\par Furthermore,
$$
|G|_{\mfM,\LTI}=
\inf\{\lambda: G^{[1]}_{\lambda}(x)\ne\emp~~~\text{for all}~~~x\in\Mc\}.
$$
\end{proposition}
\par Recall that we measure the distances in $\RT$ in the uniform norm $\|a\|_{\LTI}=\max\{|a_1|,|a_2|\}$, $a=(a_1,a_2)$.

\SECT{6. Main Theorem in $\LTI$.}{6}
\addtocontents{toc}{6. Main Theorem in $\LTI$. \hfill\thepage\par \VST}
\par {\bf 6.1 Rectangular hulls of plane convex sets.}
\medskip
\addtocontents{toc}{~~~~6.1 Rectangular hulls of plane convex sets.\hfill\thepage\par}

\msk
\par We recall that by $\Ic(\R)$ we denote the family of all closed intervals in $\R$ (bounded or unbounded). We also recall that $\RCT$ is the family of all closed rectangles in $\RT$ with sides parallel to the coordinate axes, i.e.,
$$
\RCT=\{\Pi=I_1\times I_2: I_1,I_2\in\Ic(\R)\}.
$$
We refer to every $\Pi\in\RCT$ as a {\it ``box'' or ``rectangle''}.
\par Clearly, each bounded rectangle $\Pi\in\RCT$ is a centrally symmetric set. We let $\cntr(\Pi)$ denote the center of $\Pi$.
\par Everywhere in this section we let $S$ denote a {\it non-empty convex closed subset} of $\RT$.
\begin{definition}\lbl{REN} {\em We let $\HR[S]$ denote the smallest (with respect to inclusion) rec\-tangle containing $S$. Thus,
$$
\HR[S]=\cap\{\Pi: \Pi\in\RCT, \Pi\supset S\}.
$$
\par We refer to $\HR[S]$ as a {\it ``rectangular hull``} of the set $S$.}
\end{definition}
\msk
\par We let $\Prj_i$ denote the operator of orthogonal projection onto the axis $Ox_i$, $i=1,2$, i.e.,
\bel{PR-J}
\Prj_i[x]=x_i~~~\text{for}~~~x=(x_1,x_2)\in\RT.
\ee
\par Then the rectangular hull of $S$ has the following representation:
\bel{H-DP}
\HR[S]=\Prj_1[S]\times\Prj_2[S].
\ee
\par Note also that $\Prj_i[S]=[a_i,b_i]$ where
\bel{B-PR12}
a_i=\inf\{x_i:x=(x_1,x_2)\in S\}~~~~~\text{and}~~~~~
b_i=\sup\{x_i:x=(x_1,x_2)\in S\}.
\ee
\par Property \rf{H-DP} implies the following: for every rectangle $H\in\RCT$ with center $0$, we have
\bel{NR-HS}
\HR[S+H]=\HR[S]+H.
\ee
In particular, for every $r\ge 0$ the following equality
\bel{N-HS}
\HR[S+r Q_0]=\HR[S]+r Q_0
\ee
holds. This and definition \rf{HD-R} imply the following property of rectangular hulls: Let $S_1,S_2$ be convex closed subsets of $\RT$. Then
\bel{DH-RH1}
\dhf(\HR[S_1],\HR[S_2])\le \dhf(S_1,S_2).
\ee
\par Thus, $\HR[S]$ is the only rectangle for which
\bel{HS-U}
\Prj_1[\HR[S]]=\Prj_1[S]~~~~\text{and}~~~~~
\Prj_2[\HR[S]]=\Prj_2[S].
\ee
\par We note one more obvious property characterizing the rectangular hull: $H(S)$ is the only rectangular such that
\bel{SD-HS}
\HR[S]\supset S~~~~~\text{and each side of}~~~\HR[S]~~~\text{has a common point with}~~~S.
\ee
\par Finally, we have the following obvious formula for $\HR[S]$:
\bel{HE-A}
\HR[S]=(S+Ox_1)\cap (S+Ox_2).
\ee
\par The following important property of rectangular hulls in $\RT$ has been noted in \cite[\S 6]{PY1}; see also \cite[Section 7.2]{FS-2017}.
\begin{lemma}\lbl{CH-RH} Let $S\subset\RT$ be a convex compact set. Then $\cntr(\HR[S])\in S$.
\end{lemma}
\par {\it Proof.} Suppose, given a convex compact set $S\subset\RT$, its center $\cntr(\HR[S])\notin S$.
\par Without loss of generality, we may assume that $\cntr(S)=0$. Thus, $0\notin S$. In this case the separation theorem tells us that there exists a vector $a\in\RT$ such that the inner product $\ip{a,x}>0$ for every $x\in S$.
\par Clearly, there exists a side of $\HR[S]$, say $[AB]$, such that $\ip{a,z}\le 0$ for every $z\in[AB]$. Then, $[AB]\cap S=\emp$. This contradicts property \rf{SD-HS}
of the rectangular hull proving the lemma.\bx
\msk
\par We need the following Helly-type theorem in $\RT$.
\begin{proposition} \lbl{PR-PR} Let $\mfC$ be a family of non-empty convex closed subsets of $\RT$. Suppose that either $\mfC$ is finite or at least one member of $\mfC$ is bounded. If 
\bel{C-PRJ}
\Prj_1[C_1\cap C_1']\,\cbig \Prj_1[C_2\cap C_2']\ne\emp~~~~~\text{for every}~~~~~ C_1,C_1',C_2,C_2'\in\mfC,
\ee
then there exists a point common to all of the family $\mfC$. Furthermore,
\bel{PR-CC}
\Prj_1\left[\,\bigcap_{C\in\,\mfC}\,C\right]=
\bigcap_{C,C'\in\,\mfC}
\Prj_1[C\cap C']\,.
\ee
\end{proposition}
\par {\it Proof.} Condition \rf{C-PRJ} tells us that for every $C,C'\in\mfC$ the set $C\cap C'$ is a non-empty. Clearly, $C\cap C'$ is a convex closed subset of $\RT$, so that its projection onto $Ox_1$, the set $\Prj_1[C\cap C']\subset\R$, is a closed interval in $\R$.
\par From the lemma's hypothesis it follows that either the family $\Kc=\{\Prj_1[C\cap C']:C,C'\in\mfC\}$ is finite or at least one member of $\Kc$ is bounded. Thus, $\Kc$ satisfies the hypothesis of the one dimensional Helly's Theorem formulated in Lemma \reff{H-R}, part (a). Thanks to this lemma,
\bel{V-L}
V=\bigcap_{C,C'\in\,\mfC}\Prj_1[C\cap C']\ne\emp\,.
\ee
\par Fix a point $v\in V$, and set $L=\{w\in\RT:\Prj_1[w]=v\}$. Clearly, $L$ is a straight line through $v$ orthogonal to the axis $Ox_1$.
\par Given $C\in\mfC$ we set $K(C)=C\cap L$. We know that $\Prj_1[C]\ni v$ so that $K(C)\ne\emp$. Furthermore, because $v\in V$, for every $C,C'\in\mfC$ we have $v\in\Prj_1[C\cap C']$ so that there exist $\tw\in C\cap C'$ such that $\Prj_1[\tw]=v$. Hence, $\tw\in L\cap C\cap C'=K(C)\cap K(C')$.
\par Let $\Kc=\{K(C):C\in\mfC\}$. Clearly, all members of $\Kc$ are closed intervals in $L$. We have shown that any two members of $\Kc$ have a common point, so that $\Kc$ also satisfies the hypothesis of part (a) of Lemma \reff{H-R}. This proves the existence of a point in $L$, say $u$, common to all of the family $\Kc$.
\par Thus, $u\in C\cap L$ for each $C\in\mfC$ proving that $u\in\cap\{C:C\in\mfC\}$. At the same time, $u\in L$ so that $\Prj_1[u]=v$. This shows that (i) $\cap\{C:C\in\mfC\}\ne \emp$, and (ii) the left hand side of \rf{PR-CC} contains its right hand side.
Obviously, the left hand side of \rf{PR-CC} is contained in its right hand side, proving that equality \rf{PR-CC} holds.
\smsk
\par The proof of the proposition is complete.\bx
\begin{remark} {\em Remark \reff{R-S4} enables us to slightly modify the hypothesis of Proposition \reff{PR-PR}. Namely, we can replace the requirement ``at least one member of $\mfC$ is bounded'' with ``\,there exists a finite subfamily $\widetilde{\mfC}\subset\mfC$ such that the intersection $\cap\{C:C\in\widetilde{\mfC}\}$ is non-empty and {\it bounded}\,''.
\par Indeed, suppose that there exists a subfamily  $\widetilde{\mfC}\subset\mfC$ having such a property. Let us see that in this case \rf{V-L} holds; then the remaining part of the proof holds as well.
\par We know that Proposition \reff{PR-PR} is true provided $\mfC$ is finite. Applying this proposition to $\widetilde{\mfC}$ we conclude that formula \rf{PR-CC} holds for $\widetilde{\mfC}$, i.e.,
$$
\Prj_1\left[\,\bigcap_{C\in\,\widetilde{\mfC}}\,C\right]=
\bigcap_{C,C'\in\,\widetilde{\mfC}}
\Prj_1[C\cap C']\,.
$$
\par Because the set $\cap\{C:C\in\widetilde{\mfC}\}$ is non-empty and bounded, the set
$\cap\{\Prj_1[C\cap C']:C,C'\in\,\widetilde{\mfC}\}$ is non-empty and bounded as well. Therefore, the family $\Kc=\{\Prj_1[C\cap C']:C,C'\in\mfC\}$ satisfies the hypothesis of the one dimensional Helly's Theorem formulated in Lemma \reff{H-R}, part (a), with modification given in Remark \reff{R-S4}.
\par This implies the required statement \rf{V-L} proving the proposition.\rbx}
\end{remark}
\msk
\par Proposition \reff{PR-PR} and properties \rf{HS-U}, \rf{H-DP} of rectangle hulls imply the following
\begin{corollary}\lbl{INT-RE} Let $\mfC$ be a family of convex closed subsets of $\RT$. Suppose that either $\mfC$ is finite or there exists a finite subfamily $\widetilde{\mfC}\subset\mfC$ such that the intersection $\cap\{C:C\in\widetilde{\mfC}\}$ is non-empty and bounded.
If 
\bel{P1-C}
\Prj_1[C_1\cap C_1']\,\cbig \Prj_1[C_2\cap C_2'] \ne\emp~~~~~\text{for every}~~~~~ C_1,C_1',C_2,C_2'\in\mfC,
\ee
then $\cap\{C:C\in\mfC\}\ne\emp$. Furthermore, in this case
\bel{PH-C}
\HR\left[\cap\{C:C\in\mfC\}\right] =
\cap\{\HR[C\cap C']: C,C'\in\mfC\}\,.
\ee
\end{corollary}
\par Let us formulate two useful properties of rectangles from the family $\RCT$.
\begin{lemma}\lbl{D-REC} For every $r_1,r_2\ge 0$ and every
two rectangles $\Pi_1,\Pi_2\in\RCT$ we have
$$
\dist(\,\Pi_1+r_1 Q_0,\Pi_2+r_2 Q_0)=
[\dist(\Pi_1,\Pi_2)-r_1-r_2]_+\,.
$$
\end{lemma}
\begin{lemma}\lbl{R-12} Let $\Rc_1,\Rc_2\subset\mfR(\RT)$ be two families of rectangles in $\RT$. Suppose that each family has a non-empty intersection. Then
$$
\dist\left(\,\bigcap_{\Pi\in\Rc_1}\Pi,
\bigcap_{\Pi\in\Rc_2}\Pi\right)=
\sup_{\Pi_1\in\,\Rc_1,\Pi_2\in\,\Rc_2}
\dist(\,\Pi_1,\Pi_2)\,.
$$
\end{lemma}
\par We prove both lemmas by projecting onto coordinate axes, i.e., by reduction to the one dimensional case. In this case the first lemma is elementary, while the second lemma easily follows from the one dimensional Helly's Theorem.
\par The next lemma is immediate from part (b) of Lemma \reff{H-R}.
\begin{lemma}\lbl{AAA} Let $\Kc\subset \RCT$ be a family of rectangles with non-empty intersection. Let $H\in\RCT$ be a rectangle with center $0$. Then
$$
\left(\,\bigcap_{\Pi\in\,\Kc} \Pi\right)+H
=\bigcap_{\Pi\in\,\Kc}\,
\left\{\,\Pi+H\right\}.
$$
\end{lemma}
\par The following three lemmas are certain modifications of Lemma \reff{H-IN} for the space $\LTI$.
\begin{lemma}\lbl{LTI-H-IN} Let $\Kc$ be a collection of convex closed subsets of $\RT$ with non-empty intersection, and let $\Pi\in\RCT$ be a rectangle with center $0$. Then
$$
\left(\,\bigcap_{K\in\,\Kc} K\right) +\Pi
=\bigcap_{K,K'\in\,\Kc}\,
\left\{\,\left(\,K \,\cbig\, K'\right)+\Pi\right\}.
$$
\end{lemma}
\par {\it Proof.} If the rectangle $\Pi$ is bounded then the lemma is immediate from Lemma \reff{H-IN}. If $\Pi$ is unbounded then the lemma is immediate from Proposition \reff{PR-PR} and \rf{NR-HS}. We leave the details to the interested reader.\bx
\begin{lemma}\lbl{2K-H} Let $K_1,K_2\subset\RT$ be convex closed sets with non-empty intersection. Then for every rectangle $\Pi\in\RCT$ with $\cntr(\Pi)=0$ we have
\bel{K-12}
K_1\cap K_2+ \Pi=(K_1+ \Pi)\cap (K_2+ \Pi)\cap \HR[K_1\cap K_2+\Pi].
\ee
\end{lemma}
\par {\it Proof.} Clearly, the right hand side of \rf{K-12} contains its left hand side.
\par Let us prove the converse statement. Fix a point
\bel{X-H}
x\in
(K_1+ \Pi)\cap (K_2+ \Pi)\cap \HR[K_1\cap K_2+\Pi]
\ee
and prove that $x\in K_1\cap K_2+ \Pi$.
\par Clearly, this property holds if and only if $(x+\Pi)\cap K_1\cap K_2\ne\emp$. Let us represent the rectangle  $x+\Pi$ in the form $x+\Pi=\Pi_1(x)\cap \Pi_2(x)$ where
\bel{PI-2}
\Pi_1(x)=x+Ox_1+\Pi ~~~~\text{and}~~~~~
\Pi_2(x)=x+Ox_2+\Pi.
\ee
(Recall that $Ox_1=\{x=(t,0):t\in\R\}$ and  $Ox_2=\{x=(0,t):t\in\R\}$ are the coordinate axes.) Thus,
$x\in K_1\cap K_2+ \Pi$ provided
$K_1\cap K_2\cap \Pi_1(x)\cap \Pi_2(x)\ne\emp$.
\smsk
\par Helly's Theorem \reff{H-TH} tells us that this statement is true provided any three members of the family
of sets $\Kc=\{K_1,K_2,\Pi_1(x),\Pi_2(x)\}$ have a common point. Let us see that this property holds for $x$
satisfying \rf{X-H}.
\smsk
\par Clearly, for every $i=1,2$,
$$
K_i\cap \Pi_1(x)\cap \Pi_2(x)=K_i\cap (x+\Pi)\ne\emp
$$
because $x\in K_i+\Pi$. Prove that
\bel{P-K12}
K_1\cap K_2\cap \Pi_1(x)\ne\emp.
\ee
\par Indeed, thanks to \rf{HE-A},
$$
\HR[K_1\cap K_2]=\{K_1\cap K_2+Ox_1\}\cap
\{K_1\cap K_2+Ox_2\}\subset K_1\cap K_2+Ox_1
$$
so that
$$
\HR[K_1\cap K_2]+ \Pi\subset K_1\cap K_2+Ox_1+ \Pi
$$
\par But, thanks to \rf{N-HS}, $\HR[K_1\cap K_2]+ \Pi =\HR[K_1\cap K_2+ \Pi]$, and, thanks to \rf{X-H}, $x\in \HR[K_1\cap K_2+ \Pi]$. Hence, $x\in K_1\cap K_2+Ox_1+\Pi$. Clearly, this property is equivalent to \rf{P-K12}, see \rf{PI-2}. In the same fashion we prove that $K_1\cap K_2\cap \Pi_2(x)\ne\emp$ completing the proof of the lemma.\bx
\msk
\par This lemma and Lemma \reff{LTI-H-IN} imply the following result.
\begin{lemma}\lbl{ALL-H} Let $\Kc\subset\Kc(\RT)$ be a family of convex closed subsets of $\RT$ with non-empty intersection. Then for every rectangle $\Pi\in\RCT$ with center $0$ the following equality
$$
\left(\,\bigcap_{K\in\,\Kc} K\right)+\Pi
=\left\{\bigcap_{K\in\,\Kc}\,\left(K+\Pi\right)\right\}
\bigcap
\left\{\,\bigcap_{K,K'\in\,\Kc}\,
\HR[\,K \,\cbig\, K'+\Pi]\,\right\}
$$
holds.
\end{lemma}
\par The last result of this section, Proposition \reff{FP-REC} below, presents the Finiteness Principle for Lipschitz selections for rectangles in $\RT$. Part (i) of this result is immediate from the Finiteness Principle for intervals in $\R$ given in part (i) of Lemma \reff{FP-R1}, and part (ii) is immediate from Proposition \reff{X-1DIM}. Recall that in the one dimensional case the finiteness constant $N(1,\R)=\min\{2^2,2\}=2$, see \rf{NMY-1}, and the constant $\gamma$ from Theorem \reff{MAIN-FP} equals $1$.
\begin{proposition}\lbl{FP-REC} Let $(\Mc,\rho)$ be a pseudometric space, and let $\lambda>0$. Let  $\Tc:\Mc\to\RCT$ be a set-valued mapping. Suppose that either $\Mc$ is finite or $\Tc(x)$ is bounded for some $x\in\Mc$. Let us also assume that for every $x,y\in\Mc$ the restriction $\Tc|_{\{x,y\}}$ of $\Tc$ to $\{x,y\}$ has a Lipschitz selection $g_{x,y}$ with Lipschitz  seminorm $\|g_{x,y}\|_{\Lip(\{x,y\},\LTI)}\le \lambda$. In these settings, the following statements hold:
\smsk
\par (i) The mapping $\Tc$ has a Lipschitz selection $g$ with Lipschitz  seminorm $\|g\|_{\Lip(\Mc,\LTI)}\le \lambda$;
\smsk
\par (ii) Let
$$
\Tc^{[1]}(x)=
\bigcap_{z\in\Mc}\,\left[\Tc(z)+
\lambda\,\rho(x,z)\,Q_0\right],
~~~~~~x\in\Mc,
$$
be the $\lambda$-balanced refinement of the mapping $\Tc$. Then $\Tc^{[1]}(x)\ne\emp$ for each $x\in\Mc$, and
$$
\dhf(\Tc^{[1]}(x),\Tc^{[1]}(y))\le \lambda\,\rho(x,y)
~~~~~\text{for every}~~~~~x,y\in\Mc.
$$
\end{proposition}

\par {\bf 6.2. Balanced refinements of set-valued mappings in $\LTI$.}
\medskip
\addtocontents{toc}{~~~~6.2 Balanced refinements of set-valued mappings in $\LTI$. \hfill\thepage\VST\par}

\msk
\par Theorem \reff{MAIN-RT} tells us that for the space  $\X=\LTI$ given $\lambda_1,\lambda_2,\gamma>0$ and a set-valued mapping $F:\Mc\to\Kc(\X)$, the mapping $F^{[2]}$ defined by \rf{F12-LTI} has properties \rf{F2-NEM} and \rf{HD-RT} provided $\lambda_1\ge 1$, $\lambda_2\ge 3\lambda_1$ and $\gamma\ge \lambda_2\,(3\lambda_2+\lambda_1)^2/
(\lambda_2-\lambda_1)^2$.
\par In this section we show that this result can be improved as follows.
\begin{theorem}\lbl{LTI-M} Let $\mfM=\MR$ be a pseudometric space. Let $F:\Mc\to\Kc(\RT)$ be a set-valued mapping
such that for every $\Mc'\subset\Mc$ with $\#\Mc'\le 4$, the restriction $F|_{\Mc'}$ of $F$ to $\Mc'$ has a Lipschitz selection $f:\Mc\to\LTI$ with Lipschitz seminorm $\|f\|_{\Lip(\Mc,\LTI)}\le 1$.
\par Then for every
\bel{PR-LTI}
\lambda_1\ge 1,~~~~~\lambda_2\ge 3\lambda_1,~~~~\text{and}~~~~~\gamma\ge \lambda_2\,(3\lambda_2+\lambda_1)/
(\lambda_2-\lambda_1)
\ee
properties \rf{F2-NEM} and \rf{HD-RT} hold.
\par In particular, \rf{F2-NEM} and \rf{HD-RT} hold provided $\lambda_1=1$, $\lambda_2=3$ and $\gamma=15$.
\end{theorem}
\par {\it Proof.} We mainly follow the scheme of the proof
of Theorem \reff{MAIN-RT} given in Section 3. We recall that Lipschitz extension constant $e(\mfM,\LTI)=1$, see \rf{LNI-E}.
\par Let $F:\Mc\to\Kc(\RT)$ be a set-valued mapping satisfying the hypothesis of Theorem \reff{LTI-M}. As in Section 3, this enables us to make the following
\begin{assumption}\lbl{LTI-A} For every $\Mc'\subset\Mc$, $\#\Mc'\le 4$, the restriction $F|_{\Mc'}$ of $F$ to $\Mc'$ has a $\rho$-Lipschitz selection $f_{\Mc'}:\Mc'\to \LTI$ with $\rho$-Lipschitz seminorm $\|f_{\Mc'}\|_{\Lip((\Mc',\,\rho),\LTI)}\le 1$.
\end{assumption}
\par We fix a constant $L\ge 3$ and a constant $\al\ge 1$,
and introduce a pseudometric $\ds(x,y)=\al\rho(x,y)$, $x,y\in\Mc$. Then we introduce set-valued mappings $F^{[1]}$ and  $F^{[2]}$ defined by
\bel{F-1GR}
F^{[1]}(x)=
\bigcap_{z\in\Mc}\,
\left[F(z)+\ds(x,z)\,Q_0\right],~~~~~x\in\Mc,
\ee
and
\bel{F-2GR}
F^{[2]}(x)=\bigcap_{z\in\Mc}\,
\left[F^{[1]}(z)+L\ds(x,z)\,Q_0\right],~~~x\in\Mc.
\ee
\par Recall that $F^{[1]}$ and $F^{[2]}$ are the first and the second order $(\{1,L\},\ds)$-balanced refinements of $F$ respectively. See Definition \reff{F-IT}.
\smsk
\par Lemma \reff{G-NE1} and Proposition \reff{N-EM} tell us that $F^{[1]}(x)\ne\emp$ and $F^{[2]}(x)\ne\emp$ for every $x\in\Mc$. Thus, our aim is to show that for every $\al\ge 1$, $L\ge 3$, and every $x,y\in\Mc$ the following inequality
\bel{H-LTI}
\dhf(F^{[2]}(x),F^{[2]}(y))\le \tgm(L)\,\ds(x,y)
\ee
holds with
\bel{GTL-1}
\tgm(L)=L\,\theta(L)~~~~~\text{where}~~~~
\theta(L)=(3L+1)/(L-1)
\ee
is the constant from Proposition \reff{P-F3} (for the space $X=\LTI$).
\par We prove this inequality with the help of a certain  modification of representations \rf{G-XP} and \rf{GX-1}.
To its formulation we recall that given $x,u,u',u''\in\Mc$ we set
\bel{H-D2}
T_x(u,u',u'')=
\{F(u')+\ds(u',u)Q_0\}\cap \{F(u'')+\ds(u'',u)Q_0\}+L\ds(u,x)Q_0\,.
\ee
See definition \rf{H-D}. Recall also that, thanks to \rf{G-XP},
\bel{F2H-4}
F^{[2]}(x)=\bigcap_{u,u',u''\in\Mc} T_x(u,u',u''),~~~~x\in\Mc.
\ee
\par The next lemma provides another representation of the set $F^{[2]}(x)$.
\begin{lemma} For every $x\in\Mc$ the following equalities
\bel{RH-12}
F^{[2]}(x)=F^{[1]}(x)
\bigcap\left\{
\bigcap_{u,u',u''\in\Mc}
\HR\left[T_x(u,u',u'')\right]\right\},
\ee
\bel{RH-12-A}
F^{[2]}(x)=F^{[1]}(x)
\bigcap\left\{
\bigcap_{u\in\Mc}
\HR\left[F^{[1]}(u)+L\ds(u,x)Q_0\right]\right\}
\ee
hold.
\end{lemma}
\par {\it Proof.} Let
\bel{KX}
\Kc_x=\{F(y)+\ds(x,y)Q_0:y\in\Mc\}.
\ee
Then
$$
F^{[1]}(x)=\bigcap_{K\in\Kc_x}\,K~~~~~\text{and}~~~~~
F^{[2]}(x)=\bigcap_{y\in\Mc}
\left\{\left(\bigcap_{K\in\Kc_y}\,K\right)+
L\ds(x,y)Q_0\right\}.
$$
See \rf{F-1GR} and \rf{F-2GR}.
\par Because $F^{[1]}(y)\ne\emp$, the family of sets
$\Kc_y=\{F(z)+\ds(y,z)Q_0:z\in\Mc\}$ has non-empty intersection for every $y\in\Mc$. Therefore, thanks to Lemma \reff{ALL-H},
$$
\left(\bigcap_{K\in\Kc_y}\,K\right)+L\ds(x,y)Q_0=
\left\{\bigcap_{K\in\Kc_y}\,\left(K+L\ds(x,y)Q_0\right)
\right\}
\bigcap
\left\{\bigcap_{K,K'\in\Kc_y}\,
\HR[K\cap K'+L\ds(x,y)Q_0]\right\}.
$$
Hence,
$$
F^{[2]}(x)=\left\{\bigcap_{y\in\Mc}\,
\bigcap_{K\in\Kc_y}\,\left(K+L\ds(x,y)Q_0\right)\right\}
\bigcap
\left\{\bigcap_{y\in\Mc}\,\bigcap_{K,K'\in\Kc_y}\,
\HR[K\cap K'+L\ds(x,y)Q_0]\right\}=A_1\cap A_2.
$$
Clearly, thanks to the triangle inequality,
\be
A_1&=&\bigcap_{y\in\Mc}\,
\bigcap_{K\in\Kc_y}\,\left(K+L\ds(x,y)Q_0\right)=
\bigcap_{y,z\in\Mc}\,
\left(F(z)+\ds(z,y)Q_0+L\ds(x,y)Q_0\right)
\nn\\
&\supset&
\bigcap_{z\in\Mc}\,
\left(F(z)+\ds(z,x)Q_0\right)
=F^{[1]}(x).\nn
\ee
\par On the other hand, $A_1\subset\cap\{K:K\in \Kc_x\}=F^{[1]}(x)$ so that $A_1=F^{[1]}(x)$. This equality, definition \rf{KX} and definition \rf{H-D2} imply \rf{RH-12}. Equality \rf{RH-12-A} is immediate from \rf{RH-12}, Corollary \reff{INT-RE} and Lemma \reff{AAA}.
\par The proof of the lemma is complete.\bx
\begin{lemma}\lbl{F2-N} For every $x\in\Mc$ and every
rectangle $\Pi\in\RCT$ with center $0$ we have
$$
F^{[2]}(x)+\Pi=
\bigcap_{v,u,u',u''\in\Mc}\,\,
\{\,(\HR[T_x(u,u',u'')]\cap(F(v)+\ds(x,v)Q_0))+\Pi\}.
$$
\end{lemma}
\par {\it Proof.} Let $\Kc^{(1)}=\{F(v)+\ds(x,v)Q_0:v\in\Mc\}$, and let
\bel{K2-D}
\Kc^{(2)}=\{\HR\left[T_x(u,u',u'')\right]:u,u',u''\in\Mc\}.
\ee
\par We have to prove that
$$
F^{[2]}(x)+\Pi=\cap\{(K_1\cap K_2)+\Pi:
K_1\in\Kc^{(1)},K_2\in\Kc^{(2)}\}.
$$
\par Formula \rf{RH-12} and Lemma \reff{LTI-H-IN}
tell us that
\bel{F2-K}
F^{[2]}(x)+\Pi=\cap\{(K\cap K')+\Pi:
K,K'\in \Kc^{(1)}\cup \Kc^{(2)}\}.
\ee
\par Let
\bel{A-12}
A=\cap\{(K_1\cap K_2)+\Pi:
K_1\in\Kc^{(1)},K_2\in\Kc^{(2)}\}.
\ee
Formula \rf{F2-K} tells us that
$A\subset F^{[2]}(x)+\Pi$. Prove the converse inclusion.
\smsk
\par It suffices to show that for every $K,K'\in \Kc^{(1)}\cup \Kc^{(2)}$, we have
\bel{KK-A}
A\subset (K\cap K')+\Pi.
\ee
Clearly, thanks to definition \rf{A-12}, it is true provided  $K\in\Kc^{(1)},K'\in\Kc^{(2)}$ or $K\in\Kc^{(2)},K'\in\Kc^{(1)}$.
\par Prove \rf{KK-A} for sets
$K=F(z)+\ds(z,x)Q_0$ and $K'=F(z')+\ds(z',x)Q_0$ which belong to the family $\Kc^{(1)}$. In this case, thanks to \rf{H-D2} and \rf{K2-D}, the element
$$
\tH=\HR\left[(F(z)+\ds(z,x)Q_0)\cap
(F(z')+\ds(z',x)Q_0)\right]=
\HR\left[T_x(x,z,z')\right]\in\Kc^{(2)}.
$$
Lemma \reff{2K-H} and \rf{N-HS} tell us that
$$
(K\cap K')+\Pi=(K+\Pi)\cap(K'+\Pi) \cap(\tH+\Pi).
$$
Therefore, thanks to \rf{A-12}, $(K\cap K')+\Pi\supset A$.
\par Let us prove \rf{KK-A} for $K,K'\in \Kc^{(2)}$. In this case, the sets $K$ and $K'$ are rectangles with sides parallel to the coordinate axes and with non-empty intersection, so that, thanks to Lemma \ref{AAA},
$$
(K\cap K')+\Pi=(K+\Pi)\cap(K'+\Pi).
$$
It remains to note that, thanks to \rf{A-12}, $K+\Pi\supset A$ and $K'+\Pi\supset A$ for every $K,K'\in \Kc^{(2)}$. This proves \rf{KK-A} in the case under consideration completing the proof of the lemma.\bx
\begin{lemma}\lbl{F2-NEW} For every $x\in\Mc$ and every rectangle $\Pi\in\RCT$ with center $0$ the following representations
\be
F^{[2]}(x)+\Pi&=&
\bigcap_{v,u,u',u''\in\Mc}\,\,
\{\,(T_x(u,u',u'')\cap(F(v)+\ds(x,v)Q_0))+\Pi\,\}
\nn\\
&=&
\bigcap_{v,u,u',u''\in\Mc}\,\,
\{\,(T_x(u,u',u'')\cap T_x(x,v,v))+\Pi\,\}
\nn
\ee
hold.
\end{lemma}
\par {\it Proof.} Thanks to \rf{H-D2},
$F(v)+\ds(x,v)Q_0=T_x(x,v,v)$ which proves the second equality of the lemma. Representation \rf{F2H-4} and Lemma \reff{LTI-H-IN} tell us that
$$
F^{[2]}(x)+\Pi=
\bigcap\,\,
\{\,T_x(u,u',u'')\cap T_x(v,v',v'')+\Pi\,\}
$$
where the intersection is taken over all $u,u',u'', v,v',v''\in\Mc$. Hence,
$$
F^{[2]}(x)+\Pi\subset
\bigcap_{v,u,u',u''\in\Mc}\,\,
\{\,(T_x(u,u',u'')\cap T_x(v,v',v''))+\Pi\,\}.
$$
\par On the other hand, Lemma \reff{F2-N} tells us that
\be
F^{[2]}(x)+\Pi&=&
\bigcap_{v,u,u',u''\in\Mc}\,\,
\{\,(\HR[T_x(u,u',u'')]\cap(F(v)+\ds(x,v)Q_0))+\Pi\}\nn\\
&\supset&
\bigcap_{v,u,u',u''\in\Mc}\,\,
\{\,(T_x(u,u',u'')\cap(F(v)+\ds(x,v)Q_0))+\Pi\}
\nn\\
&=&
\bigcap_{v,u,u',u''\in\Mc}\,\,
\{\,(T_x(u,u',u'')\cap T_x(x,v,v))+\Pi\}.
\nn
\ee
proving the lemma.\bx
\begin{remark} {\em Lemma \reff{F2-NEW} is a refinement of representation \rf{GX-1} for the space $\X=\LTI$.\rbx}
\end{remark}
\msk
\par Representation \rf{F2H-4} and Corollary \reff{INT-RE} imply the following formula for the rectangular hull of the set $F^{[2]}(x)$:
\bel{RH-TT}
\HR[F^{[2]}(x)]=
\bigcap\,\,
\HR\left[T_x(u,u',u'')\cap T_x(v,v',v''))\right]
\ee
where the intersection is taken over all
$u,u',u'',v,v',v''\in\Mc$.
\smsk
\par Lemma \reff{F2-NEW} enables us to prove a stronger version of representation \rf{RH-TT}.
\begin{lemma} For every $x\in\Mc$ the rectangular hull of $F^{[2]}(x)$ has the following representation:
$$
\HR[F^{[2]}(x)]=
\bigcap_{v,u,u',u''\in\Mc}\,\,
\HR\left[T_x(u,u',u'')\cap(F(v)+\ds(x,v)Q_0)\right].
$$
\par Recall that~ $T_x(u,u',u'')=
\{F(u')+\ds(u',u)Q_0\}\cap \{F(u'')+\ds(u'',u)Q_0\}+L\ds(u,x)Q_0$\,;\, see \rf{H-D2}.
\end{lemma}
\par {\it Proof.} Given $x,v,u,u',u''\in\Mc$ we put
$$
V_x[v,u,u',u'']=T_x(u,u',u'')\cap (F(v)+\ds(x,v)Q_0).
$$
\par Lemma \reff{F2-NEW} tells us that
$$
F^{[2]}(x)+Ox_i=
\bigcap_{v,u,u',u''\in\Mc}
(V_x[v,u,u',u'']+Ox_i),~~~~~i=1,2.
$$
In turn, thanks to \rf{HE-A},
$\HR[F^{[2]}(x)]=(F^{[2]}(x)+Ox_1)\cap (F^{[2]}(x)+Ox_2)$,
so that
$$
\HR[F^{[2]}(x)])=
\left\{\bigcap_{v,u,u',u''\in\Mc}
(V_x[v,u,u',u'']+Ox_1)\right\}
\bigcap
\left\{\bigcap_{v,u,u',u''\in\Mc}
(V_x[v,u,u',u'']+Ox_2)\right\}.
$$
Hence,
$$
\HR[F^{[2]}(x)])=
\bigcap_{v,u,u',u''\in\Mc}
(V_x[v,u,u',u'']+Ox_1)\cap(V_x[v,u,u',u'']+Ox_2).
$$
\par From this and \rf{HE-A}, we obtain the required representation
$$
\HR[F^{[2]}(x)])=\bigcap_{v,u,u',u''\in\Mc}
\HR[V_x[v,u,u',u'']]
$$
proving the lemma.\bx
\msk
\par We are in a position to prove inequality \rf{H-LTI}.
Our proof will follow the scheme of the proof of Proposition \reff{HD-G1}.
\smsk
\par Let $x,y\in\Mc$. We know that
\bel{Y-H2}
F^{[2]}(y)=\bigcap_{u,u',u''\in\Mc} T_y(u,u',u''),
~~~~~\text{(see \rf{F2H-4})}.
\ee
\par Let $\tau=\tgm(L)\ds(x,y)$. Lemma \reff{F2-NEW} tells us that
$$
F^{[2]}(x)+\tau Q_0=
\bigcap_{v,u,u',u''\in\Mc}\,\,
\{\,(T_x(u,u',u'')\cap(F(v)+\ds(x,v)Q_0))+\tau Q_0\}.
$$
\par Let us fix elements $u,u',u'',v\in\Mc$ and set
$$
A=(T_x(u,u',u'')\cap(F(v)+\ds(x,v)Q_0))+\tau Q_0.
$$
Prove that $A\supset F^{[2]}(y)$. Let
\bel{C123-2}
C_1=F(u')+\ds(u',u)Q_0,~~~C_2=F(u'')+\ds(u'',u)Q_0X,~~~
C=F(v)+\ds(x,v)Q_0,
\ee
and let
$$
\ve=L\,\ds(x,y)~~~~~\text{and}~~~~~r=\ds(u,x)\,.
$$
See \rf{GTL-1}. Then
$$
\tau=\tgm(L)\ds(x,y)=L\,\theta(L)\ds(x,y)=\theta(L)\,\ve\,.
$$
and
\be
A&=&\{[(F(u')+\ds(u',u)Q_0)\cap (F(u'')+\ds(u'',u)Q_0)]
\cap (F(v)+\ds(x,v)Q_0)+
L\,\ds(u,x)Q_0
\nn\\
&=&
(C_1\cap C_2+Lr Q_0)\cap C+\theta(L)\,\ve\,Q_0.
\nn
\ee
\par Let us verify condition \rf{A-PT} of Proposition \reff{P-F3}, i.e., the condition
\bel{PT-2}
C_1\cap C_2\cap(C+rQ_0)\ne\emp.
\ee
\par Let $\Mc'=\{u,u',v\}$. Clearly, $\#\Mc'\le 4$, so that, thanks to Assumption \reff{LTI-A}, there exists a $\rho$-Lipschitz selection $f_{\Mc'}:\Mc'\to\LTI$ of the restriction $F|_{\Mc'}$ with $\|f_{\Mc'}\|_{\Lip((\Mc',\,\rho),\LTI)}\le 1$.
\par Because $e(\mfM,\LTI)=1$ (see \rf{LNI-E}) and $\ds=\al\rho\ge\rho$, the mapping $f_{\Mc'}:\Mc'\to\LTI$
can be extended to a $\ds$-Lipschitz mapping $\tf:\Mc\to\LTI$ defined on all of $\Mc$ with $\ds$-Lipschitz seminorm
$$
\|\tf\|_{\Lip((\Mc,\ds),\LTI)}\le
\|f_{\Mc'}\|_{\Lip((\Mc',\,\rho),\LTI)}\le 1.
$$
\par In particular, $\tf(u')=f_{\Mc'}(u')\in F(u')$, $\tf(u'')=f_{\Mc'}(u'')\in F(u'')$, $\tf(v)=f_{\Mc'}(v)\in F(v)$,
$$
\|\tf(u')-\tf(u)\|\le \ds(u',u),~~~
\|\tf(u'')-\tf(u)\|\le \ds(u'',u)
$$
and
$$
\|\tf(x)-\tf(u)\|\le \ds(u,x)=r,~~~~~
\|\tf(x)-\tf(v)\|\le \ds(v,x).
$$
Hence, $\tf(u)\in C_1\cap C_2$ and $\tf(x)\in C$, so that
$C_1\cap C_2\cap(C+rQ_0)\ni\tf(u)$ proving \rf{PT-2}.
\par This enables us to apply Proposition \reff{P-F3} to the sets $C_1$, $C_2$ and $C$. This proposition tells us that
\begin{align}
A&=(C_1\cap C_2+LrQ_0)\cap C+\theta(L)\,\ve\,Q_0
\nn\\
&\supset [C_1\cap C_2+(Lr+\ve)Q_0]
\cap [(C_1+rQ_0)\cap C +\ve Q_0]
\cap [(C_2+rQ_0)\cap C +\ve Q_0]
\nn\\
&=S_1\cap S_2\cap S_3.
\nn
\end{align}
\par Prove that
$S_i\supset F^{[2]}(y)$ for every $i=1,2,3$. We begin with the set
\be
S_1&=&C_1\cap C_2+(Lr+\ve)Q_0\nn\\
&=&
\{F(u')+\ds(u',u)Q_0\}\cap \{F(u'')+\ds(u'',u)Q_0\}
+(L\ds(u,x)+L\ds(x,y))Q_0\,.
\nn
\ee
See \rf{C123-2}. The triangle inequality tells us that
$\ds(u,x)+\ds(x,y)\ge \ds(u,y)$
so that
$$
S_1\supset
\{F(u')+\ds(u',u)\BX\}\cap \{F(u'')+\ds(u'',u)Q_0\}
+L\ds(u,y)\BX=T_y(u,u',u'')\,.
$$
From \rf{Y-H2} we have $T_y(u,u',u'')\supset F^{[2]}(y)$ proving the required inclusion $S_1\supset F^{[2]}(y)$.
\smsk
\par Prove that $S_2\supset F^{[2]}(y)$. We have
\be
S_2&=&
(C_1+rQ_0)\cap C +\ve Q_0\nn\\
&=&
\{(F(u')+\ds(u',u)Q_0)+\ds(x,u)Q_0\}\cap \{F(v)+\ds(x,v)Q_0\}+L\ds(x,y)Q_0\,.
\nn
\ee
Therefore, thanks to the triangle inequality and \rf{Y-H2},
$$
S_2\supset
\{(F(u')+\ds(u',x)Q_0\}\cap \{F(v)+\ds(x,v)Q_0\}+L\ds(x,y)Q_0=T_y(x,u',v)\supset F^{[2]}(y).
$$
\par In the same way we show that $S_3\supset F^{[2]}(y)$.
Hence, $A\supset S_1\cap S_2\cap S_3\supset F^{[2]}(y)$.
\par Thus, we have proved that $F^{[2]}(x)+\tau Q_0\supset F^{[2]}(y)$. By interchanging the roles of $x$ and $y$ we obtain also $F^{[2]}(y)+\tau Q_0\supset F^{[2]}(x)$. These two inclusions imply inequality
$$
\dhf(F^{[2]}(x),F^{[2]}(y))\le \tau=\tgm(L)\,\ds(x,y)
$$
proving \rf{H-LTI} with $\tgm(L)=L(3L+1)/(L-1)$.
\smsk
\par We finish the proof of Theorem \reff{LTI-M} in the same fashion as we have finished the proof of Theorem \reff{MAIN-RT} (after the proof of Proposition \reff{HD-G1}). Let $\lambda_1,\lambda_2$ and $\gamma$ be parameters satisfying \rf{PR-LTI}, i.e., $\lambda_1\ge 1$, $\lambda_2\ge 3\lambda_1$ and
$\gamma\ge \lambda_2\,(3\lambda_2+\lambda_1)/
(\lambda_2-\lambda_1)$. We set $\alpha=\lambda_1$, $L=\lambda_2/\lambda_1$ which provides the required inequalities $L\ge 3$ and $\alpha\ge 1$. We also recall that $\ds=\alpha\rho=\lambda_1\,\rho$.
\par In these settings, the mappings $F^{[1]}$ and $F^{[2]}$ are the first and the second order $(\{\lambda_1,\lambda_2\},\rho)$-balanced refinements of $F$ respectively. See Definition \reff{F-IT}.
\par Thanks to Proposition \reff{N-EM}, under the above conditions on  $\alpha=\lambda_1$ and  $L=\lambda_2/\lambda_1$, the set $F^{[2]}(x)\ne\emp$ for every $x\in\Mc$. Thus, property \rf{F2-NEM} holds.
\par In turn, inequality \rf{H-LTI} tells us that
$$
\dhf(F^{[2]}(x),F^{[2]}(y))\le \tgm(L)\ds(x,y)~~~~~\text{with}~~~~~ \tgm(L)=L(3L+1)/(L-1).
$$
Hence,
\be
\dhf(F^{[2]}(x),F^{[2]}(y))&\le& L(3L+1)/(L-1)\,\ds(x,y))
=\frac{\lambda_2}{\lambda_1}\cdot \frac{(3\lambda_2+\lambda_1)}
{(\lambda_2-\lambda_1)}\,(\lambda_1\rho(x,y))\nn\\
&=&
\lambda_2\,\frac{(3\lambda_2+\lambda_1)}
{(\lambda_2-\lambda_1)}
\,\rho(x,y)\le \gamma \,\rho(x,y)
\nn
\ee
proving inequality \rf{HD-RT}.
\par In particular, this inequality holds provided, $\lambda_1=1$, $\lambda_2=3$, and
$$
\gamma=\lambda_2\,(3\lambda_2+\lambda_1)/
(\lambda_2-\lambda_1)=3(3\cdot 3+1)/(3-1)=15.
$$
\par The proof of Theorem \reff{LTI-M} is complete.\bx
\bsk\msk

\indent\par {\bf 6.3 A constructive algorithm for a nearly optimal Lipschitz selection in $\LTI$.}
\addtocontents{toc}{~~~~6.3 A constructive algorithm for a nearly optimal Lipschitz selection in $\LTI$.\hfill \thepage\par\VST}
\msk
\indent

\par The proof of Theorem \reff{LTI-M} provides a certain constructive algorithm for a Lipschitz selection of a set-valued mapping $F$ satisfying the hypothesis of this theorem. Let us briefly describe main steps of this algorithm and give an explicit formula for a nearly optimal Lipschitz selection of $F$.
\par Let $F:\Mc\to\Kc(\RT)$ and let $\lambda$ be a positive constant. We make the following
\begin{assumption}\lbl{FL-AL} For every subset $\Mc'\subset\Mc$ with $\#\Mc'\le 4$, the restriction $F|_{\Mc'}$ of $F$ to $\Mc'$ has a Lipschitz selection $f_{\Mc'}:\Mc'\to \LTI$ with Lipschitz seminorm $\|f_{\Mc'}\|_{\Lip(\Mc',\LTI)}\le \lambda$.
\end{assumption}
\par The following algorithm, given $F$ and $\lambda$ satisfying Assumption \reff{FL-AL}, constructs a Lipschitz mapping $f:\Mc\to\LTI$ with $\|f\|_{\Lip(\Mc,\LTI)}\le 15\lambda$ such that $f(x)\in F(x)$ for each $x\in\Mc$.
\par We constructs $f$ in four steps.
\msk
\par {\bf Step 1.} We construct the $\lambda$-balanced refinement of $F$, i.e., the mapping
$$
F^{[1]}(x)=
\bigcap_{z\in\Mc}\,
\left[F(z)+\lambda\,\rho(x,z)\,Q_0\right],~~~~~~x\in\Mc,
$$
\par {\bf Step 2.} We construct the second order $(\{\lambda,3\lambda\},\rho)$-balanced refinement of $F$, i.e., $3\lambda$-balanced refinement of $F^{[1]}$:
$$
F^{[2]}(x)=\bigcap_{z\in\Mc}\,
\left[F^{[1]}(z)+3\lambda\,\rho(x,z)\,Q_0\right],
~~~~~~x\in\Mc.
$$
\par From the proof of Theorem \reff{LTI-M} we know that
(i) $F^{[1]}(x)\ne\emp$ and $F^{[2]}(x)\ne\emp$ for every $x\in\Mc$, and  (ii) for every $x,y\in\Mc$
\bel{DHL-3}
\dhf(F^{[2]}(x),F^{[2]}(y))\le 15\lambda\,\rho(x,y).
\ee
\par {\bf Step 3.} We construct the rectangular hull of $F^{[2]}$, i.e., the mapping
$$
H(x)=\Hc[F^{[2]}(x)],~~~~~~x\in\Mc.
$$
\par {\bf Step 4.} We define the required mapping $f$ as the center of the rectangle $H(x)$:
\bel{LS-RF}
f(x)=\cntr H(x)=\cntr\Hc[F^{[2]}(x)],~~~~~~x\in\Mc.
\ee
\par Lemma \reff{CH-RH} tells us that $f(x)\in F^{[2]}(x)$ for each $x\in\Mc$. Because $F^{[2]}(x)\subset F^{[1]}(x)\subset F(x)$, $f(x)\in F(x)$ on $\Mc$ proving that $f$ is a selection of $F$. In turn, thanks to \rf{DH-RH1} and \rf{DHL-3},
$$
\dhf(H(x),H(y))=\dhf(\Hc[F^{[2]}(x)],\Hc[F^{[2]}(y)])\le \dhf(F^{[2]}(x),F^{[2]}(y))\le 15\lambda\,\rho(x,y).
$$
Finally, thanks to this inequality and part (ii) of Lemma \reff{SYM},
$$
\|f(x)-f(y)\|=\|\cntr H(x)-\cntr H(y)\|\le
\dhf(H(x),H(y))\le 15\lambda\,\rho(x,y)
$$
proving that $f$ is a Lipschitz selection of $F$ with $\|f\|_{\Lip(\Mc,\LTI)}\le 15\lambda$.
\smsk
\par These observations and representation \rf{RH-TT} enable us to give an explicit formula for the selection $f$.
\par In our settings formula \rf{RH-TT} looks as follows:
set $F^{[2]}(x)$:
\bel{H-F2X}
\HR[F^{[2]}(x)]=
\bigcap\,\,
\HR\left[T_x(u,u',u'')\cap T_x(v,v',v'')\right]
\ee
Here the intersection is taken over all
$u,u',u'',v,v',v''\in\Mc$, and
$$
T_x(u,u',u'')=
\{F(u')+\lambda\,\rho(u',u)Q_0\}\cap \{F(u'')+\lambda\,\rho(u'',u)Q_0\}+3\lambda\,\rho(u,x)Q_0\,.
$$
Recall that, thanks to \rf{LS-RF}, $f(x)=(f_1(x),f_2(x))=
\cntr\Hc[F^{[2]}(x)]$. Let us express the coordinates $f_1(x),f_2(x)$ in the explicit form.
\msk
\par Fix a 6-tuple $\Tc=(u,u',u'',v,v',v'')$ with $u,u',u'',v,v',v''\in\Mc$. Then, thanks to \rf{H-DP} and \rf{B-PR12},
$$
\HR\left[T_x(u,u',u'')\cap T_x(v,v',v'')\right]=
[a_1(x,\Tc),b_1(x,\Tc)]\times[a_2(x,\Tc),b_2(x,\Tc)]]
$$
where for every $j=1,2$
$$
b_j(x,\Tc)=\sup\{y_j:y=(y_1,y_2)\in
T_x(u,u',u'')\cap T_x(v,v',v'')\}
$$
and
$$
a_j(x,\Tc)=\inf\{y_j:y=(y_1,y_2)\in
T_x(u,u',u'')\cap T_x(v,v',v'')\}.
$$
From this and \rf{H-F2X} it follows that
$$
\HR[F^{[2]}(x)]=
[\al_1(x),\beta_1(x)]\times[\al_2(x),\beta_2(x)]
$$
where given $j=1,2$,
\bel{BT-12}
\beta_j(x)=\inf_{\Tc}\,b_j(x,\Tc)=
\inf_{\Tc=(u,u',u'',v,v',v'')}
\sup\{y_j:y=(y_1,y_2)\in
T_x(u,u',u'')\cap T_x(v,v',v'')\}
\ee
and
\bel{ALF-12}
\al_j(x)=\sup_{\Tc}\,a_j(x,\Tc)=
\sup_{\Tc=(u,u',u'',v,v',v'')}
\inf\{y_j:y=(y_1,y_2)\in
T_x(u,u',u'')\cap T_x(v,v',v'')\}.
\ee
\par Thus, $A(x)=(\al_1(x),\al_2(x))$ is "the smallest point", and $B(x)=(\beta_1(x),\beta_2(x))$ is "the largest point" of the rectangle $\HR[F^{[2]}(x)]$. Clearly,
its center, the point $\cntr\HR[F^{[2]}(x)]$, has the coordinates
$$
\cntr\HR[F^{[2]}(x)]=
\left(\frac{\al_1(x)+\beta_1(x)}{2},
\frac{\al_2(x)+\beta_2(x)}{2}\right).
$$
\par Therefore, according to \rf{LS-RF},
\bel{FN-FA}
f_1(x)=\frac{\al_1(x)+\beta_1(x)}{2}~~~~~~\text{and}~~~~~~
f_2(x)=\frac{\al_2(x)+\beta_2(x)}{2}.
\ee
\par This formula and formulae \rf{BT-12}, \rf{ALF-12} and
\rf{FN-FA} provide explicit formulae for a Lipschitz selection of $F$ (with Lipschitz constant at most $15\lambda$) whenever $F$ satisfies Assumption \reff{FL-AL}.
\par We can compare these formulae with corresponding explicit formulae for Lipschitz selection in one dimensional case. See \rf{FP-D1}, \rf{FM-D1} and \rf{FS-D1}. This comparison shows how grows the complexity of the Lipschitz selection problem in transition from the one dimensional to the two dimensional case.

\bsk
\par We complete the section with a refined version of the Finiteness Principle for two dimensional Banach spaces. To  its formulation, given $\gamma>0$ and a convex set $K\subset X$ {\it symmetric} with respect to a point $c\in X$, we let $\gamma\circ K$ denote the dilation of $K$ with respect to $c$ by a factor of $\gamma$.
\begin{theorem} Let $(\Mc,\rho)$ be a pseudometric space, and let $X$ be a two dimensional Banach space. Let $F$ be a set-valued mapping from $\Mc$ into the family $\Kc(X)$ of all non-empty convex compact subsets of $X$. Suppose that for every subset $\Mc'\subset\Mc$ consisting of at most four points, the restriction $F|_{\Mc'}$ of $F$ to $\Mc'$ has a Lipschitz selection $f_{\Mc'}$ with Lipschitz  seminorm $\|f_{\Mc'}\|_{\Lip(\Mc',\X)}\le 1$.
\smallskip
\par Then $F$ has a Lipschitz selection $f$ possessing the following properties: (a) $\|f\|_{\Lip(\Mc,\X)}\le \gamma_1$ where $\gamma_1>0$ is an absolute constant; (b) for every $x\in\Mc$ there exists an ellipse $\Ec_x$ centered at $x$ such that
\begin{align}\lbl{EL-F}
\Ec_x\subset F^{[2]}(x)\subset \gamma_2\circ\Ec_x.
\end{align}
\par Here $\gamma_2>0$ is an absolute constant, and $F^{[2]}$ is the second order $((\lambda_1,\lambda_2),\rho)$-balanced refinement of $F$ with $\lambda_1=4/3$ and $\lambda_2=4$. See \rf{F12-LTI}.
\end{theorem}
\par {\it Proof.} Theorem \reff{MAIN-RT} tells us that the set-valued mapping $F^{[2]}:\Mc\to \Kc(X)$ (with the parameters $\lambda_1=4/3$ and $\lambda_2=4$) is the $\gamma$-core of $F$ with $\gamma=100$. Thus, $F^{[2]}(x)\subset F(x)$ on $\Mc$, and
\begin{align}\lbl{FGR}
\dhf(F^{[2]}(x),F^{[2]}(y))\le \gamma\,\rho(x,y)
~~~~\text{for all}~~~~x,y\in\Mc.
\end{align}
\par It is shown in \cite{S-2004} that there exists a mapping $\ST:\Kc(X)\to X$ (which we call the Steiner-type point map) with the following properties:
\smsk
\par ($\bigstar 1$) $\ST(K)\in K$ for every $K\in\Kc(X)$;
\smsk
\par ($\bigstar 2$) $\|\ST(K)-\ST(K')\|\le \tgm_1\dhf(K,K')$ for every $K, K'\in\Kc(X)$;
\smsk
\par ($\bigstar 3$) for each $K\in\Kc(x)$ there exists an ellipse $\Ec^{(K)}$ centered at $\ST(K)$ such that the following inclusions
$$
\Ec^{(K)}\subset K\subset \tgm_2\circ\Ec^{(K)}
$$
hold. Here $\tgm_1$ and $\tgm_2$ are positive absolute constants.
\smsk
\par Then we set $f(x)=\ST(F^{[2]}(x))$, $x\in\Mc$, and prove that $f$ is a Lipschitz selection of $F$ satisfying \rf{EL-F}. Indeed, thanks to ($\bigstar 1$),
$f(x)\in F^{[2]}(x)\subset F(x)$ on $\Mc$ proving that $f$ is a selection of $F$. Furthermore, by property ($\bigstar 2$) and \rf{FGR}, for every $x,y\in\Mc$ we have
$$
\|f(x)-f(y)\|=\|\ST(F^{[2]}(x))-\ST(F^{[2]}(y))\|\le
\tgm_1\dhf(F^{[2]}(x),F^{[2]}(y))\le\tgm_1\cdot\gamma
\,\rho(x,y)
$$
proving that $\|f\|_{\Lip(\Mc,\X)}\le \gamma_1=\tgm\cdot\gamma=100\,\tgm.$ Finally, the existence of the ellipse $\Ec_x$ satisfying \rf{EL-F} is immediate from the property ($\bigstar 3$).
\par The proof of the theorem is complete. \bx

\par The proof of Theorem \reff{FOR-2} is complete.\bx
\end{document}